\input epsf

\documentclass[12pt,a4paper]{article}
\usepackage[dvips]{graphicx}
\usepackage{latexsym,amsmath,amsfonts,amscd,subfigure}
\usepackage{epsfig}
\usepackage{changebar}
\usepackage{indentfirst}
\usepackage{verbatim}
\usepackage{color}
\usepackage{colordvi}
\usepackage{placeins}
\usepackage{booktabs}
\usepackage{makecell}

\usepackage{algpseudocode}
\usepackage{array}

\begin{document}
\baselineskip=2pc
\begin{center}
{\Large \bf A High-Order Compact Hermite Difference Method for Double-Diffusive Convection\footnote{The research was supported by  NSFC grant 12471390.}}
\end{center}

\centerline{
Jianqing Yang\footnote{School of Mathematical Sciences, Xiamen University, Xiamen, Fujian 361005, PR China E-mail: y\_jianqing@126.com.}
and Jianxian Qiu\footnote{School of Mathematical Sciences and Fujian Provincial Key Laboratory of Mathematical Modeling and High-Performance Scientific Computing, Xiamen University, Xiamen, Fujian 361005, P.R. China. E-mail: jxqiu@xmu.edu.cn. }}

\vspace{.1in}

\baselineskip=1.6 pc

\centerline{\bf Abstract}


In this paper, a class of high-order compact finite difference Hermite scheme is presented for the simulation of  double-diffusive convection.
To maintain linear stability, the convective fluxes are split into positive and negative parts, then the compact Hermite difference methods are used to discretize the positive and  negative fluxes, respectively.  The diffusion fluxes of the governing equations are directly approximated by a high-order finite difference scheme based on the Hermite interpolation. The advantages of the proposed schemes are that the derivative values of the solutions are directly solved by the compact central difference scheme,
and the auxiliary derivative equation is no longer required. The third-order Runge-Kutta method is utilized for the temporal discretization.
Several numerical tests are presented to assess the numerical capability of the newly proposed algorithm. The numerical results are in great agreement with the benchmark solutions and some of the accurate results available in the literature. Subsequently, we apply the algorithm to solve steady and unsteady problems of double-diffusive convection and a preliminary application to the double-diffusive convection for different Raleigh numbers and aspect ratios is carried out.

\vfill {\bf Key Words:} Double-diffusive convection, High-order numerical method, Compact Hermite scheme.
{\bf AMS(MOS) subject classification:} 65M60, 35L65

\pagenumbering{arabic}

\newpage

\baselineskip=1.6pc

\section{\label{sec1:introduction}Introduction}
\setcounter{equation}{0}
\setcounter{figure}{0}
\setcounter{table}{0}
In this paper, we present a class of high-order compact finite difference Hermite scheme for double-diffusive convection (DDC) equation. The governing equations for DDC are time-dependent Navier-Stokes equations with the Boussinesq assumptions. For two dimensional case, these equations expressed in properly non-dimensional form, are given as follows
\begin{eqnarray}\label{eq1:wholeDDC}
&&\nabla^{2}\psi=-\omega, \quad \textbf{u}=\nabla^{\perp}\psi,\label{eq1:1}\\
&&\textbf{U}_{t} +\textbf{F}(\textbf{U})_{x}+\textbf{G}(\textbf{U})_{y}
=\textbf{H}(\textbf{U})_{xx}+\textbf{H}(\textbf{U})_{yy}+\textbf{S},\label{eq1:2}
\end{eqnarray}
where $\nabla^{\perp}=(\partial_{y}, -\partial_{x})$,
$$\textbf{U}=\begin{bmatrix} \omega\\ T\\ C\end{bmatrix}, \quad
\textbf{F}(\textbf{U})=\begin{bmatrix} u\omega\\ uT\\ uC \end{bmatrix}, \quad
\textbf{G}(\textbf{U})=\begin{bmatrix} v\omega\\ vT\\ vC \end{bmatrix}, \quad
\textbf{H}(\textbf{U})=\begin{bmatrix} Pr\omega\\ T\\ \frac{1}{Le}C \end{bmatrix},\quad
\textbf{S}=\begin{bmatrix} PrRa(T_{x}-\lambda C_{x})\\ 0\\ 0 \end{bmatrix}.$$
The dimensionless variables $\psi$, $\omega$, $T$ and $C$ are the stream function, vorticity, temperature and concentration, respectively. $\textbf{u}=(u,v)$ is the dimensionless velocity. $Pr$, $Le$, $Ra$ and $\lambda$ are Prandtl number, Lewis number, Rayleigh number and buoyancy ratio, respectively. They are defined by
\begin{eqnarray}\label{eq2:parameter}
&&Pr=\frac{\nu}{\kappa_{T}},\hspace{5mm}Le=\frac{\kappa_{T}}{\kappa_{C}},\hspace{5mm} Ra=\frac{g\beta_{T}(T_{h}-T_{l})W^{3}}{\kappa_{T}\nu},\hspace{5mm}\lambda=\frac{\beta_{C}(C_{h}-C_{l})}{\beta_{T}(T_{h}-T_{l})},
\end{eqnarray}
where $\beta_{T}, \beta_{C}, \kappa_{T}, \kappa_{C}, \nu$, and $g$ denote the thermal expansion coefficient, concentration expansion coefficient, thermal diffusivity, solutal diffusivity, kinematic viscosity, and acceleration due to gravity, respectively.

Double-diffusive convection occurring in a fluid with density which depends on two scalar components with different diffusion rates \cite{Turner1974,Huppert1976} is a common flow phenomena in nature. Due to its widely scientific and engineering applications in such practical problems as the action of crystal growth, soil contamination and atmospheric pollution, nuclear reactor cooling, it has been one of the major interesting research subjects. Since the complexity of double-diffusive convection, cavity physics modeling\cite{Weaver1991,Chang1993,Nishimura1998}, is extensively used in investigation of the problem. In the early years, some research findings about double-diffusive convection have been achieved in experiments\cite{Weaver1991}, analysis\cite{Garandet1992,Vasseur1995}, and numerical simulations\cite{Chang1993}.
In the recent decades, many numerical results have been reported by means of various numerical methods including the finite-element method,
the finite volume method, the lattice Boltzmann method, etc.
Nishimura et al.\cite{Nishimura1998} reported the effect of the buoyancy ratio in a rectangular cavity using a finite-element method.
Hussaina et al.\cite{Hussaina2017} solved the governing equations of double-diffusive convection by a Galerkin finite element method.

The influence of a broad range of physical parameters on heat and mass transfer in flow fields was systematically investigated.
A staggered-grid finite volume method was developed to numerically study laminar double-diffusive convection in a ventilated enclosure
by Deng et al..\cite{Deng2004}
Wang et al.\cite{Wang2014,Wang2015,Wang2016} investigated the transitions of the flow patterns of the DDC problem under the different physical parameters using a finite volume method.
Verhaeghe et al.\cite{Verhaeghe2007} combined a multi component lattice Boltzmann scheme with a finite difference scheme to simulate the DDC problem. It was found the DDC problem exhibited a rich variety of flow structures. Chen et al.\cite{Chen2017} simulated the double-diffusive convection by a lattice Boltzmann method. They provided an alternative numerical tool for modeling complex heat and mass transfer in fluid-saturated porous media. Later, they proposed new lattice Boltzmann methods for the more complex natural convection, such as the fluids with large Prandtl number\cite{Chen2021,Chen2023}.
We can note that the governing equations of the two-dimensional DDC problem including the vorticity equation, the temperature equation and the concentration equation are essentially convection-diffusion equations. The finite difference method more easily achieves high-order accuracy solutions for convection-diffusion equations. Consequently, it has attracted significant research attention, particularly in developing high-order compact schemes.  Krishnamoorthy et al.\cite{Krishnamoorthy2024} proposed a finite difference method for the double-diffusive natural convection using using a staggered grid system. The focus was on exploring the influence of Soret and Dufour effects and on heat transfer characteristics. Hansda et al.\cite{Hansda2024} chose the fourth-order compact finite difference approximation for solving the governing equations of Double-diffusive. Chuhan et al.\cite{Chuhan2024} numerically investigated double-diffusive mixed convection in a closed rectangular cavity. The finite volume method was used to solve the governing equations, and the temporal terms were discretized using a second-order backward finite difference scheme.  Yang et al.\cite{Yang2024} proposed a meshfree framework based on the reproducing kernel collocation method for nonlinear analysis of the double-diffusive convection in a porous enclosure, and the forward difference method was adopted for temporal discretization.  The local and global entropy generation characteristics of double-diffusive convection was reported in.\cite{Kumar2023}. Karki et al.\cite{Karki2024} considered a cavity containing twin-adiabatic blocks.  The influence of the location of twin-adiabatic blocks on double-diffusive convection and entropy generation was discussed.

Compared to explicit finite difference schemes, compact schemes can achieve comparable accuracy with narrower stencils, thereby reducing the computational burden of boundary discretization. Lele\cite{Lele1992} proposed the symmetric compact schemes and systematically analyzed their resolution characteristics in earlier years.
Peter et al.\cite{Peter1998} proposed a three-point combined compact difference scheme.
It becomes apparent that the nonphysical oscillation would occur when the symmetric compact central schemes are applied to discretize the nonlinear convection term. The compact upwind schemes with inherent dissipation can dissipate the nonphysical oscillation. Therefore, the upwind schemes are preferable in the direct numerical simulation of strong flow problems.
Liang et al.\cite{Liang2009} used the fifth-order upwind compact scheme for approximating convection terms and the fourth-order symmetrical compact scheme for viscous terms to simulate the double-diffusive convection in a rectangular cavity.  Their study initially focused on the transition of flow structures with increasing buoyancy ratio from 0 to 1. Subsequently, they employed the same methodology to systematically examine the influence of buoyancy ratio ($(1\leq N \leq2)$) on heat and mass transfer and flow characteristics \cite{Liang2019}.
Qin et al. \cite{Qin2014} developed a weighted-average-based combined upwind compact scheme with fourth-order and fifth-order accuracy to systematically analyze the influence of Prandtl and Lewis numbers on flow structures in double-diffusive convection systems.
Yu et al. \cite{Yu2017} used the same finite different method to report some results about the double-diffusive convection under a uniform horizontal magnetic field and heat source.
Yan et al.\cite{Yan2022} developed a high-order upwind compact scheme for simulating double-diffusive convection. Their numerical results are excellent agreement with the data reported in existing literature.
While considerable advances have been achieved in simulating DDC problems, the development of effective physics models and advanced numerical methods remains an area of significant research interest.
So far, most scholars mainly propose new numerical methods for Non-conservative form of the governing equations of the two-dimensional double-diffusive convection. While the research on the conservative form of the governing equations is relatively weak.

In this paper, we propose a class of high-order compact Hermite difference scheme for solving the conservative form of the two-dimensional double-diffusive convection problem. The scheme combines a Hermite interpolation polynomial approximation and a central compact difference scheme. It is well known, the Hermite interpolation requires not only the function values at the interpolation nodes, but also the values of their derivatives. Most scholars take the partial derivative for variable on both sides of the governing equation and obtain a new equation about the partial derivative of variable. They obtain the derivative values of the variable by simultaneously solving the original equation and its derivative equation\cite{Qiu2004,Qiu2005,ZhaoZ2020}.
This treatment will incur additional computational costs, for there is the auxiliary derivative equation is required to evaluate the derivative values of the variable. In this paper, similar to \cite{MaWu2018} in which a new type of Hermite weighted essentially non-oscillatory schemes based on compact difference schemes was proposed for solving both one and two dimensional hyperbolic conservation laws, instead of solving the auxiliary derivative equation, the derivative values of the solutions are directly solved by the compact central difference scheme, the efficiency of the new algorithm is great improved.

The paper is organized as follows. In Section 2  we would introduce the Hermite interpolation and compact difference scheme for double-diffusive convection equation. Numerical experiments are performed for the validation of accuracy and efficiency of the newly proposed algorithm in Section 3. Concluding remarks are included in Section 4.

\section{\label{sec2:Algorithms}Numerical methods}
\setcounter{equation}{0}
\setcounter{figure}{0}
\setcounter{table}{0}
In this section, we describe the compact Hermite difference scheme for solving the governing equations (\ref{eq1:1})-(\ref{eq1:2}) of  two-dimensional laminar double-diffusive convection problem in detail.

\subsection{\label{sec2.1:Problem} Description of compact Hermite difference scheme}
For the two-dimensional case, the computation domain is divided by uniform meshes $I_{i,j}=[x_{i-\frac{1}{2}},x_{i+\frac{1}{2}}]\times[y_{j-\frac{1}{2}},y_{j+\frac{1}{2}}],\ i=1,2,...,N_{x},\ j=1,2,...,N_{y} $, and $(x_{i},y_{j})$ is the center of $I_{i,j}$. $\Delta x=x_{i+1/2}-x_{i-1/2}$, $\Delta y=y_{j+1/}-y_{j-1/2}$ are spatial step sizes in $x$ and $y$ direction, respectively. The semi-discretization high order finite difference scheme for Eq.\eqref{eq1:2} can be written in conservation form as following
\begin{equation}\label{eq3:semi-discrete}
\begin{aligned}
\frac{d\textbf{U}_{i,j}}{dt} = & -\frac{1}{\Delta x}
\left(\widehat{\textbf{F}}_{i+\frac{1}{2},j}-\widehat{\textbf{F}}_{i-\frac{1}{2},j}\right)
-\frac{1}{\Delta y}\left(\widehat{\textbf{G}}_{i,j+\frac{1}{2}}-\widehat{\textbf{G}}_{i,j-\frac{1}{2}}\right)\\
& +\textbf{\emph{M}}_{i,j}(\textbf{U})+\textbf{\emph{R}}_{i,j}(\textbf{U})+ \textbf{S}_{i,j},
\end{aligned}
\end{equation}
where $\widehat{\textbf{F}}_{i\pm\frac{1}{2},j}$ and $\widehat{\textbf{G}}_{i,j\pm\frac{1}{2}}$ are the numerical fluxes for convection terms, and
$\textbf{\emph{M}}_{i,j}(\textbf{U})$ and $\textbf{\emph{R}}_{i,j}(\textbf{U})$ are approximations of the diffusion terms $\textbf{H}(\textbf{U})_{xx}$ and $\textbf{H}(\textbf{U})_{yy}$, respectively, at $(x_{i},y_{j})$.  $\textbf{S}_{i,j}$ is the approximation value of $\textbf{S}(x_{i},y_{j})$. Those terms would be reconstructed by Hermite finite difference method combine with compact central finite difference method.

Now we would like to describe the procedure to reconstruct the $\widehat{\textbf{F}}_{i\pm\frac{1}{2},j}$,  $\textbf{\emph{M}}_{i,j}(\textbf{U})$,  and $\textbf{S}_{i,j}$ based on the compact Hermite difference. The procedure to reconstruct the $\widehat{\textbf{G}}_{i,j\pm\frac{1}{2}}$ and $\textbf{\emph{R}}_{i,j}(\textbf{U})$ is similar to those for $\widehat{\textbf{F}}_{i\pm\frac{1}{2},j}$,and $\textbf{\emph{M}}_{i,j}(\textbf{U})$, respectively.

{\bf Step 1. The procedure to reconstruct $\widehat{\textbf{F}}_{i\pm\frac{1}{2},j}$.}

For the stability of the present algorithm, the convection flux $\textbf{F}(\textbf{U})$ is split into two parts
\begin{equation}\label{eq5:split}
\textbf{F}(\textbf{U})=\textbf{F}^{+}(\textbf{U})+\textbf{F}^{-}(\textbf{U}),
\end{equation}
where $\textbf{F}^{\pm}(\textbf{U})$ are computed by the global Lax-Friedrichs flux splitting which is
$$\textbf{F}^{\pm}(\textbf{U})=\frac{1}{2}(\textbf{F}(\textbf{U})\pm\alpha\textbf{U}), $$
and  $\alpha=\max_{U}\rho(\textbf{J})$, $\rho(\textbf{J})$ is spectral radius of the Jacobian matrix $\textbf{J}=\frac{d\textbf{F}}{d\textbf{U}}$. The numerical flux is also split into two parts
 \begin{equation}
 \label{fluxes}\widehat{\textbf{F}}_{i+\frac{1}{2},j}=\widehat{\textbf{F}}^{+}_{i+\frac{1}{2},j}
+\widehat{\textbf{F}}^{-}_{i+\frac{1}{2},j}.
 \end{equation}
 Following \cite{so2,js}, we define the function $\Phi(x,y)$ such that
 \begin{equation}
 \label{eq4:definition}
\textbf{F}^+\left(\textbf{U}(x,y)\right)=\frac{1}{\Delta x}\int^{x+\Delta x/2}_{x-\Delta x/2}\Phi(m,y)dm,
\end{equation}
 then we have
 $$\frac{\partial \textbf{F}^+\left(\textbf{U}(x_i,y_j)\right)}{\partial x}=\frac 1{\Delta x}(\Phi(x_{i+1/2},y_j)-\Phi(x_{i-1/2},y_j)).$$

 If $\widehat{\textbf{F}}^+_{i+\frac{1}{2},j}$ is the high order
 approximation to $\Phi(x_{i+1/2},y_j)$ and Lipschitz continuous concerning their arguments, then $\frac1{\Delta x}(\widehat{\textbf{F}}^+_{i+\frac{1}{2},j}-\widehat{\textbf{F}}^+_{i+\frac{1}{2},j})$ is the high order approximation of the $\frac{\partial \textbf{F}^+\left(\textbf{U}(x_i,y_j)\right)}{\partial x}.$

 Now we give the detail of procedure to reconstruct $\widehat{\textbf{F}}^+_{i+\frac{1}{2},j}$.  For convenience,  let $\textbf{F}_{i,j}$  and $(\textbf{F}_{x})_{i,j}$  represent the values of $\textbf{F}(\textbf{U})$ and $\textbf{F}(\textbf{U})_{x}$, at the point $(x_{i},y_{j})$, respectively. The derivative values of $(\textbf{F}_{x})_{i,j}$ would be computed by the central compact scheme which is described in Step 3.
 From (\ref{eq4:definition}), we have
 $$\overline{\Phi}_{i,j}=\frac{1}{\Delta x}\int^{x_i+\Delta x/2}_{x_i-\Delta x/2}\Phi(m,y_j)dm={\textbf{F}}^+_{i,j},$$

 $$\overline{(\Phi_x)}_{i,j}=\frac{1}{\Delta x}\int^{x_i+\Delta x/2}_{x_i-\Delta x/2}\frac{\partial \Phi(m,y_j)}{\partial m}dm=(\textbf{F}_{x})^+_{i,j}.$$

There is a class of fifth-degree Hermite interpolation polynomial $p(x)$ on the stencil $S_0=\{x_{i-1}, x_{i}, x_{i+1}\}$ such that
\begin{equation}\label{eq8:condition}
\begin{aligned}
& \frac{1}{\Delta x}\int_{I_{i+l}}p(x)dx=\textbf{F}^{+}_{i+l,j},
           \quad \frac{1}{\Delta x}\int_{I_{i+l}}p(x)_{x}dx=(\textbf{F}_{x})^{+}_{i+l,j}, l=-1, 0, 1,\\
\end{aligned}
\end{equation}
where $p(x)$ is the  Hermite interpolation polynomial of $\Phi(x,y_j)$ in \eqref{eq4:definition}, and we take $\widehat{\textbf{F}}^+_{i+\frac{1}{2},j}=p(x_{i+\frac{1}{2}})$, which is the sixth order approximation to $\Phi(x_{i+\frac{1}{2}},y_j)$
\begin{equation}\label{eq9:F-hat}
\begin{aligned}
\widehat{\textbf{F}}^+_{i+\frac{1}{2},j}=\frac{11}{60}\textbf{F}^{+}_{i-1,j}+\frac{19}{30}\textbf{F}^{+}_{i,j}+\frac{11}{60}\textbf{F}^{+}_{i+1,j}
                       +\frac{\Delta x}{20}\left((\textbf{F}_{x})^{+}_{i-1,j}+10(\textbf{F}_{x})^{+}_{i,j}-(\textbf{F}_{x})^{+}_{i+1,j}\right).
\end{aligned}
\end{equation}
The reconstruction for $\widehat{\textbf{F}}^{-}_{i+\frac{1}{2},j}$ is mirror symmetric of that for $\widehat{\textbf{F}}^{+}_{i+\frac{1}{2},j}$ with respect to $x_{i+\frac{1}{2}}$. The reconstruction for $\widehat{\textbf{G}}^{\pm}_{i,j+\frac{1}{2}}$ is the same as that for $\widehat{\textbf{F}}^{\pm}_{i+\frac{1}{2},j}$ with $i$ and  $j$ interchanged.

{\bf Step 2. The procedure to reconstruct $\textbf{\emph{M}}_{i,j}(\textbf{U}).$ }

For $\textbf{H}(\textbf{U})_{xx}$, we directly present a sixth-order approximation $\textbf{\emph{M}}_{i,j}(\textbf{U})$ based on Hermite interpolation on the stencil $S_{1}=\{x_{i-2}, x_{i-1}, x_{i+1}, x_{i+2}\}$. Again, for convenience, we denote $\textbf{H}_{i,j}$, $(\textbf{H}_{x})_{i,j}$ and $(\textbf{H}_{xx})_{i,j}$  for  the values of $\textbf{H}(\textbf{U})$, $\textbf{H}(\textbf{U})_{x}$, $\textbf{H}(\textbf{U})_{xx}$ at the point $(x_{i},y_{j})$, respectively.
A class of seventh-degree polynomial $q(x,y)$ can be obtained on this
stencil such that
\begin{equation}\label{eq10:condition}
\begin{aligned}
& q(x_{i+l},y_{j})=\textbf{H}_{i+l,j},\ l=-2,-1,0, 1,2,
           \quad q(x_{i+l},y_{j})_{x}=(\textbf{H}_{x})_{i+l,j}, \ l=-1, 0, 1.\\
\end{aligned}
\end{equation}
Then $q(x,y)$ is the eighth order approximation to $\textbf{H}(\textbf{U})$, and $q(x,y)_{xx}$ is the sixth order approximation to $\textbf{H}(\textbf{U})_{xx}$. We take $\emph{\textbf{M}}_{i,j}=q(x_{i},y_{j})_{xx}$, which is the sixth order approximation to $(\textbf{H}_{xx})_{i,j}$. By performing algebraic operations, we have
\begin{equation}\label{eq11:Hxx}
\begin{aligned}
\emph{\textbf{M}}_{i,j}=q(x_{i},y_{j})_{xx}&=\frac{1}{36\Delta x^{2}}(
                    \textbf{H}_{i-2,j}+80\textbf{H}_{i-1,j}-162\textbf{H}_{i,j}+80\textbf{H}_{i+1,j}\\
                   &+\textbf{H}_{i+2,j}+24\Delta x(\textbf{H}_{x})_{i-1,j}-24\Delta x(\textbf{H}_{x})_{i+1,j}).
\end{aligned}
\end{equation}

Similar to the procedure for reconstruction of $\textbf{\emph{M}}_{i,j}(\textbf{U})$, we have
\begin{equation}\label{eq12:Hyy}
\begin{aligned}
\emph{\textbf{R}}_{i,j}&=\frac{1}{36\Delta y^{2}}(
                    \textbf{H}_{i,j-2}+80\textbf{H}_{i,j-1}-162\textbf{H}_{i,j}+80\textbf{H}_{i,j+1}\\
                   &+\textbf{H}_{i,j+2}+24\Delta y(\textbf{H}_{y})_{i,j-1}-24\Delta y(\textbf{H}_{y})_{i,j+1}).
\end{aligned}
\end{equation}

{\bf Step 3. The procedure to reconstruct the derivatives.}

Following the idea of Hermite interpolation, the schemes of \eqref{eq9:F-hat}, \eqref{eq11:Hxx} and \eqref{eq12:Hyy} require the values of $(\textbf{F}_{x})_{i,j}$, $(\textbf{H}_{x})_{i,j}$ and $(\textbf{H}_{y})_{i,j}$ and its neighboring points. We choose the central compact difference scheme to solve the derivative values instead of solving the derivative equation.  Let $\Upsilon$ represents $\textbf{F}$ or $\textbf{H}$.
The central compact difference scheme is based on the symmetric Taylor series expansions, one can obtain central difference approximations of the first derivative. The standard fourth-order central compact difference (CCD4) scheme\cite{Lele1992} can be expressed as
\begin{equation}\label{eq13:CCD4}
\begin{aligned}
& \frac{1}{6}(\Upsilon_{x})_{i-1,j}+\frac{2}{3}(\Upsilon_{x})_{i,j}+\frac{1}{6}(\Upsilon_{x})_{i+1,j}=
\frac{1}{2\Delta x}(\Upsilon_{i+1,j}-\Upsilon_{i-1,j}),\\
\end{aligned}
\end{equation}
and the sixth-order central compact difference (CCD6) scheme\cite{Lele1992} is
\begin{equation}\label{eq14:CCD6}
\begin{aligned}
& \frac{1}{3}(\Upsilon_{x})_{i-1,j}+(\Upsilon_{x})_{i,j}+\frac{1}{3}(\Upsilon_{x})_{i+1,j}=
\frac{7}{9\Delta x}(\Upsilon_{i+1,j}-\Upsilon_{i-1,j})+\frac{1}{36\Delta x}(\Upsilon_{i+2,j}-\Upsilon_{i-2,j}).\\
\end{aligned}
\end{equation}
The approximation for $\Upsilon_{y}$ is the same as that for $\Upsilon_{x}$ with $i$ and $j$ interchanged.
The final order of the compact Hermite difference scheme is determined depending on the chosen CCD scheme for the derivatives of $\textbf{F}_{x}$, $\textbf{F}_{y}$, $\textbf{H}_{x}$ and $\textbf{H}_{y}$. So we denote the compact Hermite difference scheme as CHD4 when the CCD4 scheme is applied in the computing. When the CCD6 scheme is used to compute, we denote the compact Hermite difference scheme as CHD6.

Additionally, we further observe that the components of $\textbf{S}$ contain first partial derivatives $T_{x}$ and $C_{x}$. The terms $T_{x}$ and $C_{x}$ can be computed using the CCD4 or CCD6 scheme. We also use the CCD4 or CCD6 scheme to compute the velocity $\textbf{u}$,  that they are
\begin{equation}\label{eq15:velocity-uv}
\begin{aligned}
& \frac{1}{6}u_{i,j-1}+\frac{2}{3}u_{i,j}+\frac{1}{6}u_{i,j+1}=
\frac{1}{2\Delta y}(\psi_{i,j+1}-\psi_{i,j-1}),\\
& \frac{1}{6}v_{i-1,j}+\frac{2}{3}v_{i,j}+\frac{1}{6}v_{i+1,j}=
-\frac{1}{2\Delta x}(\psi_{i+1,j}-\psi_{i-1,j}),\\
\end{aligned}
\end{equation}
and
\begin{equation}\label{eq16:velocity-uv}
\begin{aligned}
& \frac{1}{3}u_{i,j-1}+u_{i,j}+\frac{1}{3}u_{i,j+1}=
\frac{7}{9\Delta y}(\psi_{i,j+1}-\psi_{i,j-1})+\frac{1}{36\Delta y}(\psi_{i,j+2}-\psi_{i,j-2}),\\
&\frac{1}{3}v_{i-1,j}+v_{i,j}+\frac{1}{3}v_{i+1,j}=
-\frac{7}{9\Delta x}(\psi_{i+1,j}-\psi_{i-1,j})-\frac{1}{36\Delta x}(\psi_{i+2,j}-\psi_{i-2,j}).\\
\end{aligned}
\end{equation}

The equation of stream function $\psi$ in Eq.\eqref{eq1:1} is a Poisson-type equation.
A fourth-order compact symmetric scheme\cite{Tian2007} is used to solve the equation of stream function.
It is
\begin{equation}\label{eq17:Possion}
\begin{aligned}
&2\left(\frac{5}{\Delta x^{2}}-\frac{1}{\Delta y^{2}}\right)
(\psi_{i+1,j}-\psi_{i-1,j})
+2\left(\frac{5}{\Delta y^{2}}-\frac{1}{\Delta x^{2}}\right)
(\psi_{i,j+1}-\psi_{i,j-1})\\
&+ \left(\frac{1}{\Delta x^{2}}+\frac{1}{\Delta y^{2}}\right)
(\psi_{i+1,j+1}+\psi_{i+1,j-1}+\psi_{i-1,j+1}+\psi_{i-1,j-1}-20\psi_{i,j})\\
&=-(8\omega_{i,j}+\omega_{i+1,j}+\omega_{i-1,j}+\omega_{i,j+1}+\omega_{i,j-1}).\\
\end{aligned}
\end{equation}
The successive Over-Relaxation (SOR) method is employed to solve Eq.\eqref{eq17:Possion}, with a geometric multigrid for saving in computing time.

\subsection{\label{sec2.2:Boundaryscheme}Numerical scheme for the boundary condition}
To ensure the closure of the compact Hermite difference scheme, some special methods need to be implemented near the boundary for non-periodic problems.
As mentioned earlier, $\Upsilon$ can be $\textbf{F}$ or $\textbf{H}$.
The values of $\Upsilon$ near the left boundary are obtained by a quintic polynomial $l(x,y)$.  $l(x,y)$ satisfies
\begin{eqnarray}\label{eq18:condition}
&&l(x_{i+m},y_{j})=\Upsilon_{i+m,j}, \ l(x_{i+m}y_{j})_{x}=(\Upsilon_{x})_{i+m,j},\ m=1,2,3.
\end{eqnarray}
The values of $\Upsilon$ and $\Upsilon_{x}$ near the left boundary are obtained by $\Upsilon_{i,j}=l(x_{i},y_{j}), (\Upsilon_{x})_{i,j}=l(x_{i},y_{j})_{x}$. One can get
\begin{equation}\label{eq19:LeftBoundry}
\begin{aligned}
&\Upsilon_{i,j}=-18\Upsilon_{i+1,j}+9\Upsilon_{i+2,j}+10\Upsilon_{i+3,j}-
                \Delta x \left(9(\Upsilon_{x})_{i+1,j}+18(\Upsilon_{x})_{i+2,j}+3(\Upsilon_{x})_{i+3,j} \right),\\
&(\Upsilon_{x})_{i,j}=\frac{1}{\Delta x}\left(57\Upsilon_{i+1,j}-24\Upsilon_{i+2,j}-33\Upsilon_{i+3,j}\right)
                     + 24(\Upsilon_{x})_{i+1,j}+57(\Upsilon_{x})_{i+2,j}+10(\Upsilon_{x})_{i+3,j}.
\end{aligned}
\end{equation}
We can construct the values of $\Upsilon$ and $\Upsilon_{x}$ near the right boundary using the same idea. $r(x,y)$ is a quintic polynomial satisfying
\begin{eqnarray}\label{eq20:condition}
&&r(x_{i-m},y_{j})=\Upsilon_{i-m,j}, \ r(x_{i-m},y_{j})=(\Psi_{x})_{i-m,j},\ m=1,2,3.
\end{eqnarray}
The values of $\Upsilon_{i,j}$, $(\Upsilon_{x})_{i,j}$ are $r(x_{i},y_{j})$, $r(x_{i},y_{j})_{x}$, respectively. They are explicitly expressed as
\begin{equation}\label{eq21:RightBoundry}
\begin{aligned}
&\Upsilon_{i,j}=-18\Upsilon_{i-1,j}+9\Upsilon_{i-2,j}+10\Upsilon_{i-3,j}-
             \Delta x(9(\Upsilon_{x})_{i-1,j}+18(\Upsilon_{x})_{i-2,j}+3(\Upsilon_{x})_{i-3,j}),\\
&(\Upsilon_{x})_{i,j}=\frac{1}{\Delta x}(57\Upsilon_{i-1,j}-24\Upsilon_{i-2,j}-33\Upsilon_{i-3,j})+
                24(\Upsilon_{x})_{i-1,j}+57(\Upsilon_{x})_{i-2,j}+10(\Upsilon_{x})_{i-3,j}.\\
\end{aligned}
\end{equation}

For the CCD4 scheme, we use the fourth-order boundary schemes. The left and right boundary schemes\cite{Yu2012} are
\begin{equation}\label{eq22:LeftBoundry}
\begin{aligned}
& (\Upsilon_{x})_{i,j}+3(\Upsilon_{x})_{i+1,j}=\frac{1}{6\Delta x}(-17\Upsilon_{i,j}+9\Upsilon_{i+1,j}+ 9\Upsilon_{i+2,j}-\Upsilon_{i+3,j}),\\
& (\Upsilon_{x})_{i,j}+3(\Upsilon_{x})_{i-1,j}=-\frac{1}{6\Delta x}(-17\Upsilon_{i,j}+9\Upsilon_{i-1,j}+ 9\Upsilon_{i-2,j}-\Upsilon_{i-3,j}).
\end{aligned}
\end{equation}
The sixth-order boundary schemes\cite{Ma2019} are
\begin{equation}\label{eq23:LeftBoundry}
\begin{aligned}
&(\Upsilon_{x})_{i,j}+5(\Upsilon_{x})_{i+1,j}=\frac{1}{\Delta x}(-\frac{197}{60}\Upsilon_{i,j}-\frac{5}{12}\Upsilon_{i+1,j}+
                   5\Upsilon_{i+2,j}-\frac{5}{3}\Upsilon_{i+3,j}+\frac{5}{12}\Upsilon_{i+4,j}-\frac{1}{20}\Upsilon_{i+5,j}),\\
&(\Upsilon_{x})_{i,j}+5(\Upsilon_{x})_{i-1,j}=-\frac{1}{\Delta x}(-\frac{197}{60}\Upsilon_{i,j}-\frac{5}{12}\Upsilon_{i-1,j}+
                   5\Upsilon_{i-2,j}-\frac{5}{3}\Upsilon_{i-3,j}+\frac{5}{12}\Upsilon_{i-4,j}-\frac{1}{20}\Upsilon_{i-5,j}),
\end{aligned}
\end{equation}
for the CCD6 scheme. The same way can be used to compute the values near the boundary in the $y$ direction.

\subsection{\label{sec2.3:Temporal} Temporal discretization}
The semi-discrete scheme is equivalent to the ordinary differential equation system after spatial discretization with above method. We see
Eq. \eqref{eq3:semi-discrete} as
\begin{equation}\label{eq27:ordinaryEq}
\begin{aligned}
&\frac{d\textbf{U}}{dt}=R(\textbf{U}),
\end{aligned}
\end{equation}
where $R(\textbf{U})$ is discrete operators of spatial derivatives. For the discretization in time, the explicit third-order Runge-Kutta scheme is used. It is given as
\begin{eqnarray*}
  &&\textbf{U}^{(1)}=\textbf{U}^{n}+\Delta tR(\textbf{U}^{n}),\\
  &&\textbf{U}^{(2)}=\frac{3}{4}\textbf{U}^{n}+\frac{1}{4}(\textbf{U}^{(1)}+\Delta tR(\textbf{U}^{(1)})),\\
  &&\textbf{U}^{n+1}=\frac{1}{3}\textbf{U}^{n}+\frac{2}{3}(\textbf{U}^{(2)}+\Delta tR(\textbf{U}^{(2)})).
\end{eqnarray*}
$R(\textbf{U}^{n})$ is the discrete results of the right source of Eq. \eqref{eq27:ordinaryEq} at the time, $n$.

\subsection{\label{sec2.4:algorithm} The algorithm for the double-diffusive system}
Now we outline the new algorithm in brief. Assuming $\psi^{n}$, $\textbf{u}^{n}$ and $\textbf{U}^{n}$ are known. As previously introduced, $\textbf{U}^{n}=(\omega^{n}, T^{n}, C^{n})^T$. The solutions for $\omega^{n+1}$, $T^{n+1}$, $C^{n+1}$, $\psi^{n+1}$ and $\textbf{u}^{n+1}$ are obtained through the following iterative procedure.
\begin{table}[htbp]
\centering
\setlength{\tabcolsep}{1.5pt}
\begin{tabular}{
    >{\raggedright\arraybackslash}p{1.5cm}
    >{\raggedright\arraybackslash}p{14cm}
  }
\hline
\multicolumn{2}{l}{\textbf{Algorithm for the double-diffusive system} } \vspace{2mm}\\
\hline
Step 1: &Solutions of $\textbf{U}^{(n+1,k+1)}$ ($(n,k)$ denotes the $k$th iteration step at the time $n$).\\
        &Initialize $\psi^{(n+1,0)}=\psi^{n}$, $\textbf{u}^{(n+1,0)}=\textbf{u}^{n}$ and $\textbf{U}^{(n+1,0)}=\textbf{U}^{n}$.\\
        &(i) Obtain the boundary values of $\omega^{(n+1,k+1)}$, $T^{(n+1,k+1)}$, $C^{(n+1,k+1)}$.\\
        &Obtain $\omega^{(n+1,k+1)}$ on the boundaries using the scheme \eqref{eq26:vorticity}.
        $T^{(n+1,k+1)}$ and $C^{(n+1,k+1)}$ are obtained on the horizontal walls using the schemes \eqref{eq24:temperature} and \eqref{eq25:concentration}, respectively. \\
        &(ii) Obtain $(T_{x})_{i,j}^{(n+1,k)}$ and $(C_{x})_{i,j}^{(n+1,k)}$ in the term $PrRa(T_{x}-\lambda C_{x})$ using the schemes
          \eqref{eq13:CCD4} and \eqref{eq22:LeftBoundry} (\eqref{eq14:CCD6} and \eqref{eq23:LeftBoundry}).\\
        &(iii) Obtain $R(\textbf{U}^{(n+1,k)})$  of the equation \eqref{eq27:ordinaryEq} in the inner region.\\
        &The steps (i) to (iii) are implemented at each Runge-Kutta stage.
          After finishing the third-order Runge-Kutta steps, one can get $\textbf{U}^{(n+1,k+1)}$.\\
\hline
Step 2: & Solution of $\psi^{(n+1,k+1)}$.\\
        & obtain the stream-function $\psi^{(n+1,k+1)}$ using the scheme\eqref{eq17:Possion}.\\
\hline
Step 3:& Solution of $\textbf{u}^{(n+1,k+1)}$.\\
       & obtian the velocity $\textbf{u}^{(n+1,k+1)}$ using the scheme \eqref{eq15:velocity-uv} or \eqref{eq16:velocity-uv}.\\
       \hline
Step 4:& Check the convergence of $\textbf{U}$.\\
       & Set {$Err=\max\left\{\left|\frac{\textbf{U}^{(n+1,k+1)}-\textbf{U}^{(n+1,k)}}{\textbf{U}^{(n+1,k+1)}}\right|\right\}$.}\\
       & \qquad If $Err>\epsilon$, where $\epsilon$ is the tolerance limit. Set $k=k+1$, then go back to Step1--Step3.\\
       & \qquad If $Err<\epsilon$, the convergence solutions of $\textbf{U}^{n+1}=\textbf{U}^{(n+1,k+1)}$, $\psi^{n+1}=\psi^{(n+1,k+1)}$,\\
       & \qquad and $\textbf{u}^{n+1}=\textbf{u}^{(n+1,k+1)}$ are obtained. So far, the inner iteration process is end.\\
       \hline
Step 5:& Determine the stability of the solution.\\
       &Let  {$\Delta E_{\max}=\max[(u^{n+1}_{i,j}-u^{n}_{i,j})^{2}+(v^{n+1}_{i,j}-v^{n}_{i,j})^{2}]^{\frac{1}{2}}$.} \\
       &The solution is steady when $\Delta E_{\max}\leq10^{-10}$.\\ \hline
\end{tabular}
\end{table}

\section{\label{sec3:Numerical experiments}Numerical Experiments}
\setcounter{equation}{0}
\setcounter{figure}{0}
\setcounter{table}{0}
To comprehensively evaluate the spatial performance of the proposed algorithm, we conduct a series of systematic numerical experiments. The reliability and accuracy of our numerical method are rigorously validated through detailed comparisons with benchmark results reported in the literature.
\subsection{\label{sec3.1:Accuracy test} Accuracy tests}
To verify the accuracy and reliability of the proposed algorithm, we apply it to solve the convection-diffusion equation and nonlinear Burgers' equation with analytical solutions. The convergence order is defined by
\begin{equation}\label{eq28:rate}
\begin{aligned}
&Order=\frac{\log (L_{2}(h_{1})/L_{2}(h_{2})}{\log(h_{1}/h_{2})},
\end{aligned}
\end{equation}
where $L_{2}=\sqrt{\Delta x \sum_{i}(\varphi_{i}^{num}-\varphi_{i}^{exact})^{2}}$ for the one dimensional case, $h_{1}$ and $h_{2}$ are two distinct mesh sizes.
We also define $L_{\infty}=\max_{i}(\varphi^{num}_{i}-\varphi_{i}^{exact})$. $\varphi^{num}$ and $\varphi^{exact}$ are the computed solution and the exact solution, respectively. For the two dimensional case, the definitions of $L_{2}$, $L_{\infty}$  are similar to the definitions of $L_{2}$, $L_{\infty}$ of the one dimensional case.

\subsubsection{\label{sec3.1.1:1D}One dimensional convection-diffusion equation}
We solve the one dimensional convection-diffusion equation
\begin{equation}\label{eq28:CD-1Dcase}
\frac{\partial u}{\partial t}+\frac{\partial u}{\partial x}=\frac{\partial^{2} u}{\partial x^{2}},\ x\in[0,2\pi],\ t\in(0,T].
\end{equation}
The exact solution is $u(x,t)=e^{-t}sin(x-t)$. The initial condition is directly obtained from its exact solution and the period boundary condition is used in the present problem.
\begin{table}[htbp]
\centering
\caption{One dimensional convection-diffusion equation: $L_{{2}}$ and $L_{{\infty}}$ errors and convergence orders at $T=1$, $\Delta t=\Delta x^{2}$. }\label{table1:CD-1Dcase}
\begin{tabular}{lllllllll}
\hline
& \multicolumn{4}{c}{CHD4}& \multicolumn{4}{c}{CHD6}\vspace{1mm}\\
  \cline{2-9}
$N_{x}$ & $L_{{2}}$-error&Rate&$L_{{\infty}}$-error&Rate
  &$L_{{2}}$-error&Rate&$L_{{\infty}}$-error&Rate\\
\hline
$20$  &5.326(-5)&-   &8.648(-5)& -  &5.895(-6) &-   &1.250(-5)  & - \\
$40$  &2.062(-6)&4.69&3.082(-6)&4.81&4.729(-8) &6.96&8.892(-8)  &7.13\\
$60$  &3.336(-7)&4.49&4.994(-7)&4.49&2.612(-9) &7.14&4.729(-9)  &7.24\\
$80$  &9.473(-8)&4.38&1.434(-7)&4.34&3.298(-10)&7.19&5.811(-10)&7.29\\
$100$ &3.628(-8)&4.30&5.544(-8)&4.26&6.629(-11)&7.19&1.136(-10)&7.32\\
$120$ &1.671(-8)&4.25&2.574(-8)&4.21&1.894(-11)&6.87&3.126(-11)&7.07\\
\hline
\end{tabular}
\end{table}
Table \ref{table1:CD-1Dcase} shows the $L_{2}$, $L_{{\infty}}$ errors and convergence orders for the CHD4 and  CHD5 schemes.
It can be seen that the results computed by the CHD4 and CHD5 schemes can achieve at least fourth-order and sixth-order accuracy in space, respectively, which are in good agreement with the theoretical predictions.

\subsubsection{\label{sec3.1.2:1D}One dimensional Burgers' equation}
Solve the one dimensional Burgers' equation\cite{Gao2014,Xie2008,Jia2019,Sheng2023}
\begin{equation}\label{eq29:Burgers}
\frac{\partial u}{\partial t}+\frac{\partial (\frac{1}{2}u^{2})}{\partial x}=\varepsilon\frac{\partial^{2} u}{\partial x^{2}},\ x\in[a,b],\ t\in(t_{0},T].
\end{equation}
The analytical solution to the equation is
\begin{equation}\label{eq30:BurgersSolution1}
u(x,t)=\frac{2\pi\varepsilon e^{-\pi^{2}\varepsilon t}sin(\pi x)}
{\gamma+e^{-\pi^{2}\varepsilon t}cos(\pi x)},
\end{equation}
where $x\in[0,1],\ t\in(0,1]$.
\begin{table}[htbp]
\centering
\caption{\label{table2:BurgersSolution1}%
The exact solution \eqref{eq30:BurgersSolution1} of the one dimensional Burgers' equation: $L_{{2}}$ and $L_{{\infty}}$ errors and the convergence orders at $T=1$, $\Delta t=\Delta x^{2}$, $\varepsilon=0.01$, $\gamma=2$.}
\begin{tabular}{lllllllll}
\hline
& \multicolumn{4}{c}{CHD4}& \multicolumn{4}{c}{CHD6}\vspace{1mm}\\
  \cline{2-9}
$N_{x}$ & $L_{{2}}$-error  & Rate   &  $L_{{\infty}}$-error   &  Rate
        & $L_{{2}}$-error  & Rate   &  $L_{{\infty}}$-error   &  Rate\\
\hline
$20$  &1.443(-6) &-   &4.964(-6)& -      &7.345(-7)  &-   &2.570(-6)  & -       \\
$40$  &9.017(-8) &4.00&2.901(-7)&4.10    &5.611(-9)  &7.03&2.045(-8)  &6.97     \\
$60$  &1.477(-8) &4.46&4.260(-8)&4.73    &6.868(-10) &5.18&2.440(-9)  &5.24     \\
$80$  &4.054(-9) &4.50&1.085(-8)&4.75    &1.130(-10) &6.27&3.989(-10) &6.29    \\
$100$ &1.498(-9) &4.46&4.013(-9)&4.46    &2.600(-11) &6.58&9.145(-11) &6.60   \\
$120$ &6.686(-10)&4.42&1.795(-9)&4.41    &7.633(-12) &6.72&2.676(-11) &6.74  \\
\hline
\end{tabular}
\end{table}
\begin{table}[htbp]
\centering
\caption{\label{table3:BurgersSolution1}
The exact solution \eqref{eq30:BurgersSolution1} of the one dimensional Burgers' equation: $L_{{2}}$ and $L_{{\infty}}$ errors and the convergence order at $T=1$, $\Delta t=\Delta x^{2}$, $\varepsilon=0.0001$, $\gamma=2$.}
\begin{tabular}{lllllllll}
\hline
& \multicolumn{4}{c}{CHD4}& \multicolumn{4}{c}{CHD6}\vspace{1mm}\\
  \cline{2-9}
$N_{x}$ & $L_{{2}}$-error  &  Rate   &  $L_{{\infty}}$-error   &  Rate
        &  $L_{{2}}$-error  &  Rate   &  $L_{{\infty}}$-error   &  Rate\\
\hline
$20$  &7.079(-10) &-   &3.162(-9)& -       &4.320(-10)  &-   &1.935(-9)   & -       \\
$40$  &1.068(-10) &2.73&6.701(-10)&2.24    &6.238(-12)  &6.11&3.990(-11)  &5.60     \\
$60$  &2.439(-11) &3.64&1.895(-10)&3.11    &1.447(-12)  &3.60&1.129(-11)  &3.11     \\
$80$  &7.662(-12) &4.02&6.861(-11)&3.53    &3.289(-13)  &5.15&2.921(-12)  &4.70    \\
$100$ &2.936(-12) &4.30&2.892(-11)&3.87    &9.124(-14)  &5.75&8.824(-13)  &5.37   \\
$120$ &1.292(-12) &4.50&1.359(-11)&4.14    &2.989(-14)  &6.12&3.057(-13)  &5.81  \\
\hline
\end{tabular}
\end{table}
\begin{table}[htbp]
\centering
\caption{\label{table4:BurgersSolution1}
The exact solution \eqref{eq30:BurgersSolution1} of the one dimensional Burgers' equation: $L_{{2}}$ error at $\Delta t=0.0001$, $\Delta x=0.025$, $\varepsilon=1$, $\gamma=2$ for various times.}
\begin{tabular}{llllll}
\hline
& \multicolumn{1}{l}{CHD4}& \multicolumn{1}{l}{CHD6}& \multicolumn{1}{l}{HAMQ \cite{Gao2014}} & \multicolumn{1}{l}{QBFE\cite{Xie2008}}& \multicolumn{1}{l}{HOC\cite{Jia2019}}\vspace{1mm}\\
  \cline{2-6}
 T  &$L_{{2}}$-error&$L_{{2}}$-error&$L_{{2}}$-error& $L_{{2}}$-error&$L_{{2}}$-error \\
\hline
$0.2$  &1.161(-7)  &5.583(-9)  &2.16(-4) &4.01(-4) &1.030(-7)  \\
$0.4$  &5.800(-9)  &7.776(-10) &2.77(-4) &4.98(-4) &8.041(-9)  \\
$0.6$  &3.462(-9)  &1.083(-10) &3.74(-4) &6.88(-4) &2.079(-10) \\
$0.8$  &8.551(-10) &1.509(-11) &4.05(-4) &8.01(-4) &9.751(-11) \\
$1.0$  &1.708(-10) &2.102(-12) &5.19(-4) &8.96(-4) &3.110(-11) \\
\hline
\end{tabular}
\end{table}
The initial and boundary conditions are obtained from the exact solution. Table \ref{table2:BurgersSolution1} gives the $L_{{2}}$ and $L_{{\infty}}$ errors and convergence orders with $\varepsilon=0.01$, $\gamma=2$. We can see that the CHD4 and CHD5 schemes
are capable of attaining precision exceeding fourth-order and sixth-order accuracy, respectively. Under mesh refinement, they exhibit excellent numerical stability.  When $\varepsilon=0.01$, $\gamma=2$, the convergence orders of CHD4 and CHD5 schemes still guarantee the theoretical accuracy as shown in Table \ref{table3:BurgersSolution1}. When $\varepsilon=1$, $\gamma=2$, $\Delta t=0.0001$, $\Delta x=0.025$, the $L_{{2}}$ error for various times is shown in Table \ref{table4:BurgersSolution1}. One can see that the numerical results of the CHD4 scheme are 3-6 orders of magnitude lower than that of the Refs. \cite{Gao2014,Xie2008} with the increase of time. As time progresses, while CHD4 exhibits larger $L_{{2}}$ errors than HOC\cite{Jia2019}, CHD6 demonstrates 1-2 orders of magnitude smaller errors compared to HOC\cite{Jia2019}.

We also consider another analytical solution of the equation \eqref{eq29:Burgers}, that is
\begin{equation}\label{eq31:BurgersSolution2}
u(x,t)=\frac{x}{t+t\sqrt{\frac{t}{t_{0}}}e^{\frac{x^{2}}{4\varepsilon t}}},\ t_{0}= e^{\frac{1}{8\varepsilon}},
\end{equation}
where $x\in[0,1.2],\ t\in(1,T]$.
The relevant initial and boundary conditions are
\begin{equation}
\begin{aligned}
&u(x,1)=\frac{x}{1+e^{(\frac{x^{2}}{4\varepsilon}-\frac{1}{16\varepsilon})}},\\
&u(0,t)=0,\ u(1.2,t)=\frac{b}{t+t\sqrt{\frac{t}{t_{0}}}e^{\frac{b^{2}}{4\varepsilon t}}}.
\end{aligned}
\end{equation}

\begin{table}[htbp]
\centering
\caption{\label{table5:BurgersSolution2}%
The exact solution \eqref{eq31:BurgersSolution2} of the Burgers' equation: $L_{{2}}$ and $L_{{\infty}}$ errors and the convergence orders at $T=2$, $\Delta t=\Delta x^{2}$, $\varepsilon=0.05$. }
\begin{tabular}{lllllllll}
\hline
& \multicolumn{4}{c}{CHD4}& \multicolumn{4}{c}{CHD6}\vspace{1mm}\\
  \cline{2-9}
$N_{x}$ & $L_{{2}}$-error & Rate   &  $L_{{\infty}}$-error   &  Rate
        & $L_{{2}}$-error &  Rate  &  $L_{{\infty}}$-error   & Rate\\
\hline
$20$  &2.498(-6) &-   &5.749(-6)& -     &7.784(-7)  &-   &1.662(-6) & -    \\
$40$  &1.332(-7) &4.23&2.867(-7)&4.33   &5.063(-9)  &7.26&1.565(-8) &6.73  \\
$60$  &2.462(-8) &4.16&5.116(-8)&4.25   &2.703(-10) &7.23&1.019(-9) &6.74  \\
$80$  &7.530(-9) &4.12&1.534(-8)&4.19   &3.463(-11) &7.14&1.411(-10)&6.87  \\
$100$ &3.021(-9) &4.09&6.087(-9)&4.14   &7.000(-12) &7.16&3.012(-11)&6.92   \\
$120$ &1.437(-9) &4.08&2.870(-9)&4.12   &1.893(-12) &7.17&8.486(-12)&6.95   \\
\hline
\end{tabular}
\end{table}
\begin{table}[htbp]
\centering
\caption{\label{table6:BurgersSolution2}%
The exact solution \eqref{eq31:BurgersSolution2} of the Burgers' equation: $L_{{2}}$ and $L_{{\infty}}$ errors and the convergence orders at $T=2$, $\Delta t=\Delta x^{2}$, $\varepsilon=0.005$.}
\begin{tabular}{lllllllll}
\hline
& \multicolumn{4}{c}{CHD4}& \multicolumn{4}{c}{CHD6}\vspace{1mm}\\
  \cline{2-9}
$N_{x}$ &  $L_{{2}}$-error  &  Rate   &  $L_{{\infty}}$-error   &  Rate
        &  $L_{{2}}$-error  &  Rate   &  $L_{{\infty}}$-error   &  Rate\\
\hline
$20$  &7.427(-4)&-   &2.184(-3)& -     &5.126(-4) &-   &1.516(-3) & -       \\
$40$  &9.387(-5)&2.98&3.679(-4)&2.57   &1.911(-5) &4.75&8.259(-5) &4.20     \\
$60$  &1.996(-5)&3.82&8.983(-5)&3.48   &1.642(-6) &6.05&9.077(-6) &5.49     \\
$80$  &6.295(-5)&4.01&2.677(-5)&4.21   &2.658(-7) &6.33&1.517(-6) &6.22    \\
$100$ &2.534(-6)&4.08&1.168(-5)&3.72   &6.276(-8) &6.47&3.551(-7) &6.51   \\
$120$ &1.200(-6)&4.10&5.558(-6)&4.07   &1.924(-8) &6.49&1.069(-7) &6.59   \\
\hline
\end{tabular}
\end{table}
\begin{table}[htbp]
\centering
\caption{\label{table7:BurgersSolution2}
The exact solution \eqref{eq31:BurgersSolution2} of the one dimensional Burgers' equation: $L_{{2}}$ and $L_{{\infty}}$ errors at $\varepsilon=0.005$ for various times.}
\begin{tabular}{l@{\hspace{2pt}}llllllll}
\hline
 & \multicolumn{2}{c}{CHD4}          & \multicolumn{2}{c}{HOC\cite{Jia2019}}
 & \multicolumn{2}{c}{CHD6}          & \multicolumn{2}{c}{FD6\cite{Sheng2023}}\vspace{1mm}\\
   \multicolumn{1}{l}{\scalebox{0.7}{($\Delta x,\Delta t$)}}
 & \multicolumn{2}{c}{\footnotesize(0.006,0.001)}  & \multicolumn{2}{c}{\footnotesize(0.005,0.001)}
 & \multicolumn{2}{c}{\footnotesize(0.01,0.01)}    & \multicolumn{2}{c}{\footnotesize(0.01,0.01)}\vspace{1mm}\\
  \cline{2-9}
 T  &           $L_{{2}}$-error  &  $L_{\infty}$-error  &
                $L_{{2}}$-error  &  $L_{\infty}$-error  &
                $L_{{2}}$-error  &  $L_{\infty}$-error  &
                $L_{{2}}$-error  &  $L_{\infty}$-error   \\
\hline
$1.7$  &2.187(-7)  &1.070(-6)  &8.04(-6) &2.56(-6) &1.132(-7)  &6.312(-7) &3.68(-7) &1.68(-6)      \\
$2.5$  &8.659(-8)  &3.764(-7)  &9.72(-6) &3.15(-6) &2.784(-8)  &1.362(-7) &3.41(-7) &1.24(-6)  \\
$3.0$  &5.602(-8)  &2.342(-7)  &9.90(-6) &3.12(-6) &1.436(-8)  &6.965(-8) &3.27(-7) &1.11(-6)   \\
$3.5$  &3.878(-8)  &1.539(-7)  &9.87(-6) &3.02(-6) &8.257(-9)  &4.133(-8) &3.15(-7) &1.01(-6)  \\
\hline
\end{tabular}
\end{table}
\begin{figure}[t]
	\centering
	\subfigure[]{
		\includegraphics[width=0.35\textwidth]{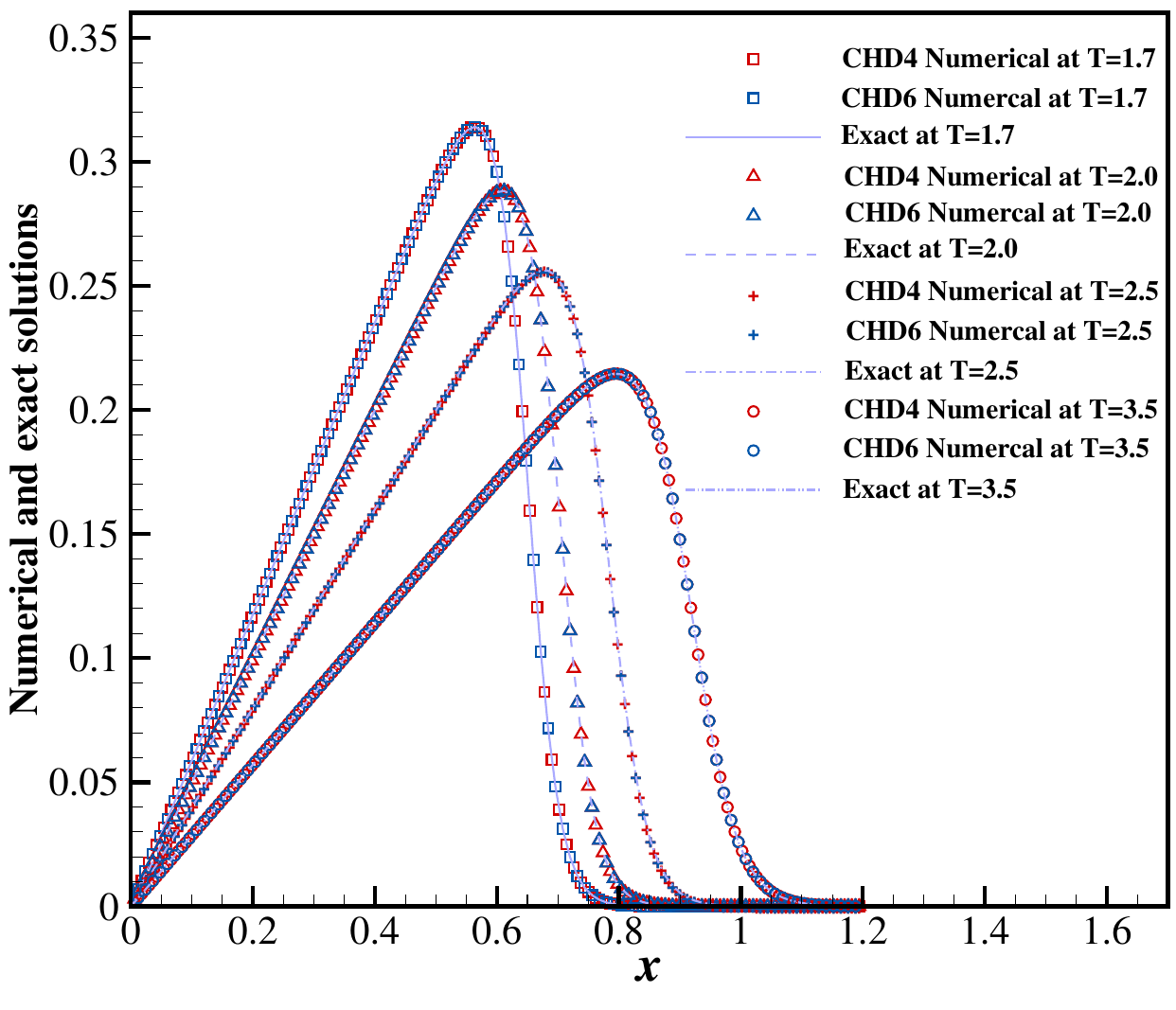}
	}
	\hspace{0.2cm}
	\subfigure[]{
		\includegraphics[width=0.35\textwidth]{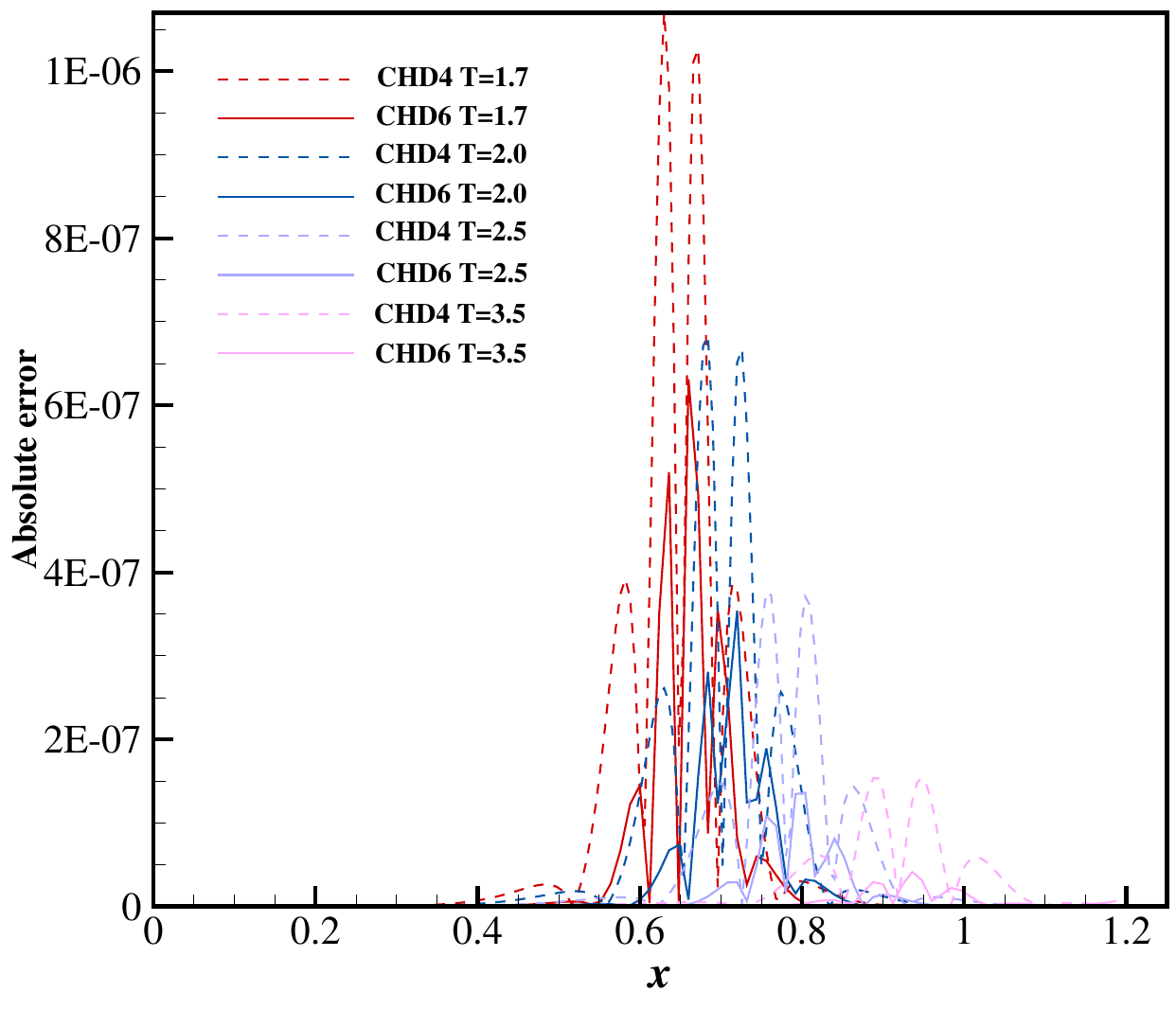}

	}
	\caption{\label{fig2:Burgers}Numerical results of the Burgers' equation with the exact solution \eqref{eq31:BurgersSolution2}
at T=1.7, 2, 2.5, 3.5 in the case $\varepsilon=0.005$, ($\Delta x$, $\Delta t$)=(0.006,0.001) for CHD4 and ($\Delta x$, $\Delta t$)=(0.01,0.01) for CHD6.}
\end{figure}
Table \ref{table5:BurgersSolution2} and Table \ref{table6:BurgersSolution2} show the $L_{{2}}$ and $L_{{\infty}}$ errors and the convergence orders for the CHD4 and CHD5 schemes with their boundary schemes when $\varepsilon=0.05$, $\varepsilon=0.005$, respectively. It can be seen that, with increasing mesh number, the results using the proposed schemes can achieve the analysis accuracy in space. It is not hard to see that, when $T=2$ and $\varepsilon$ is sufficiently small, the solutions of the equation exhibit a sharp transition near $x=0.6$, leading to compromised accuracy in solutions captured by our schemes.
Table \ref{table6:BurgersSolution2} shows the $L_{{2}}$ and $L_{{\infty}}$ errors for various times. We can seen that, the CHD4 scheme with  ($\Delta x$, $\Delta t$)=(0.006,0.001) is more accurate than the HOC scheme\cite{Jia2019} with ($\Delta x$, $\Delta t$)=(0.005,0.001),
the CHD6 scheme is more accurate than the FD6\cite{Sheng2023} under the same mesh sizes. Moreover, the CHD4 and HOC schemes have the same theoretical accuracy, as do the CHD6 and FD6 schemes. The numerical solutions and the absolute error of the numerical solutions corresponding to  Table \ref{table6:BurgersSolution2} are shown in Fig.\ref{fig2:Burgers}. These results well reflect the advantages of the high-order accuracy scheme in the calculation.

\subsubsection{\label{sec3.1.3:2D} Two dimensional convection-diffusion equation}
Consider the two dimensional convection-diffusion equation:
\begin{equation}\label{eq32:2DCDeq}
\frac{\partial u}{\partial t}+\frac{\partial (pu)}{\partial x}+\frac{\partial (qu)}{\partial y}
=\frac{1}{Re}(\frac{\partial^{2} u}{\partial x^{2}}+\frac{\partial^{2} u}{\partial y^{2}}),\ (x,y)\in[0,\pi]\times[0,\pi],\ t\in(0,T].
\end{equation}
The analytical solution to the equation is
\begin{equation}
u(x,y,t)=2e^{\frac{-2t}{Re}}cos(x)cos(y),
\end{equation}
with
\begin{equation}
\begin{aligned}
&p(x,y,t)=-e^{\frac{-2t}{Re}}cos(x)sin(y),\\
&q(x,y,t)=e^{\frac{-2t}{Re}}sin(x)cos(y).
\end{aligned}
\end{equation}
The initial and boundary conditions are taken from the exact solution.
\begin{table}[htbp]
\centering
\caption{\label{table8:2DCDeq}%
The two dimensional convection-diffusion equation: $L_{{2}}$ and $L_{{\infty}}$ errors and the convergence orders at $Re=1$, $T=0.5$.}
\begin{tabular}{lllllllll}
\hline
& \multicolumn{4}{c}{CHD4}& \multicolumn{4}{c}{CHD6}\vspace{1mm}\\
  \cline{2-9}
$N_{x}\times N_{y}$ &   $L_{{2}}$-error  &  Rate   & $L_{{\infty}}$-error   &  Rate
                    &   $L_{{2}}$-error  &  Rate   &  $L_{{\infty}}$-error   &  Rate\\
\hline
$10\times10$  &3.749(-5)&-   &8.396(-5)& -     &8.193(-6)  &-    &2.445(-5)  & -       \\
$20\times20$  &1.454(-6)&4.69&2.817(-6)&4.90   &3.217(-8)  &7.99 &1.070(-7)  &7.84      \\
$30\times30$  &2.436(-7)&4.41&4.736(-7)&4.40   &8.767(-10) &8.88 &2.832(-9)  &8.96     \\
$40\times40$  &7.204(-8)&4.23&1.406(-7)&4.22   &1.218(-10) &6.86 &2.815(-10) &8.02    \\
$50\times50$  &2.856(-8)&4.15&5.585(-8)&4.14   &3.812(-11) &5.21 &7.466(-11) &5.95   \\
$60\times60$  &1.353(-8)&4.10&2.647(-8)&4.09   &1.418(-11) &5.42 &2.613(-11) &5.76   \\
\hline
\end{tabular}
\end{table}
\begin{table}[htbp]
\centering
\caption{\label{table9:2DCDeq}%
The two dimensional convection-diffusion equation: $L_{{2}}$ and $L_{{\infty}}$ errors and the convergence orders at $Re=10$, $T=0.5$.}
\begin{tabular}{lllllllll}
\hline
&\multicolumn{4}{c}{CHD4}& \multicolumn{4}{c}{CHD6}\vspace{1mm}\\
  \cline{2-9}
$N_{x}\times N_{y}$ &   $L_{{2}}$-error  &  Rate   & $L_{{\infty}}$-error    &  Rate
                    &   $L_{{2}}$-error  &  Rate   &  $L_{{\infty}}$-error   &  Rate\\
\hline
$10\times10$  &1.174(-4)&-   &7.087(-4)& -     &2.706(-4)  &-   &8.772(-4)  & -       \\
$20\times20$  &1.157(-6)&7.55&3.711(-6)&7.58   &1.506(-5)  &7.49&5.336(-6)  &7.36      \\
$30\times30$  &1.600(-7)&4.88&3.846(-7)&5.59   &4.735(-8)  &8.53&1.638(-7)  &8.59     \\
$40\times40$  &4.939(-8)&4.09&1.054(-7)&4.50   &3.790(-9)  &8.78&1.296(-8)  &8.82    \\
$50\times50$  &1.977(-8)&4.10&4.136(-8)&4.19   &5.163(-10) &8.93&1.743(-9)  &8.99   \\
$60\times60$  &9.394(-9)&4.08&1.957(-8)&4.10   &9.889(-11) &9.06&3.280(-10) &9.16   \\
\hline
\end{tabular}
\end{table}

To balance accuracy and efficiency, The grid is systematically coarsened by a factor of $2$ in each spatial direction relative to the one dimensional reference mesh. Table \ref{table8:2DCDeq} and Table \ref{table9:2DCDeq} show the $L_{{2}}$ and $L_{{\infty}}$ errors and the convergence orders for the CHD4 and CHD5 schemes when $Re=1$, $Re=10$, respectively. With the mesh refinement, the CHD4 scheme can preserve fourth-order accuracy when $Re=1$ and $Re=10$. Although the accuracy of the CHD6 scheme degrades, it still achieves better than fifth-order accuracy with the mesh refinement when
$Re=1$. When $Re=10$, the CHD6 scheme can achieves better than sixth-order accuracy.

\subsection{\label{sec3.2:application} Application: Double-diffusive Convection}
\begin{figure}[htbp] 
	\centering
	\includegraphics[width=0.2\textwidth]{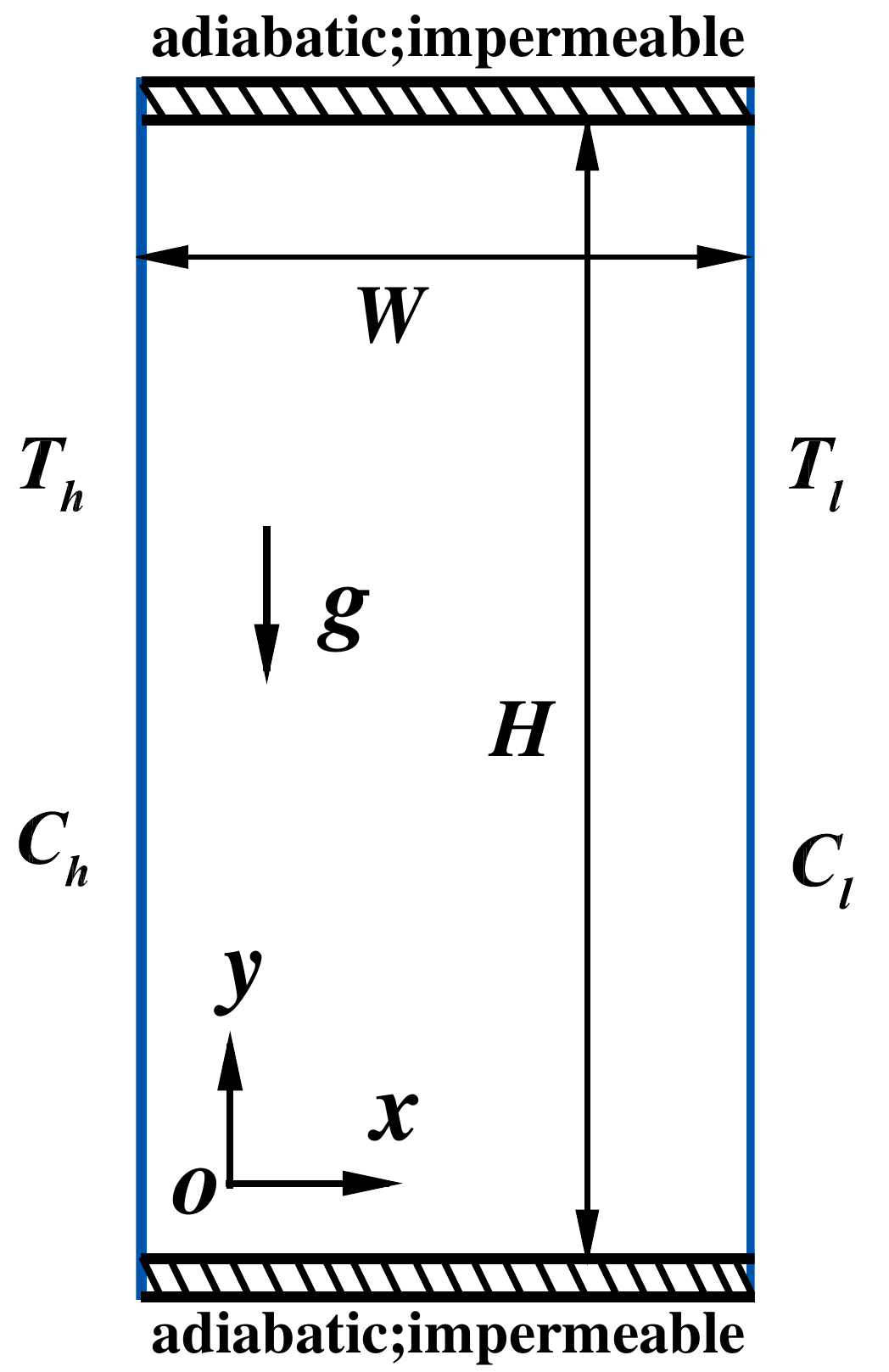}
	\caption{Enclosure flow configuration, coordinates, and boundary conditions.}\label{fig:enclosure}
\end{figure}

In this section, we consider the double-diffusive convection in a rectangular cavity with an aspect ratio $A=H/W$, where $H$ and $W$ are the height and the width of the cavity, respectively, as shown in Fig.\ref{fig:enclosure}. The top and bottom of the cavity are considered to be adiabatic and impermeable for mass transfer. Uniform temperature and concentration differences are imposed across the vertical walls, where the temperature of left wall $T_{h}$ is higher than that of the right wall $T_{l}$, similarly, the concentration of the left wall $C_{h}$ is higher than that of the right wall $C_{l}$. The associated dimensionless boundary conditions can be written as
\begin{eqnarray}\label{eq2:Boundary}
&&\psi=u=v=0,\  T=0.5,\quad                   C=0.5, \quad\                          at \quad x=0,     \quad\ 0\leq y\leq A, \label{eq2:1}\\
&&\psi=u=v=0,\  T=-0.5,\                      C=-0.5,\                               at \quad x=1,     \quad\ 0\leq y\leq A, \label{eq2:2}\\
&&\psi=u=v=0,\  \frac{\partial T}{\partial y}=0,\quad  \frac{\partial C}{\partial y}=0,\ \quad  at \quad  y=0,A, \ 0\leq x\leq1.\label{eq2:3}
\end{eqnarray}

The boundary conditions \eqref{eq2:3} of temperature and concentration on the top and bottom walls are Neumann boundary conditions.
A approximation formula \cite{Qin2014} of the temperature boundary condition $\frac{\partial T}{\partial y}=0$ is considered. It is
\begin{equation}\label{eq24:temperature}
\begin{aligned}
&T_{w}=\frac{1}{25}(48T_{1}-36T_{2}+16T_{3}-3T_{4})+O(\Delta y^{5}),
\end{aligned}
\end{equation}
where $1, 2, 3, 4$ denote the first four neighboring internal points on normal through the boundary $w$ ($y=1$ or $y=A$). Considering
Similarly, the discretization scheme of the concentration boundary condition $\frac{\partial C}{\partial y}=0$ can also be given by
\begin{equation}\label{eq25:concentration}
\begin{aligned}
&C_{w}=\frac{1}{25}(48C_{1}-36C_{2}+16C_{3}-3C_{4})+O(\Delta y^{5}).
\end{aligned}
\end{equation}

Due to the lack of the physical boundary conditions for the vorticity $\omega$ on all the walls, numerical boundaries of vorticity are
updated by the scheme\cite{Spotz1998}\\
\begin{equation}\label{eq26:vorticity}
\begin{aligned}
&\frac{h}{21}(6\omega_{w}+4\omega_{1}-\omega_{2})+O(h^{4})=
         \frac{1}{14h}(15\psi_{w}-16\psi_{1}+\psi_{2})+\nu_{w}.
\end{aligned}
\end{equation}
where $h$ represents $\Delta x$ or $\Delta y$. $\nu_{w}$ is the tangential wall velocity which equals zero for no slip boundary condition.

As we all know, the compositional buoyancy force is primarily dominated in the flow when $\lambda>1$, whereas the thermal buoyancy force is primarily dominated for $\lambda<1$. The thermal and compositional buoyancy forces are exactly balanced at $\lambda=1$, meanwhile, the flow dynamics exhibit maximal complexity.
It is conventional to describe the heat and mass transfer characteristics in a cavity using the average Nusselt number and Sherwood number along the hot wall. They are introduced by
\begin{equation}\label{Nusselt}
  Nu_{av}=\frac{1}{A}\int^{A}_{0}(\frac{\partial T}{\partial x})\mid_{x=0}dy, \quad
  Sh_{av}=\frac{1}{A}\int^{A}_{0}(\frac{\partial C}{\partial x})\mid_{x=0}dy,
\end{equation}
respectively.

At the initial instant $t=0$, the temperature and concentration at the vertical sidewalls are abruptly altered, and they are maintained thereafter, which is a set of complex competitive systems.

\subsubsection{\label{sec3.2.1} Verification of code}
To verify numerical accuracy and computational efficiency, we present the results of double-diffusive convection in a rectangular cavity using the CHD4 scheme. In this section, the values of the parameters in \eqref{eq2:parameter} and aspect ratio $A$ are listed in Table \ref{table:VerifyCode}.
\begin{table}[htbp]
\centering
\caption{\label{table:VerifyCode}%
Range of control parameters in subsection 3.2.1.}
\begin{tabular}{llllll}
\hline
Flow state                 & $Pr$    & $Le$    &$Ra$       &  $A$    & $\lambda$ \\
\hline
\scalebox{0.9}{Period flow}& $1.0$   & $2.0$   &$10^{5}$   &  $2.0$  & $1.0$ \\

\scalebox{0.9}{Steady flow}& $1.0$   & $2.0$   &$10^{5}$   &  $2.0$  & $0.8$ \\
\cline{2-6}
                           & $1.0$   & $2.0$   &$10^{5}$   &  $2.0$  & $1.3$ \\
\hline
\end{tabular}
\end{table}

{\bf Part 1. Comparative analysis of periodic solutions. }

To validate the numerical accuracy and computational efficiency of the proposed algorithm, we present the results of the period flow of double-diffusive convection in a cavity. The values of control parameters are listed in Table \ref{table:VerifyCode}.

We compare the present solutions with spectral solutions, \cite{Morega1996} finite element solutions \cite{Nishimura1998} and other finite difference results in the study. \cite{Yan2022}
Table \ref{table:chooseMesh} shows the dimensionless period of oscillation, the stream function extremum, $\mid\psi_{\max}\mid$, $\mid\psi_{\min}\mid$ at different grid sizes.
It can be seen that our results are agree well with ones in the literature,\cite{Morega1996,Nishimura1998,Yan2022}  and a further refinement of grid sizes from $30\times40$ to $40\times80$ does not lead to a significant effect on the calculated solutions by the present algorithm.
\begin{table}[htbp]
\centering
\caption{\label{table:chooseMesh}%
Comparison of the different numerical method for $Pr=1$, $Le=2$, $Ra=10^{5}$, $A=2$ and $\lambda=1$.}
\begin{tabular}{llllllll}
\hline
Method                  & Present      & Present    & Present    &FD\cite{Yan2022}& FD\cite{Yan2022}
                                                                 &FEM\cite{Nishimura1998}&  SM\cite{Morega1996}\\
\hline
$N_{x}\times N_{y}$     & $30\times40$ &$30\times60$&$40\times80$ & $30\times40$  &$30\times60$   &$30\times40$ &$40\times80$ \\
Period                  & $0.0490$     &$0.0492$    &$0.0492$     &$0.0489$       &$0.0492$       &$0.0497$      &$0.0494$\\
Max$\mid\psi_{max}\mid$ & $26.866$     &$26.838$    &$26.829$     &$26.72$        &$26.77$        &$26.7$        &$26.8$\\
Min$\mid\psi_{max}\mid$ & $12.682$     &$12.662$    &$12.693$     &$12.64$        &$12.65$        &$12.9$        &$12.7$\\
Max$\mid\psi_{min}\mid$ & $5.514$      &$5.568$     &$5.553$      &$5.57$         &$5.58$         &$5.76$        &$5.52$\\
Min$\mid\psi_{min}\mid$ & $0.313$      &$0.316$     &$0.329$      &$0.327$        &$0.321$        &$0.351$       &$0.333$\\
\hline
\end{tabular}
\end{table}

\begin{figure}[t]
	\centering
	\subfigure[]{
		\includegraphics[width=0.11\textwidth]{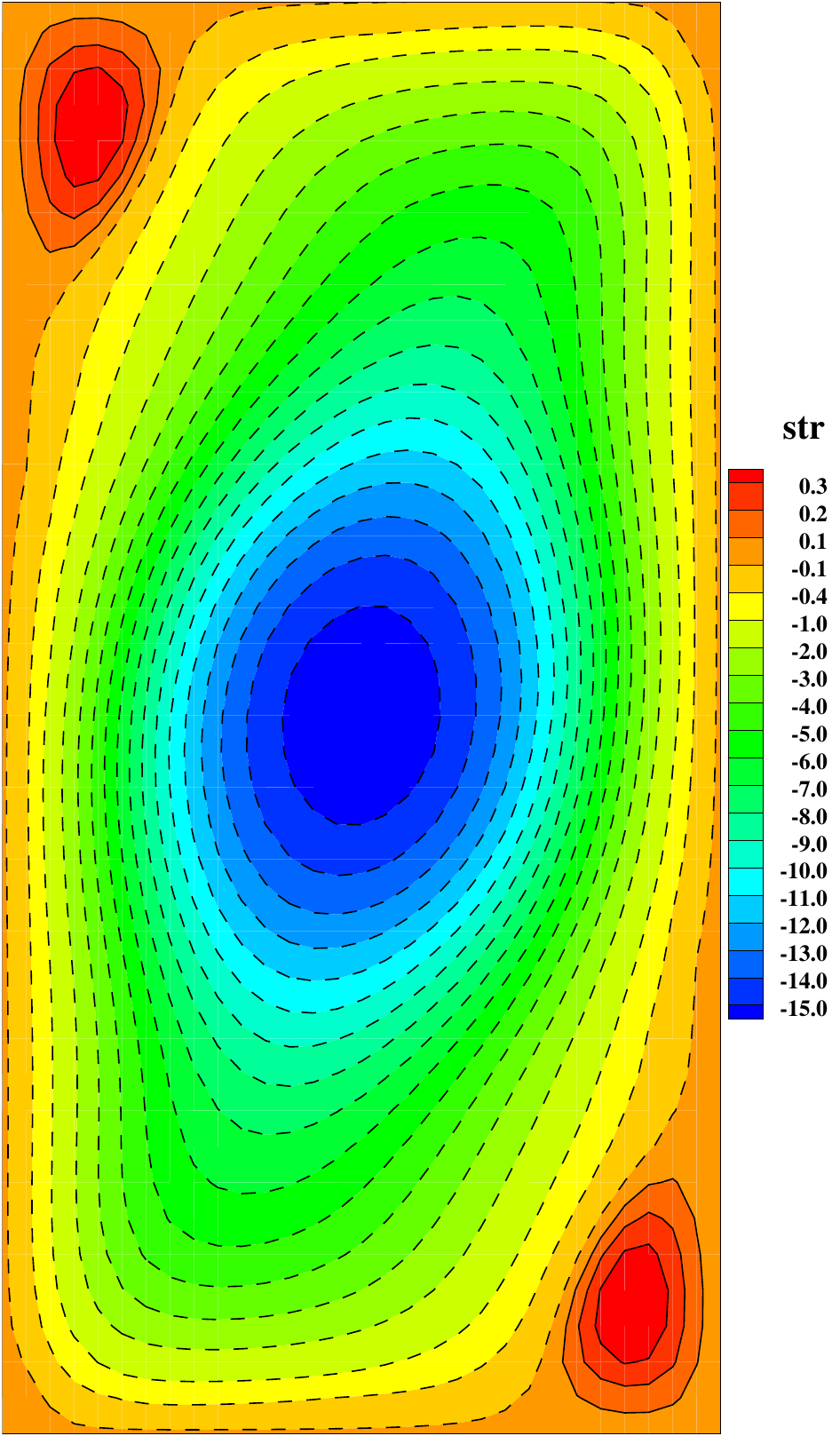}}
	\subfigure[]{
        \includegraphics[width=0.11\textwidth]{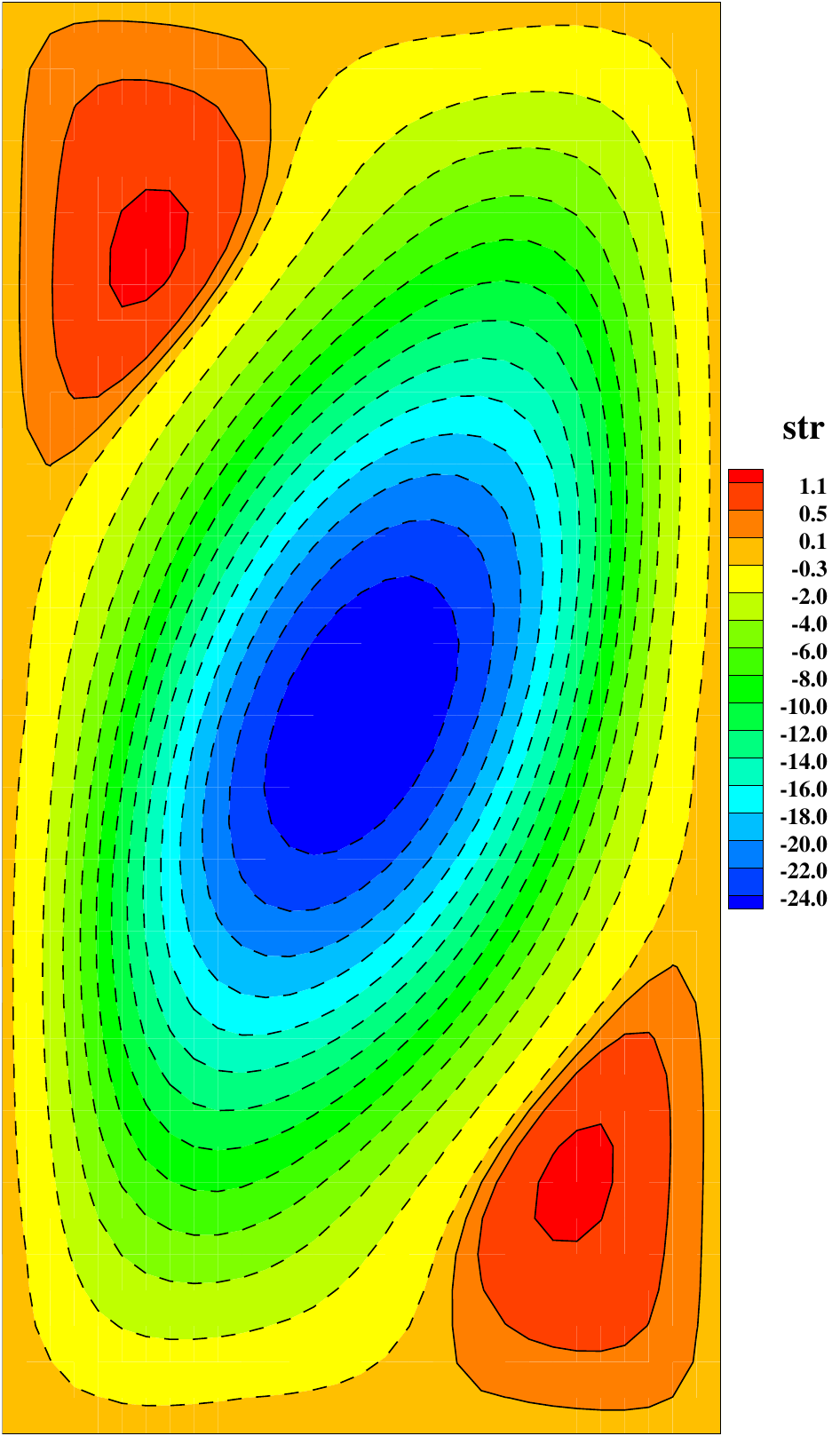}}
	\subfigure[]{
		\includegraphics[width=0.11\textwidth]{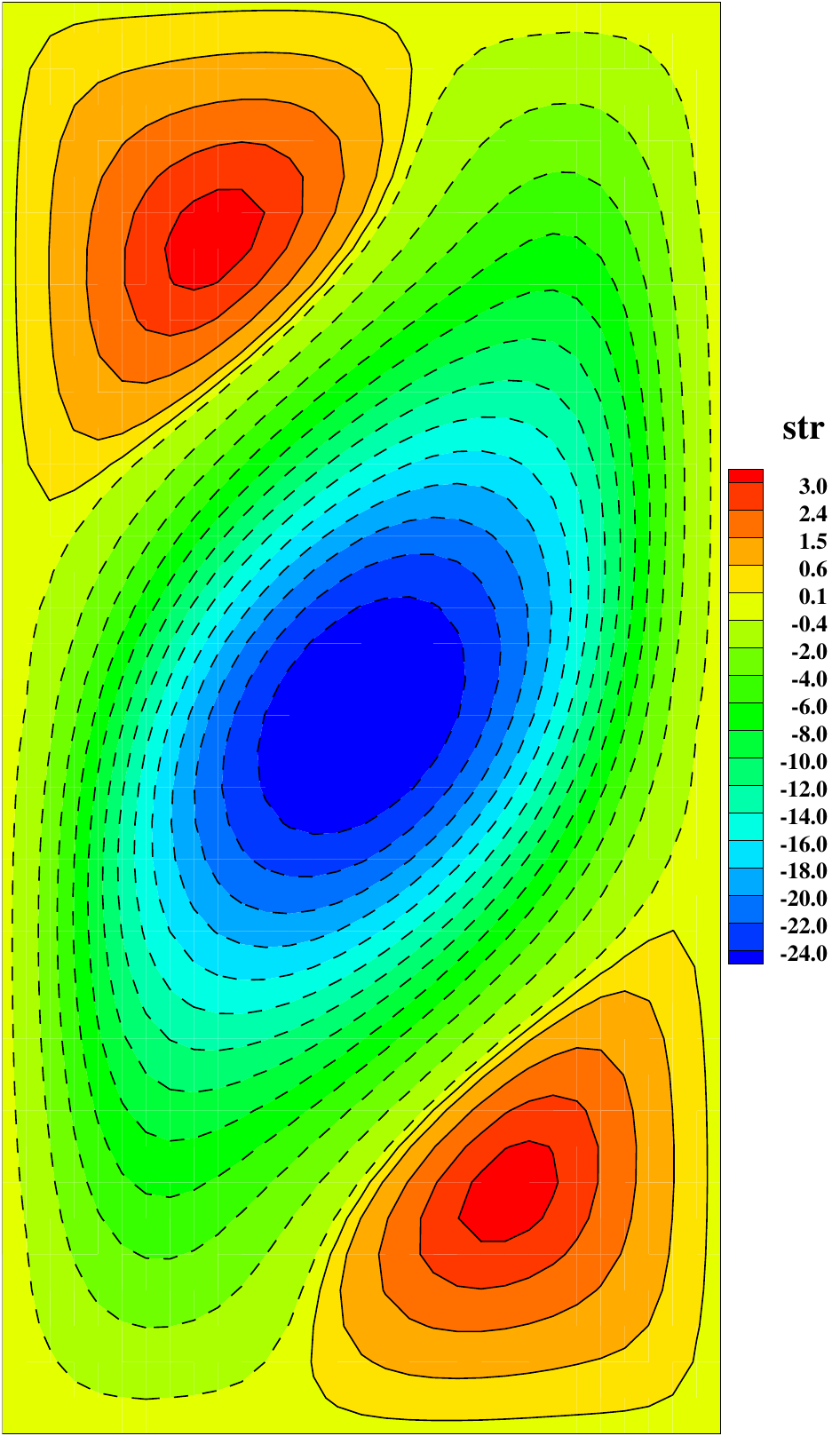}}
	\subfigure[]{
		\includegraphics[width=0.11\textwidth]{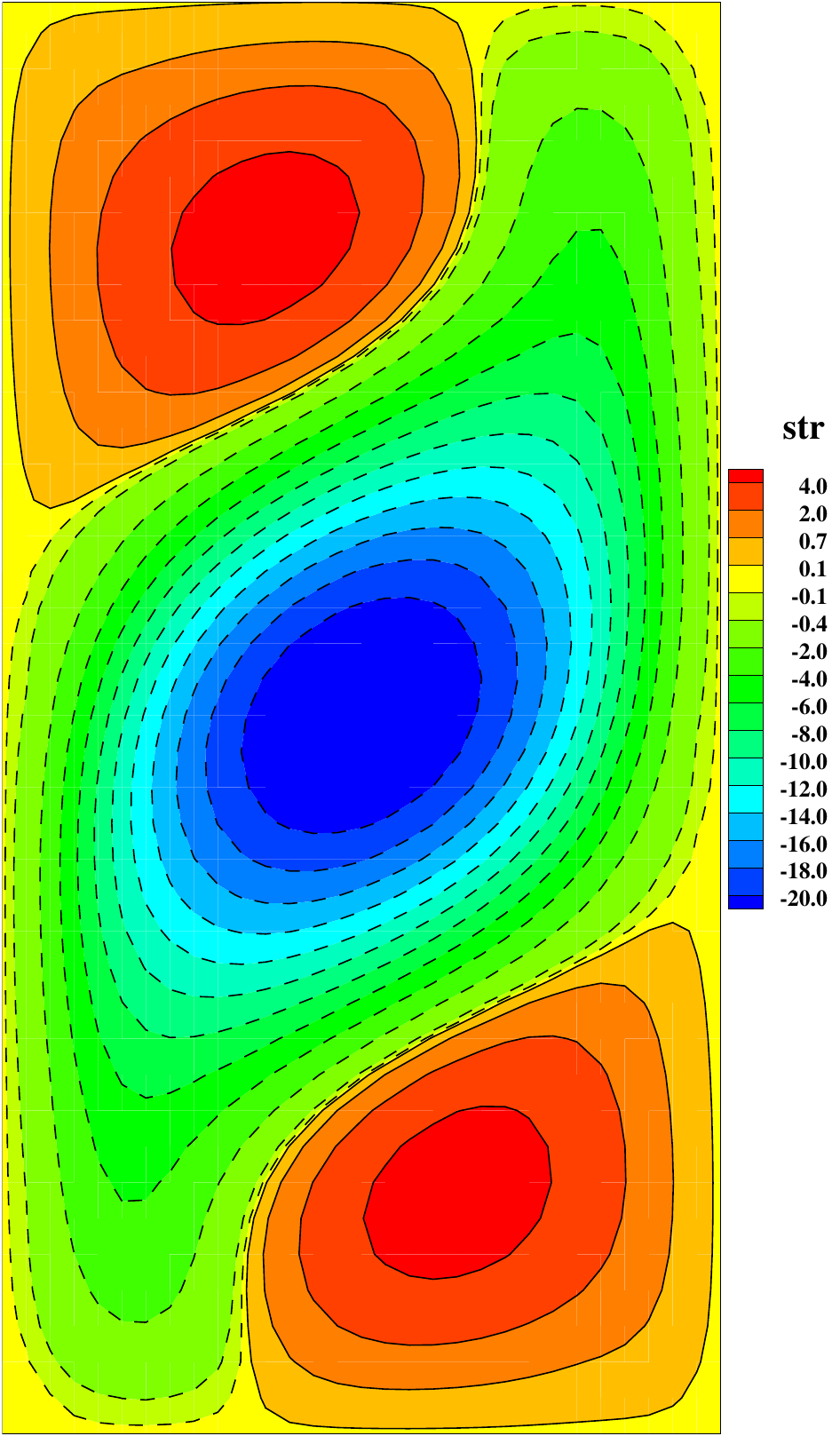}}
    \subfigure[]{
		\includegraphics[width=0.11\textwidth]{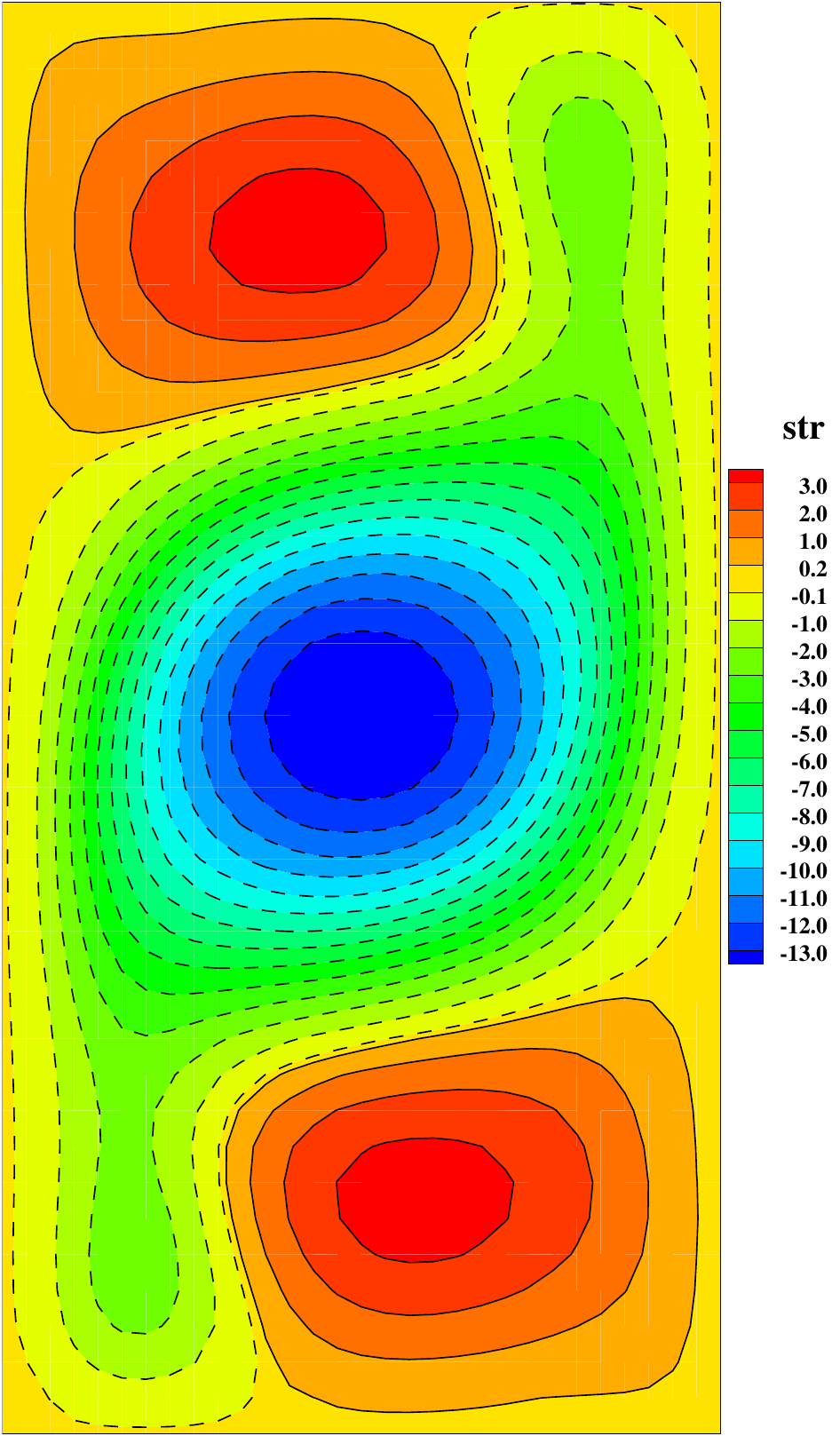}}
	\subfigure[]{
        \includegraphics[width=0.11\textwidth]{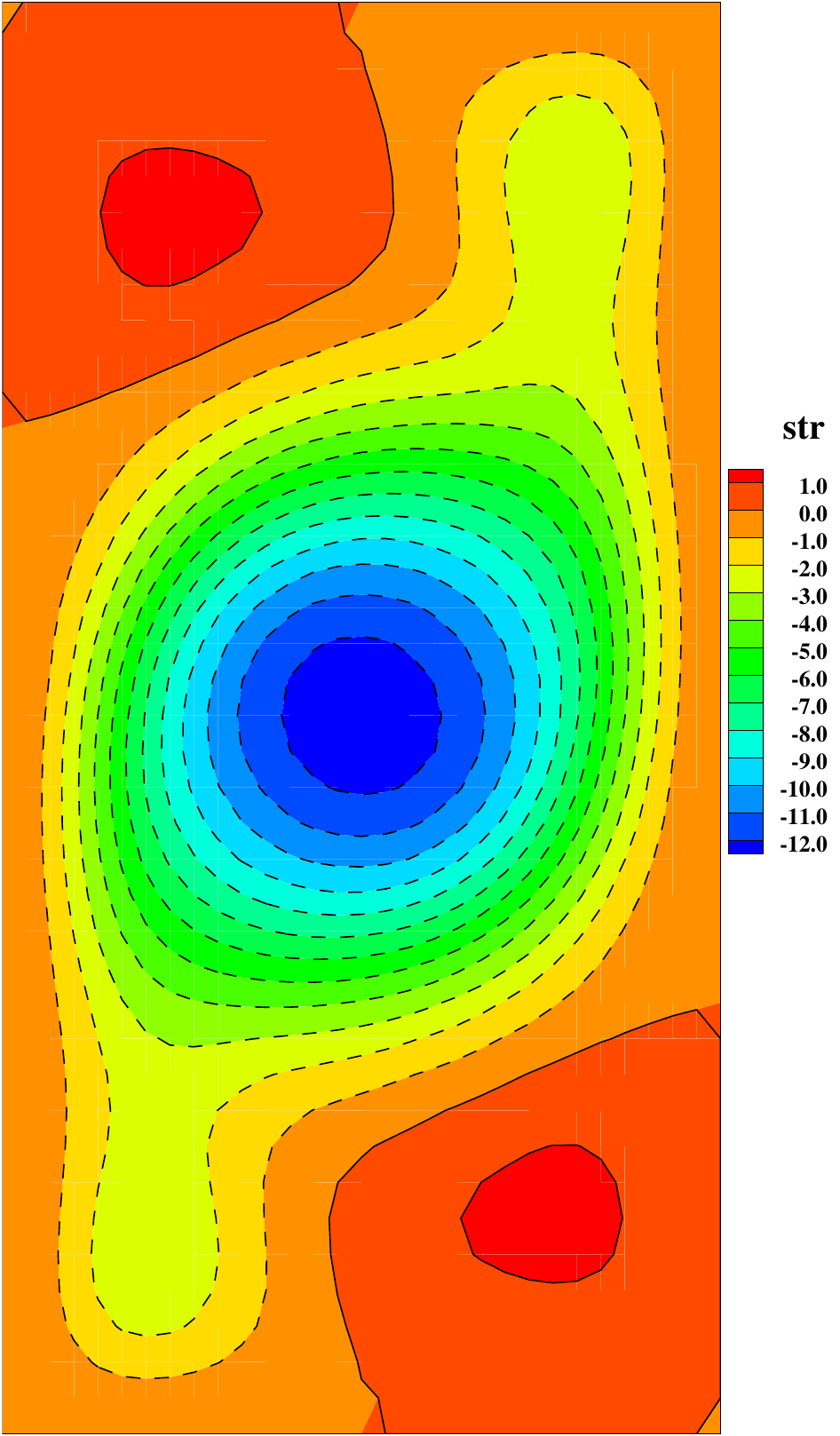}}
	\subfigure[]{
		\includegraphics[width=0.11\textwidth]{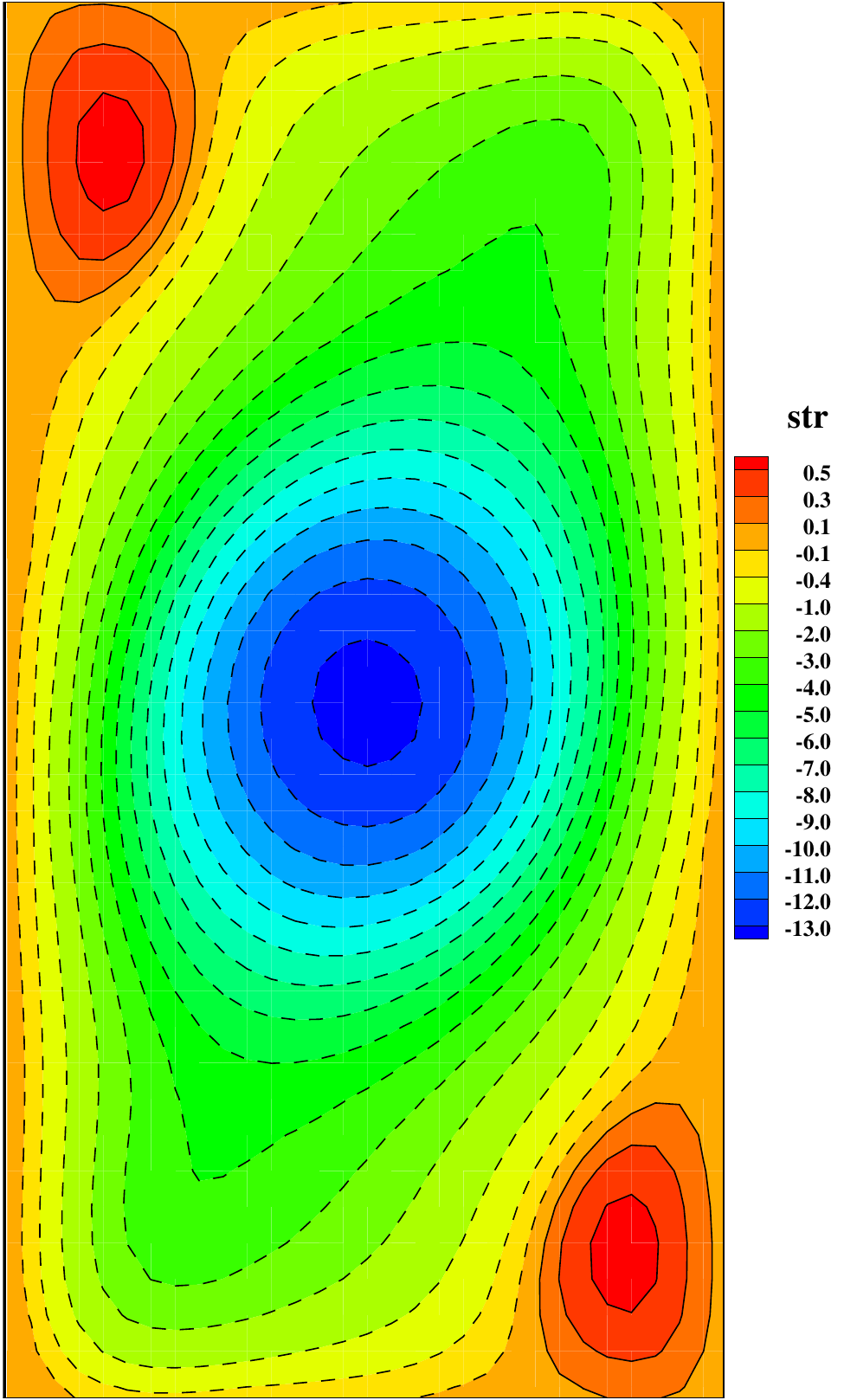}}
	\subfigure[]{
		\includegraphics[width=0.11\textwidth]{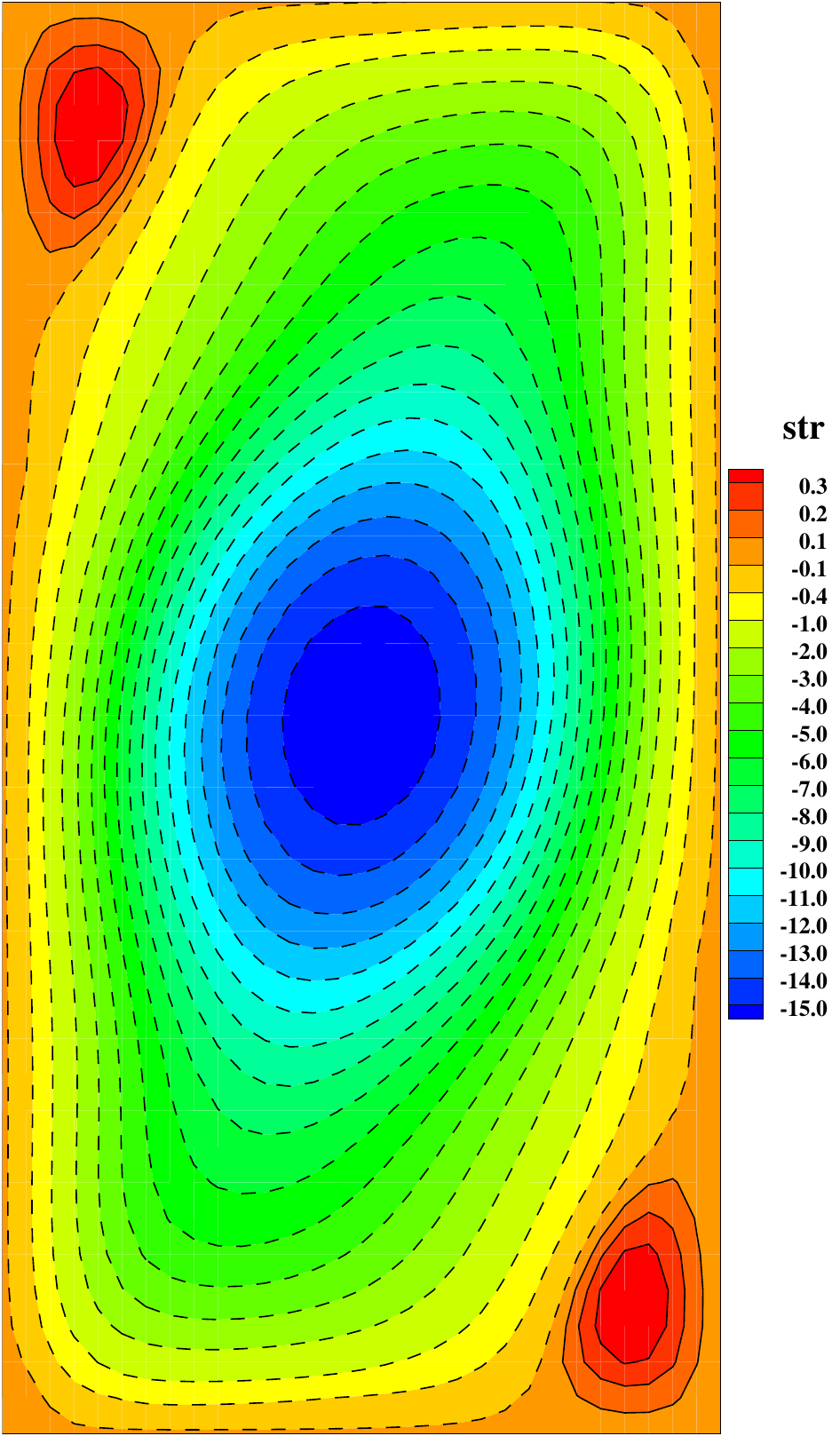}}
	\caption{\label{fig:fq1Str} Time history of streamline contours for $Pr=1$, $Le=2$, $Ra=10^{5}$, $A=2$, $\lambda=1.0$. }
\end{figure}

Fig. \ref{fig:fq1Str}(a) to (g) shows a whole periodic evolution of the flow field. In the all contour lines, the dashed line denotes that the rotation of the vortex is clockwise and the value is negative, and the solid line denotes that the rotation of the vortex is counter-clockwise and the value is positive.

In Fig. \ref{fig:fq1Str}(a), one can see that a large primary vortex occupies the middle part of the cavity and the top left corner and bottom right corner of the cavity appear small secondary vortices, respectively.
Then the two secondary vortices move from the top left and bottom right corners along the vertical wall towards the center and squeeze the primary vortex, as shown in Fig. \ref{fig:fq1Str}(b)-(e).
In Fig. \ref{fig:fq1Str}(f), the two secondary vortices back to the top left and bottom right corners, respectively.
A complete periodic motion is formed in Fig. \ref{fig:fq1Str}(g). A new cycle starts in Fig. \ref{fig:fq1Str}(h).

{\bf Part 2. Comparative analysis of steady solutions.}

\begin{figure}[t]
	\centering
	\subfigure[$\lambda=0.8$]{
		\includegraphics[width=0.4\textwidth]{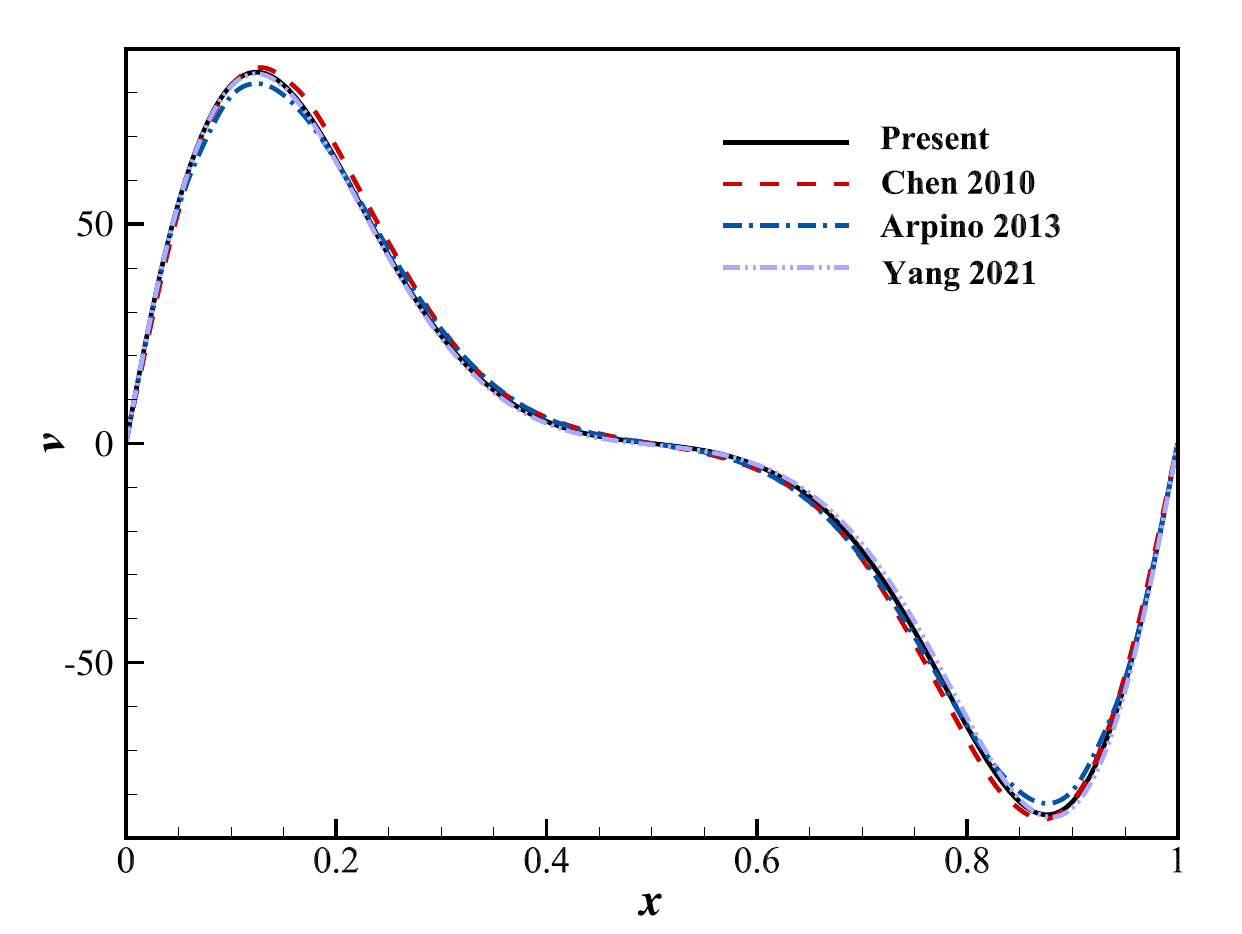}
	}
	\hspace{0.2cm}
	\subfigure[$\lambda=1.3$]{
		\includegraphics[width=0.4\textwidth]{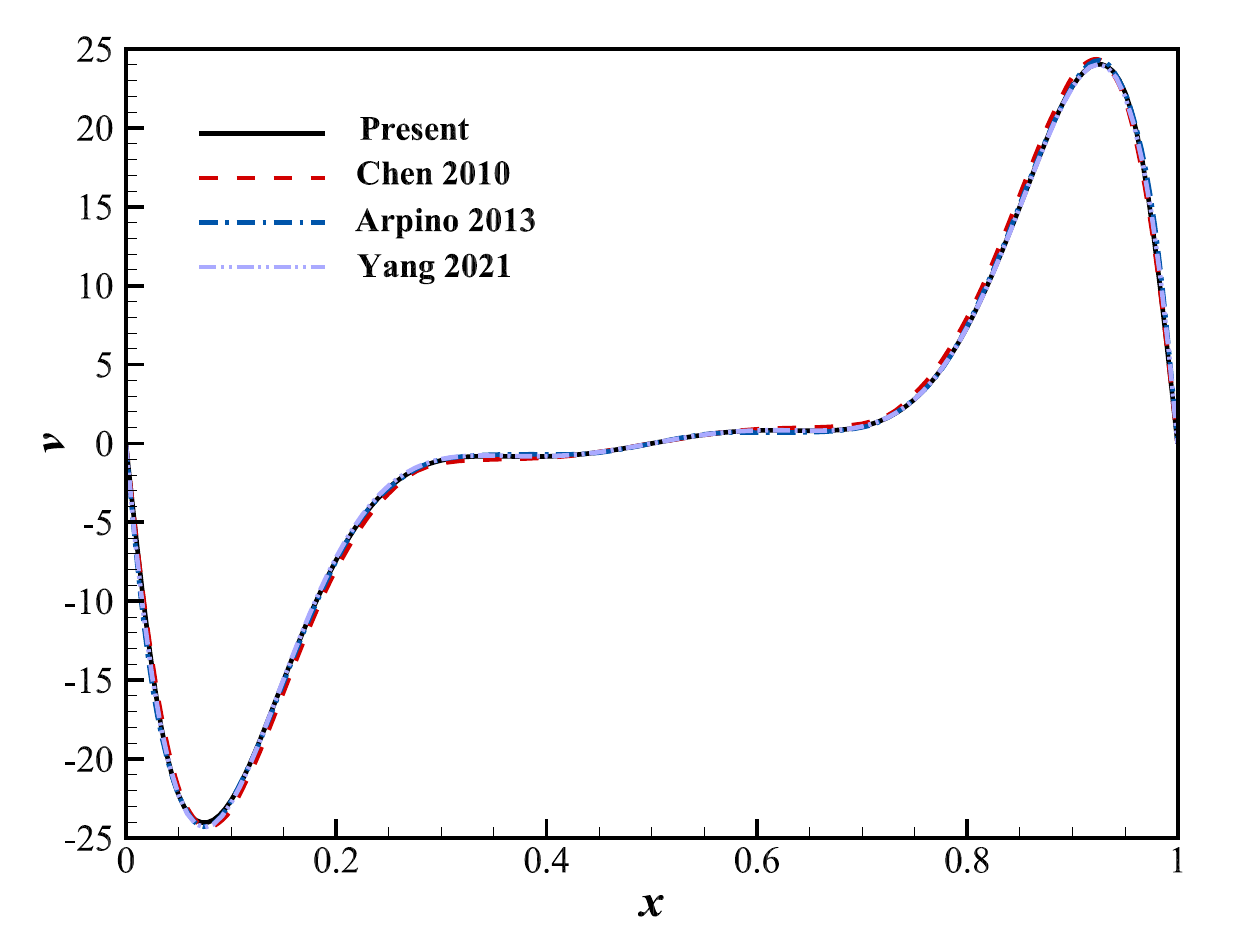}

	}
	\caption{\label{fig:Velocity}Vertical velocity profiles at mid-height of the cavity for (a)$\lambda=0.8$ and (b)$\lambda=1.3$.}
\end{figure}

\begin{figure}[t]
	\centering
	\subfigure[$\lambda=0.8$]{
		\includegraphics[width=0.4\textwidth]{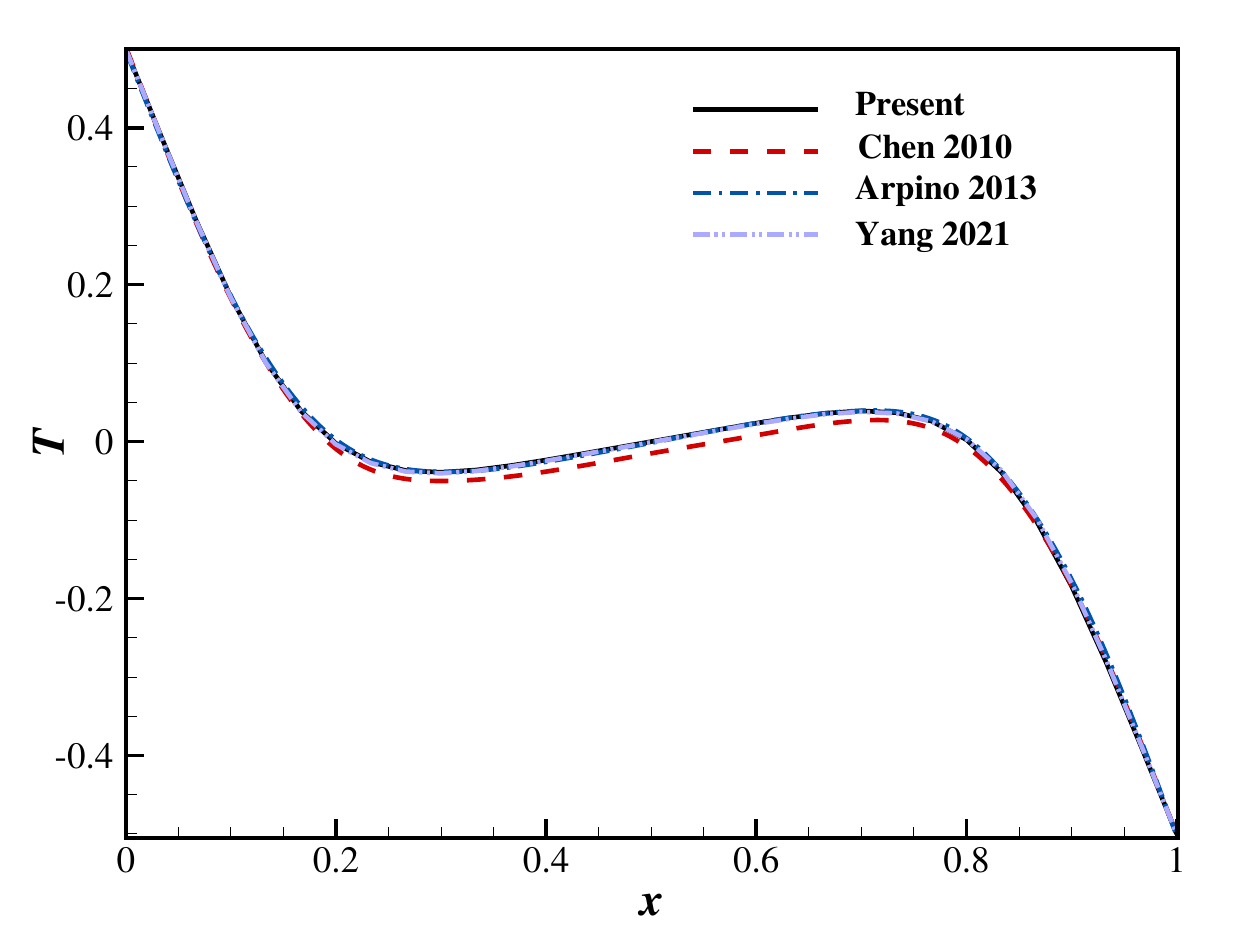}
	}
	\hspace{0.2cm}
	\subfigure[$\lambda=1.3$]{
		\includegraphics[width=0.4\textwidth]{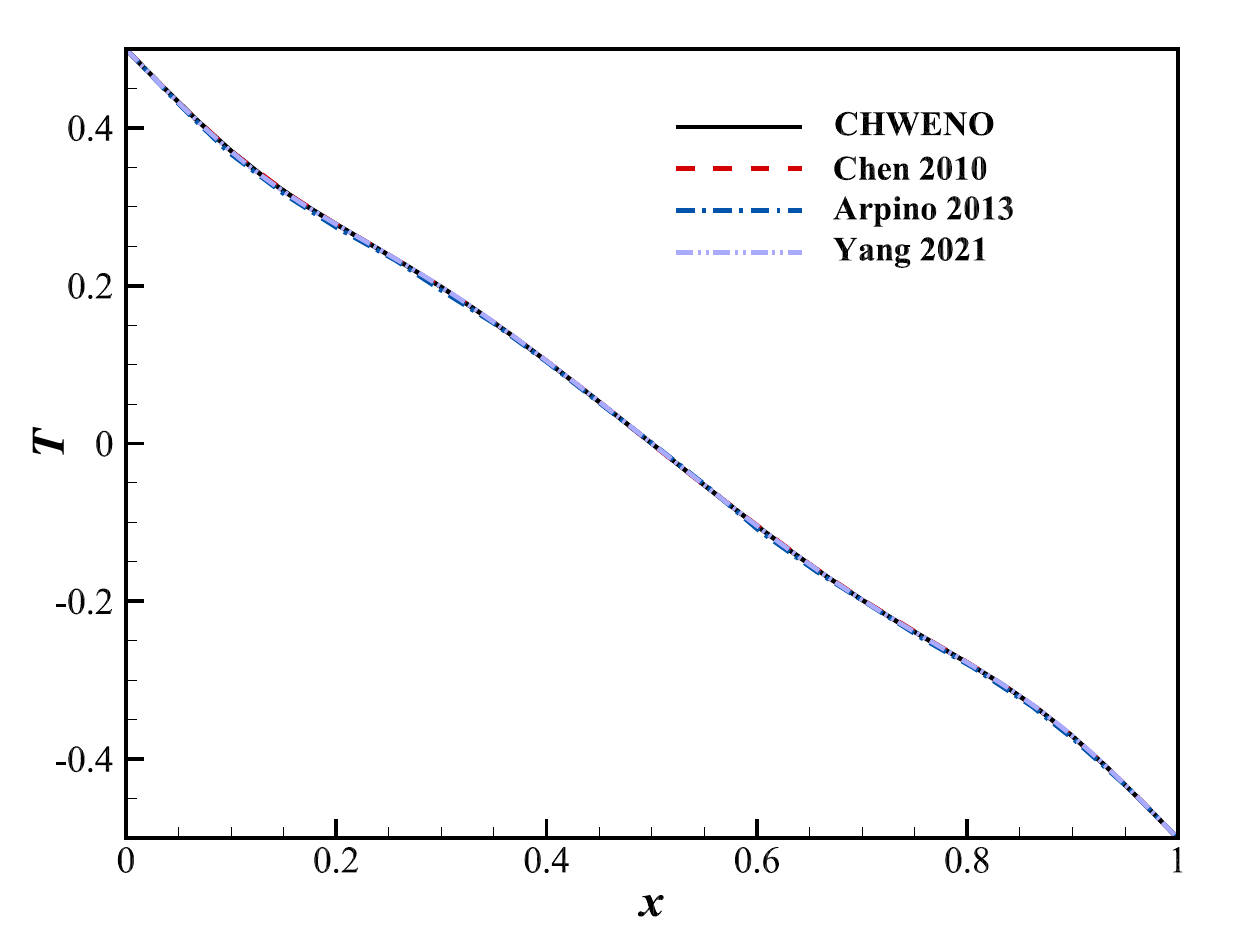}

	}
	\caption{\label{fig:Temperature}Temperature profiles at mid-height of the cavity for (a)$\lambda=0.8$ and (b)$\lambda=1.3$.}
\end{figure}

\begin{figure}[t]
	\centering
	\subfigure[$\lambda=0.8$]{
		\includegraphics[width=0.4\textwidth]{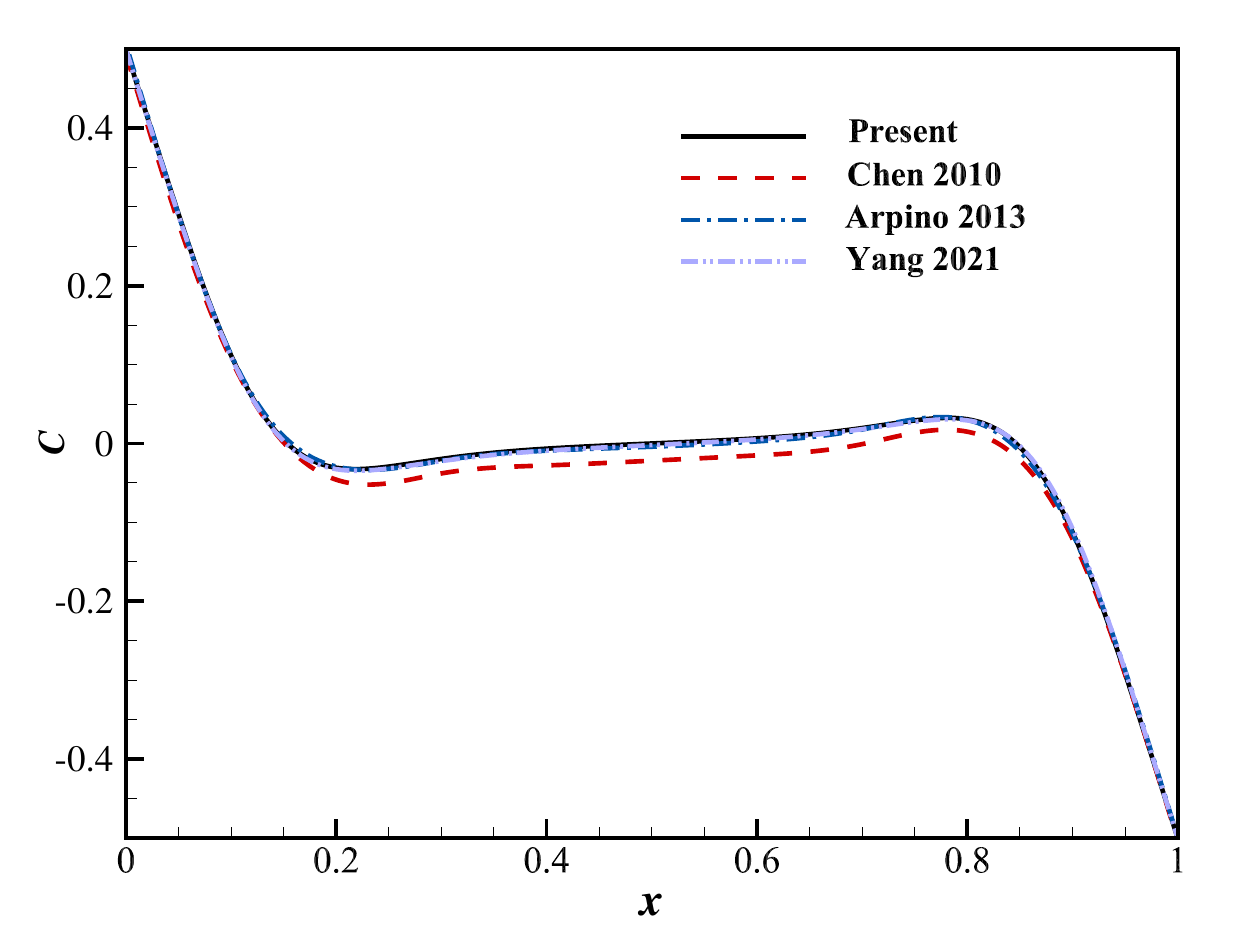}
	}
	\hspace{0.2cm}
	\subfigure[$\lambda=1.3$]{
		\includegraphics[width=0.4\textwidth]{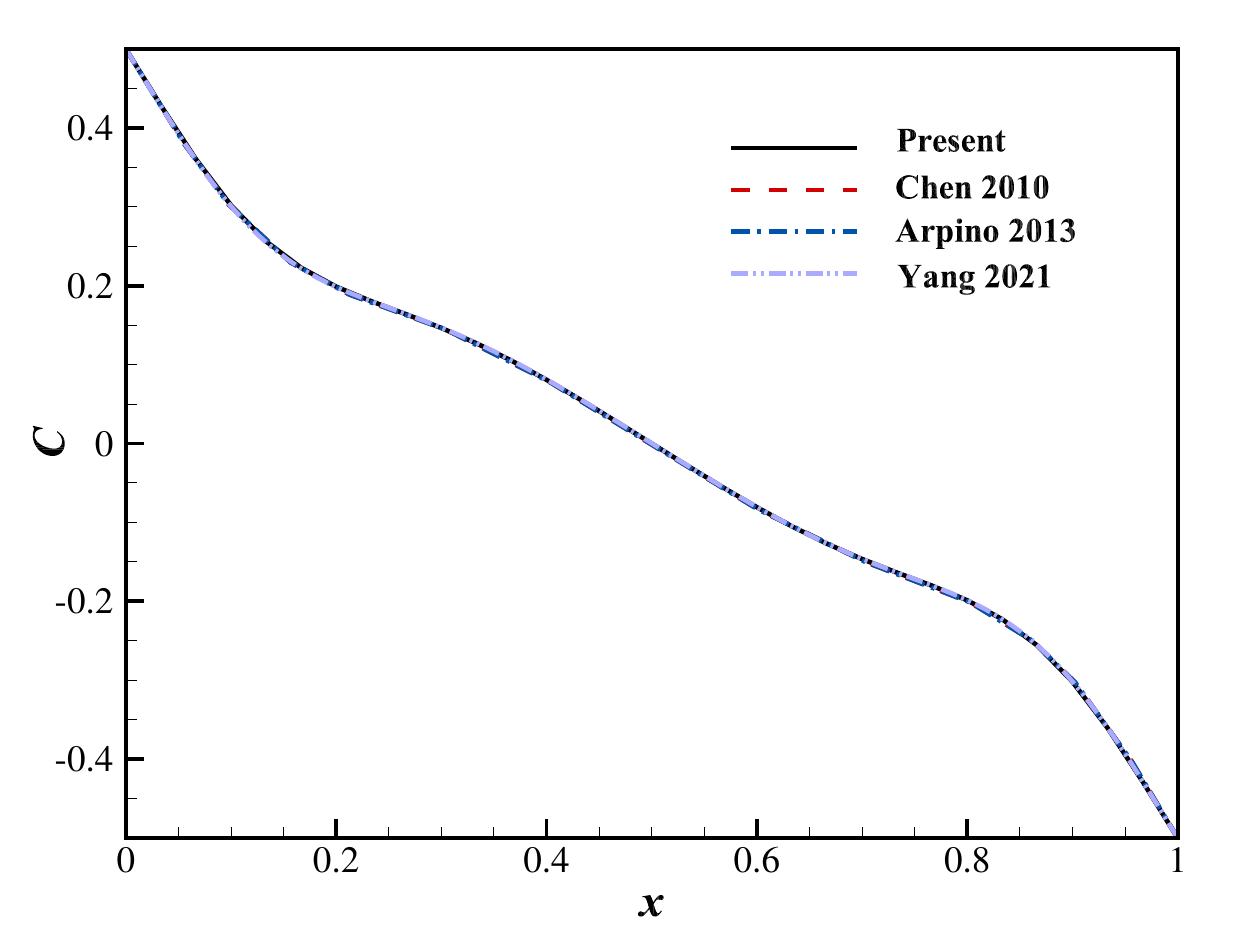}

	}
	\caption{\label{fig:Concentration}Concentration profiles at mid-height of the cavity for (a)$\lambda=0.8$ and (b)$\lambda=1.3$.}
\end{figure}

We further verify the numerical accuracy and efficiency for the steady flow of double-diffusive convection. Table \ref{table:VerifyCode} gives the values of control parameters.
Fig. \ref{fig:Velocity}, Fig. \ref{fig:Temperature} and Fig. \ref{fig:Concentration} present comparative results for velocity $v(x,\frac{A}{2})$, temperature $T(x,\frac{A}{2})$  and concentration $C(x,\frac{A}{2})$ of the cavity, respectively. These results demonstrate excellent agreement with solutions reported in previous studies.\cite{Arpino2013,Chen2010,Zhao2021}

When $\lambda=0.8$, the flow field is temperature dominated. As shown in Fig. \ref{fig:Velocity}(a), the vertical velocity profiles exhibit positive and negative peaks adjacent to the left and right walls, respectively.
Whereas the trend of the vertical velocity is opposite when $\lambda=1.3$ as shown in Fig. \ref{fig:Velocity}(b).
At this case, the flow field is concentration dominated.

In Fig. \ref{fig:Temperature} and Fig. \ref{fig:Concentration}, it can be seen that temperature and concentration profiles are almost flat in the center of the cavity for $\lambda=0.8$, however they exhibit a linear trend for $\lambda=1.3$.
It is worth emphasizing that the literature values for $v$, $T$, $C$ shown in Fig. \ref{fig:Velocity}, Fig. \ref{fig:Temperature} and
Fig. \ref{fig:Concentration} by digitizing Figures of the studies \cite{Arpino2013,Chen2010,Zhao2021} using Engauge Digitizer 11.3.

\begin{figure}[htbp]
	\centering
	\subfigure[$\lambda=0.8$]{
		\includegraphics[width=0.11\textwidth]{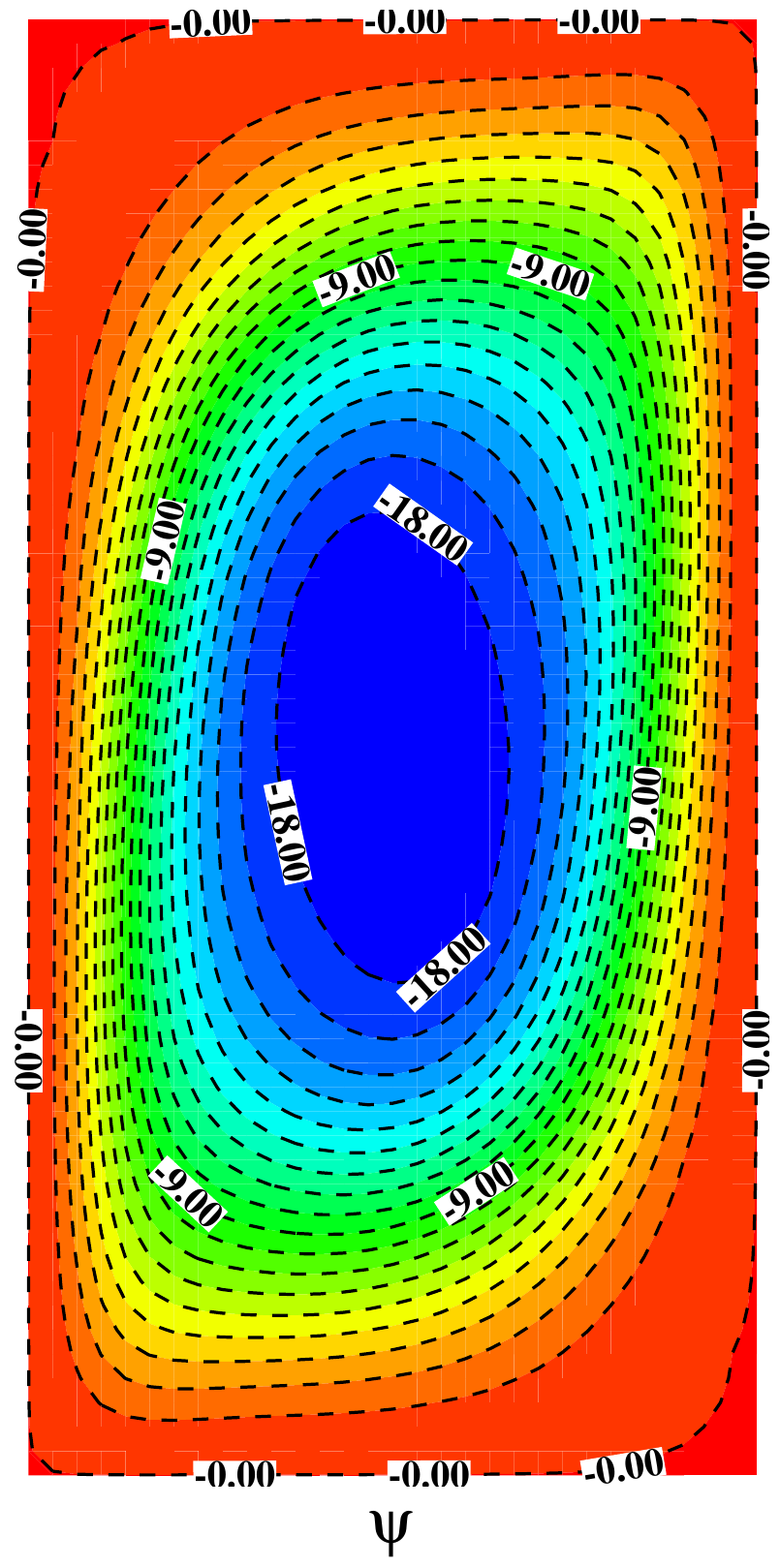}
        \includegraphics[width=0.11\textwidth]{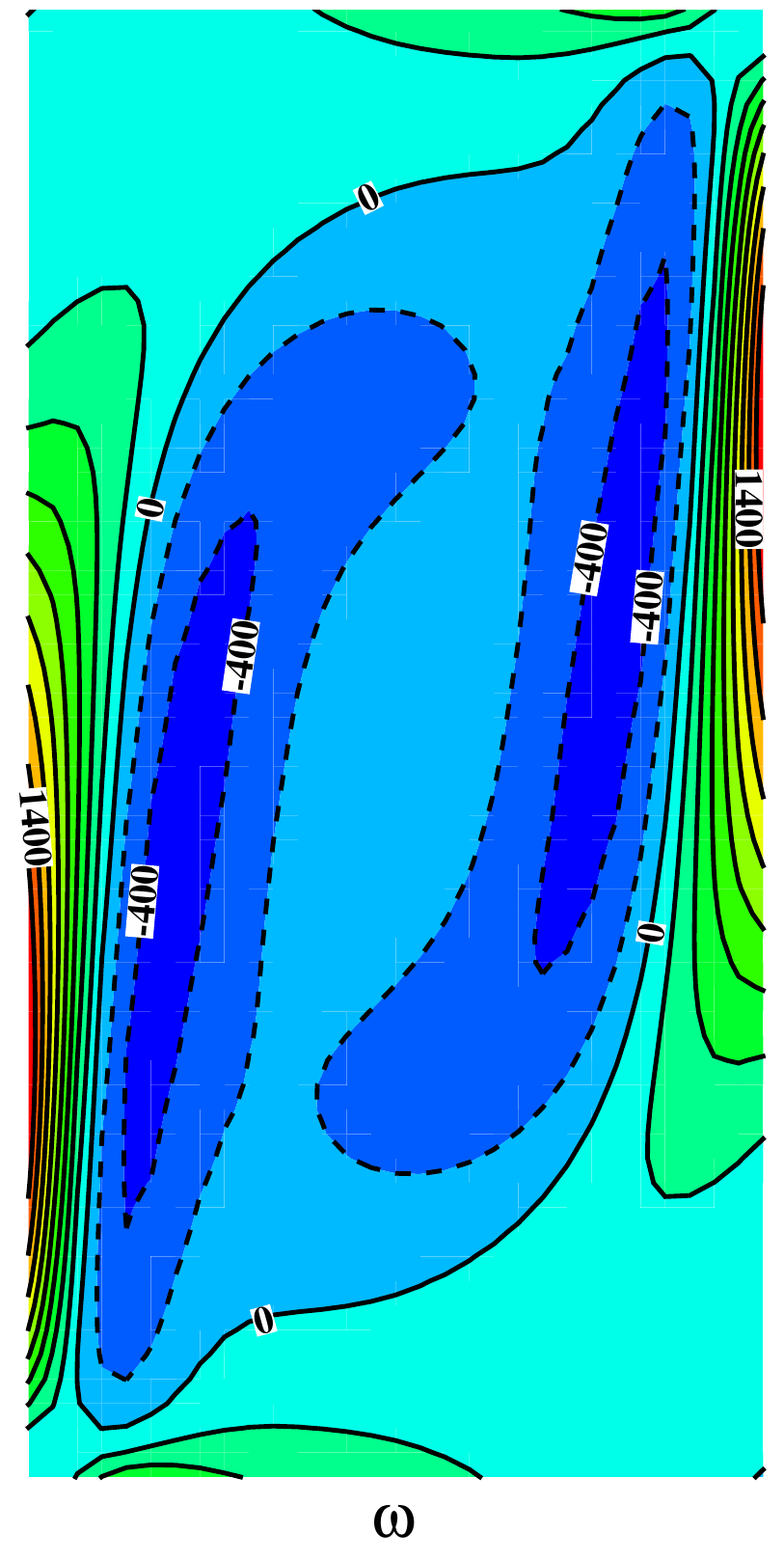}
		\includegraphics[width=0.11\textwidth]{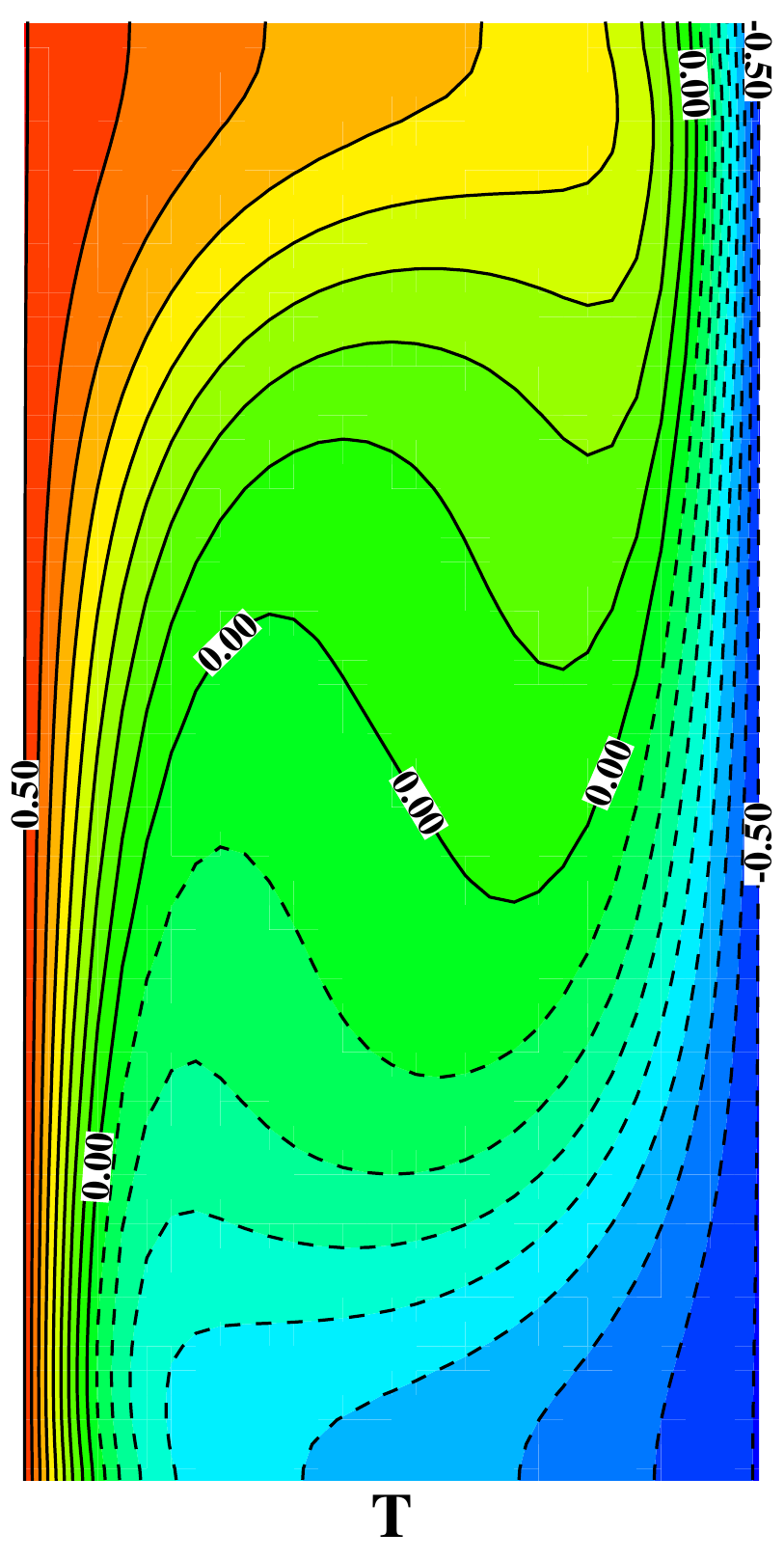}
		\includegraphics[width=0.11\textwidth]{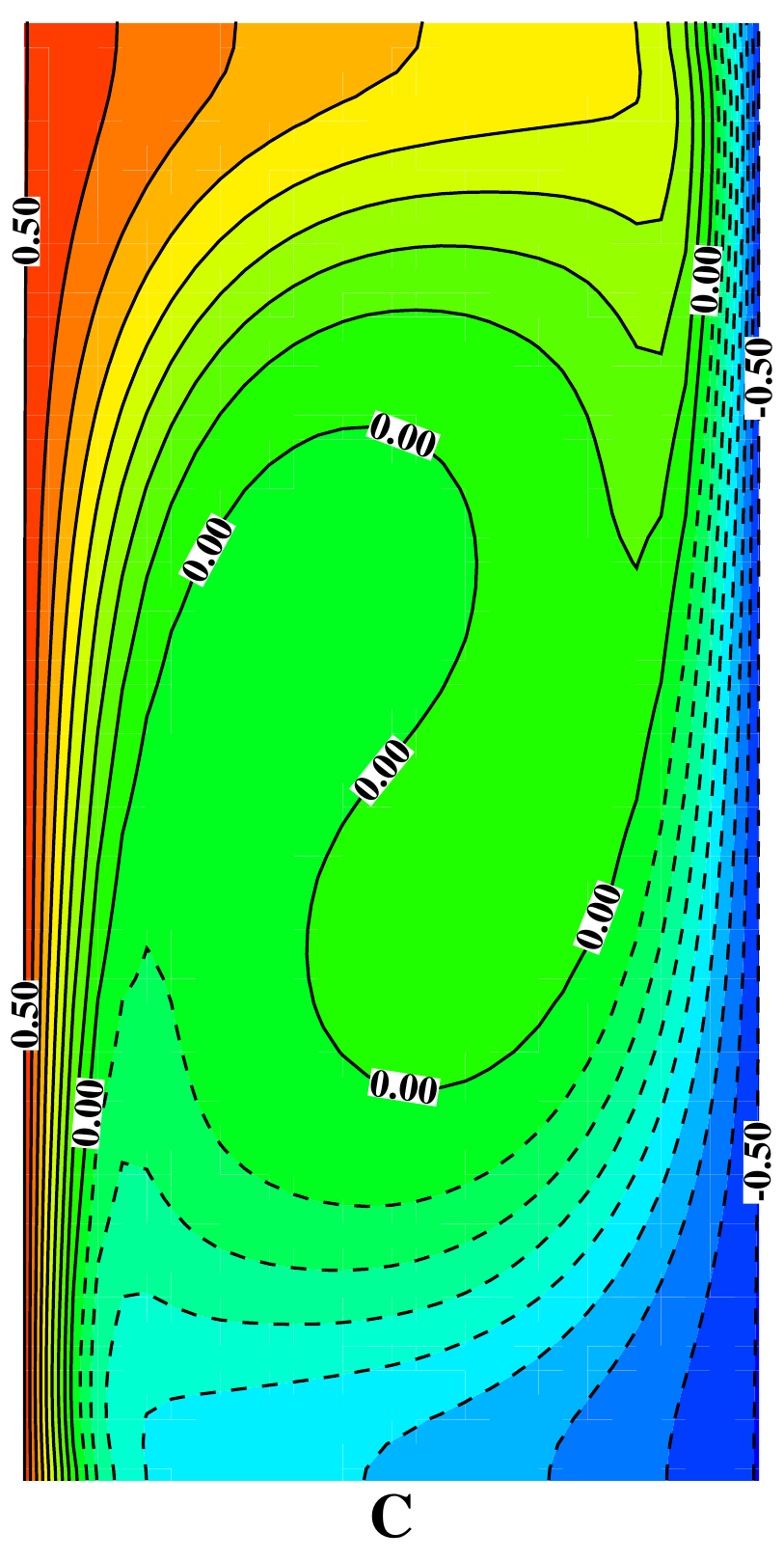}
	}
	\hspace{0.2cm}
	\subfigure[$\lambda=1.3$]{
		\includegraphics[width=0.11\textwidth]{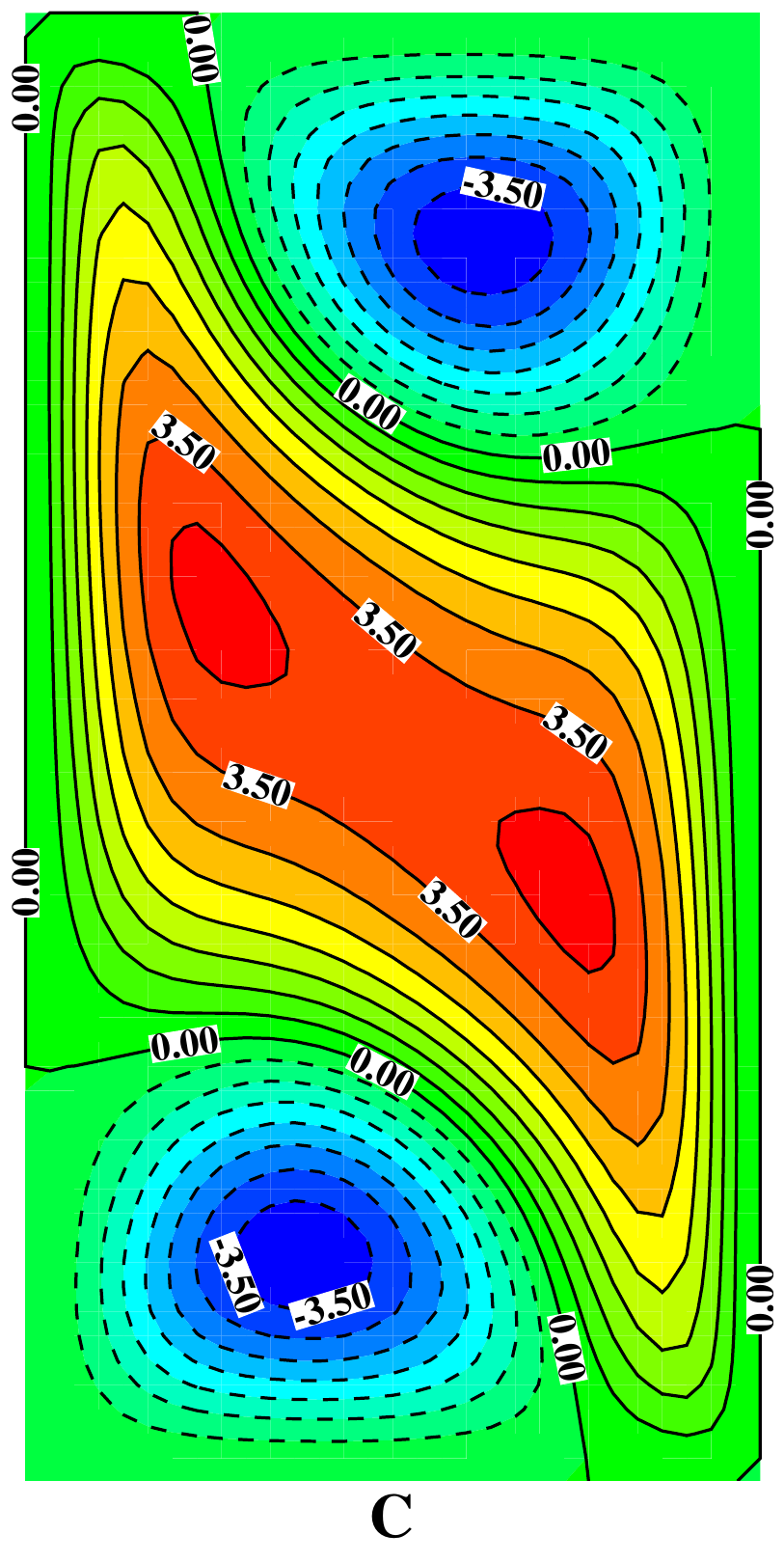}
        \includegraphics[width=0.11\textwidth]{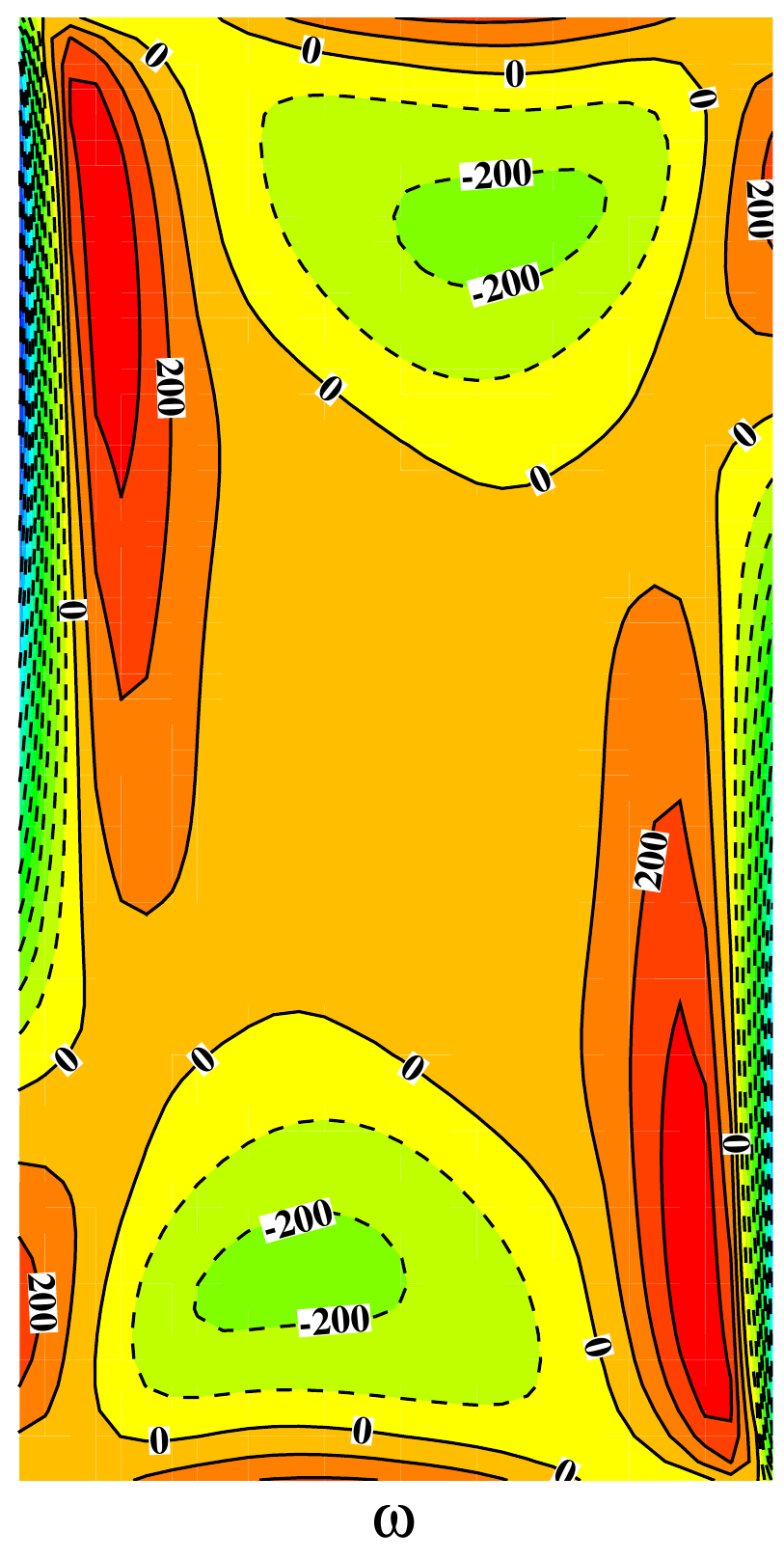}
		\includegraphics[width=0.11\textwidth]{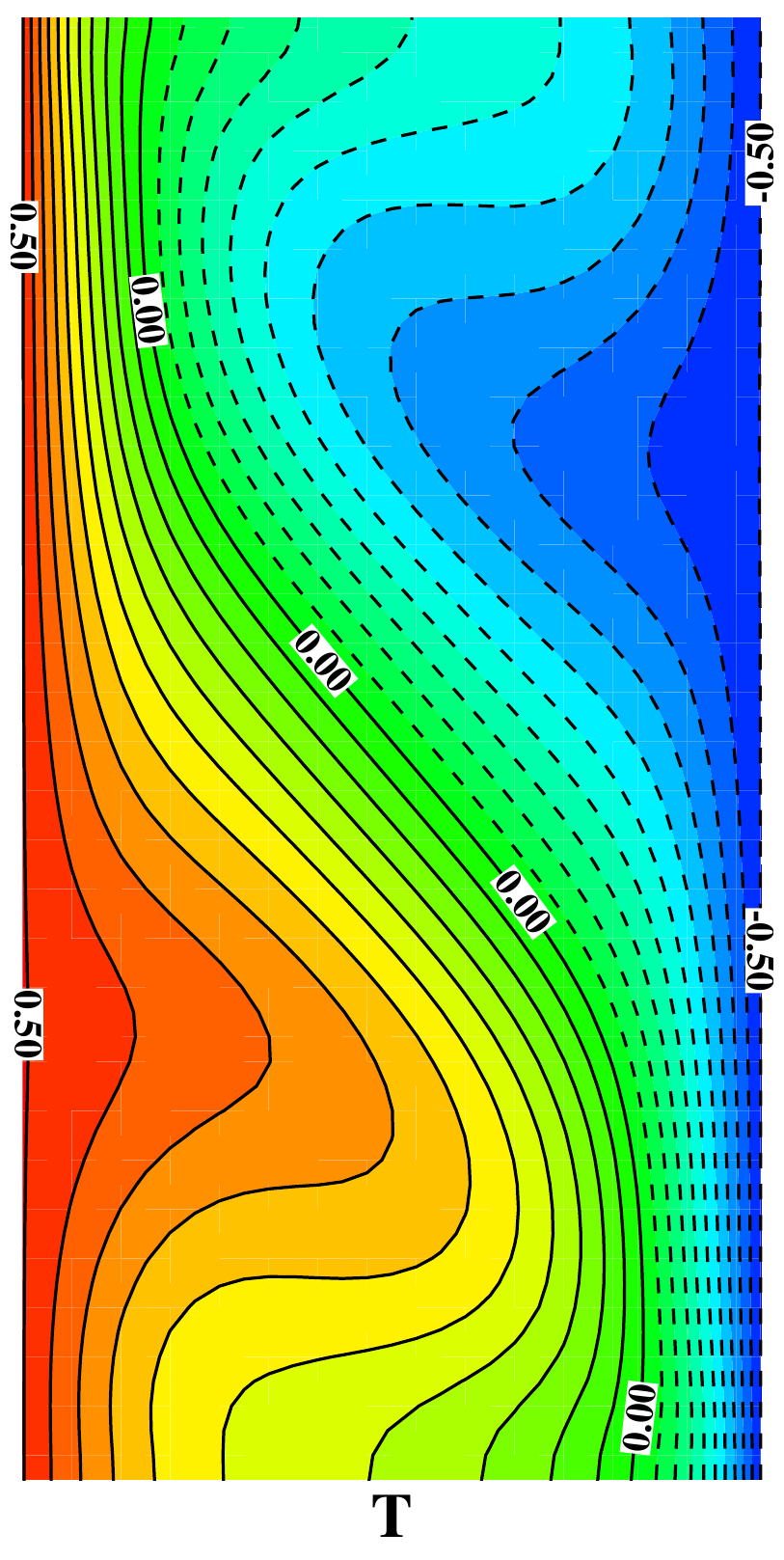}
		\includegraphics[width=0.11\textwidth]{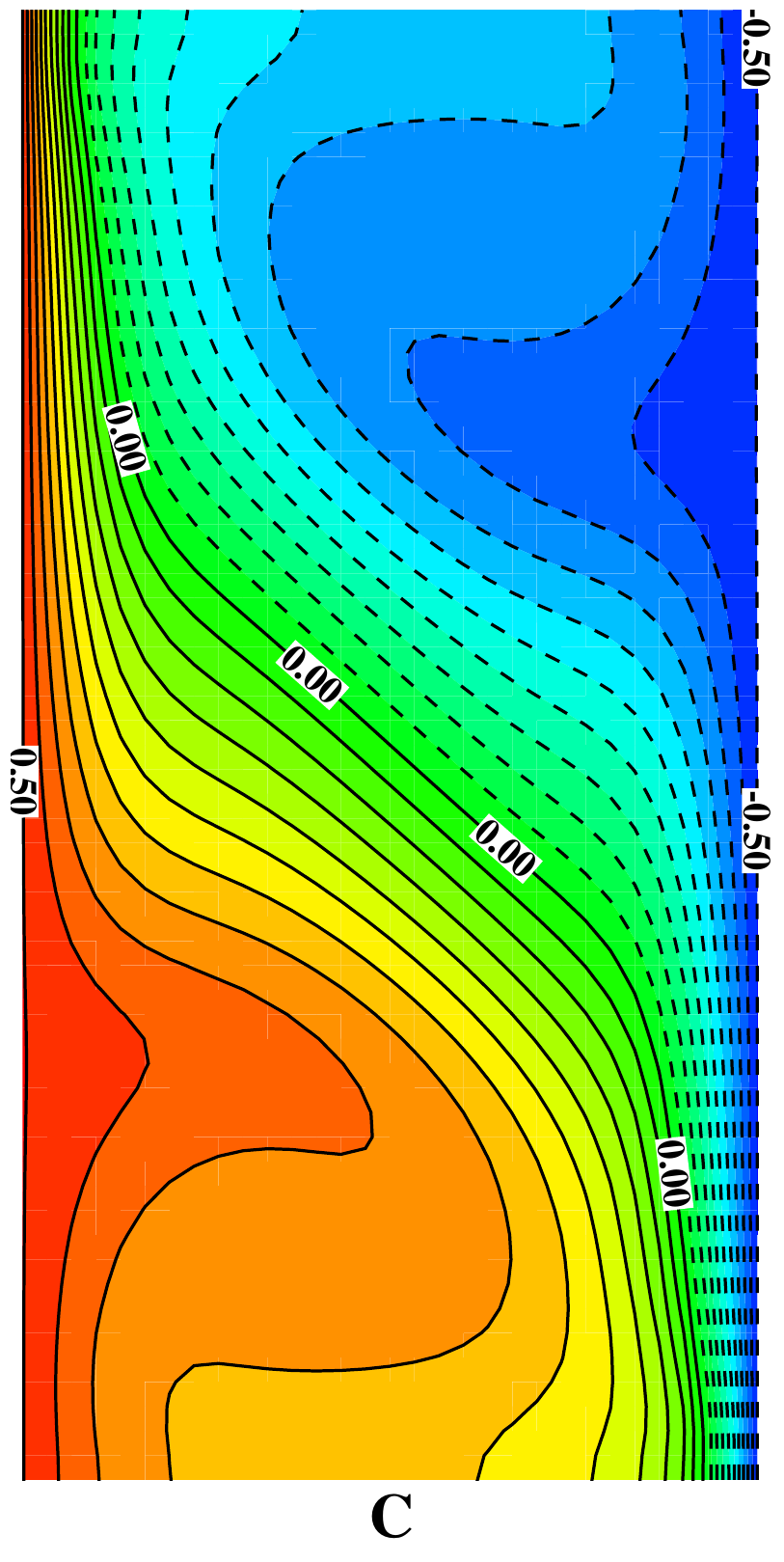}
	}
	\caption{\label{fig:fq0p8andfq1p3} Solutions for $Pr=1$, $Le=2$, $Ra=10^{5}$, $A=2$: (a) $\lambda=0.8$ and (b) $\lambda=1.3$. $\psi$, $\omega$, $T$, and $C$ denote contour lines of the stream function, vorticity, temperature, and concentration, respectively.}
\end{figure}

To have a better understanding the steady flow, The contour line of stream function, vorticity, temperature and concentration are shown in Fig. \ref{fig:fq0p8andfq1p3}.
We can see that a main single vortex occupies the middle part of the cavity and its rotation is clockwise when $\lambda=0.8$.
The isotherms affected by the vertical velocity profile exhibit significant distortion in the core region far from the vertical walls.
Due to the solutal diffusivity being $0.8$ times thermal diffusivity, the concentration contours are also distorted in the core.
When $\lambda=1.3$, there exists a counter-clockwise vortex in the middle part of the cavity and a pair of clockwise vortices exists near the corners of the cavity.
The temperature contours near the contour value $T=0$  are parallel to each other. The concentration contours have the Similar properties.
It can be found that our results are match well with the solutions in the studies.\cite{Nishimura1998}

{\bf Part 3. Grid independence }

To choose a suitable grid size for the simulation of double-diffusive convection, a grid independence study is essential.
Table \ref{table:grid-independence} shows the results of double-diffusive convection in a cavity with $\lambda=1.3$, $Pr=1$, $Le=2$, $Ra=10^{5}$, $A=2$ for different grid meshes.
$u_{\max}$ and $v_{\max}$ are maximum horizontal and vertical velocities, respectively.  $\mid\psi_{\mbox{mid}}\mid$ is the absolute value of steam function at the mid node of the cavity.
\begin{table}[h]
\centering
\caption{\label{table:grid-independence}%
Results of the gride independence study at $Ra=10^{5}$, $Pr=1$, $Le=2$ $\lambda=1.3$ and $A=2$.}
\begin{tabular}{l@{\hspace{4pt}}l@{\hspace{8pt}}l@{\hspace{8pt}}l@{\hspace{8pt}}l@{\hspace{8pt}}
                                l@{\hspace{8pt}}l@{\hspace{8pt}}l@{\hspace{8pt}}l@{\hspace{8pt}}l@{\hspace{8pt}}l@{\hspace{8pt}}}
\hline
$N_{x}\times N_{y}$  &  $Nu_{av}$   &  $\epsilon(\%)$  &  $Sh_{av}$         & $\epsilon(\%)$ &  $u_{\max}$  &  $\epsilon(\%)$
                     & $v_{\max}$   &  $\epsilon(\%)$  &  $\mid\psi_{\mbox{mid}}\mid  $ &  $\epsilon(\%)$  \\
\hline
\textcolor[rgb]{0.00,0.18,0.51}{CHD4}&&&&&&&&&&  \\
$20\times40$	 &1.98331&3.572 &	3.14141	&8.709	&19.28084	&3.149	&22.76448	&4.960	&3.82979	&0.604\\
$30\times60$     &1.91776&0.149 &	2.91170	&0.760	&18.74871	&0.302	&23.88516	&0.281	&3.81604	&0.243\\
$40\times80$	 &1.91178&0.163 &	2.87780	&0.413	&18.70551	&0.071	&24.03571	&0.347	&3.80864	&0.049\\
$50\times100$    &1.91247&0.127 &	2.87861	&0.385	&18.67574	&0.088	&23.95500	&0.010	&3.80668	&0.003\\
$60\times120$	 &1.91340&0.079 &	2.88247	&0.252	&18.69516	&0.016	&23.83735	&0.481	&3.80645	&0.009\\
$100\times200$	 &1.91491&-	    &	2.88974 &-		&18.69220 	&-	    &23.95251 	&-	    &3.80679    &-	\\
\hline
\textcolor[rgb]{0.00,0.25,0.50}{CHD6}&&&&&&&&&& \\
$20\times30$&	1.95245&	1.939&	3.06257&	5.927&	18.77368&	0.435&	22.90105&	4.388&	3.85699&	1.319\\
$30\times60$&	1.91488&	0.023&	2.88796&	0.112&	18.62524&	0.359&	23.86177&	0.377&	3.81488&	0.213\\
$40\times80$&	1.91452&	0.042&	2.88160&	0.332&	18.69857&	0.033&	24.01446&	0.260&	3.80618&	0.015\\
$50\times100$&	1.91506&	0.014&	2.88701&	0.145&	18.68036& 	0.064&	23.94339&	0.037&	3.80540&    0.036\\
$60\times120$&	1.91518&    0.007&	2.88919&    0.070&	18.69932&	0.037&	23.83172&	0.503&	3.80581&	0.025\\
$100\times200$&	1.91532&	    -&  2.89120&		-&  18.69232&		-&  23.95218&		-&  3.80677&	 -   \\
\hline
\end{tabular}
\end{table}

It can be seen that, compared with the results on the grid $100\times200$, the maximum relative errors (such as the $\epsilon$ of $Nu_{av}$ is $\epsilon=\mid\frac{Nu^{N_{x}\times N_{y}}_{av}-Nu^{100\times200}_{av}}{Nu^{100\times200}_{av}}\mid\times100\%$ of $Nu_{av}$, $Sh_{av}$, $u_{\max}$, $v_{\max}$ and $\mid\psi_{\mbox{mid}}\mid$ on the grid $40\times80$ change less than $1.00\%$.
In addition, we observe that the numerical results of the CHD4 and CHD6 schemes are almost identical from the Table \ref{table:grid-independence}.
Taking into account the CPU time, the grid mesh $40\times80$ using the CHD4 scheme is chosen for the simulation of double-diffusive convection when the flow is steady flow. It is worth noting that a finer grid mesh using the CHD4 scheme is required to obtain better accuracy when the flow is oscillatory flow.
	
\subsubsection{\label{sec3.2.2} Effect of $Ra$}
We consider these cases with buoyancy ratios $\lambda=0.2, 0.8, 1.3, 1.8$ in a rectangular cavity with $A=2$, $Pr=1$, and $Le=2$. The range of Rayleigh number is $Ra=10^{4}$--$10^{7}$. Table \ref{table:RaFlow} presents the flow state of double-diffusive convection for each operating condition, where S, MP and CH represent the steady, periodic and chaotic flows, respectively.
All numerical simulations are performed with the grid size of $40\times80$, and finer grid size may be required for periodic or chaotic flows.
\begin{table}[htbp]
\centering
\caption{\label{table:RaFlow}%
Flow pattern for different Rayleigh numbers at $A=2$, $Pr=1$, and $Le=2$.}
\begin{tabular}{ccccc}
\hline
$Ra$         &  $\lambda=0.2$  &  $\lambda=0.8$  &  $\lambda=1.3$  &  $\lambda=1.8$  \\
\hline
$10^{4}$  & S & S  & S  & S \\
$10^{5}$  & S & S  & S  & S \\
$10^{6}$  & S & MP & S  & S \\
$10^{7}$  & S & CH & CH & CH \\
\hline
\end{tabular}
\end{table}

\begin{figure}[h]
	\centering
	\subfigure[$Nu_{av}$ vs $Ra$]{
		\includegraphics[width=0.35\textwidth]{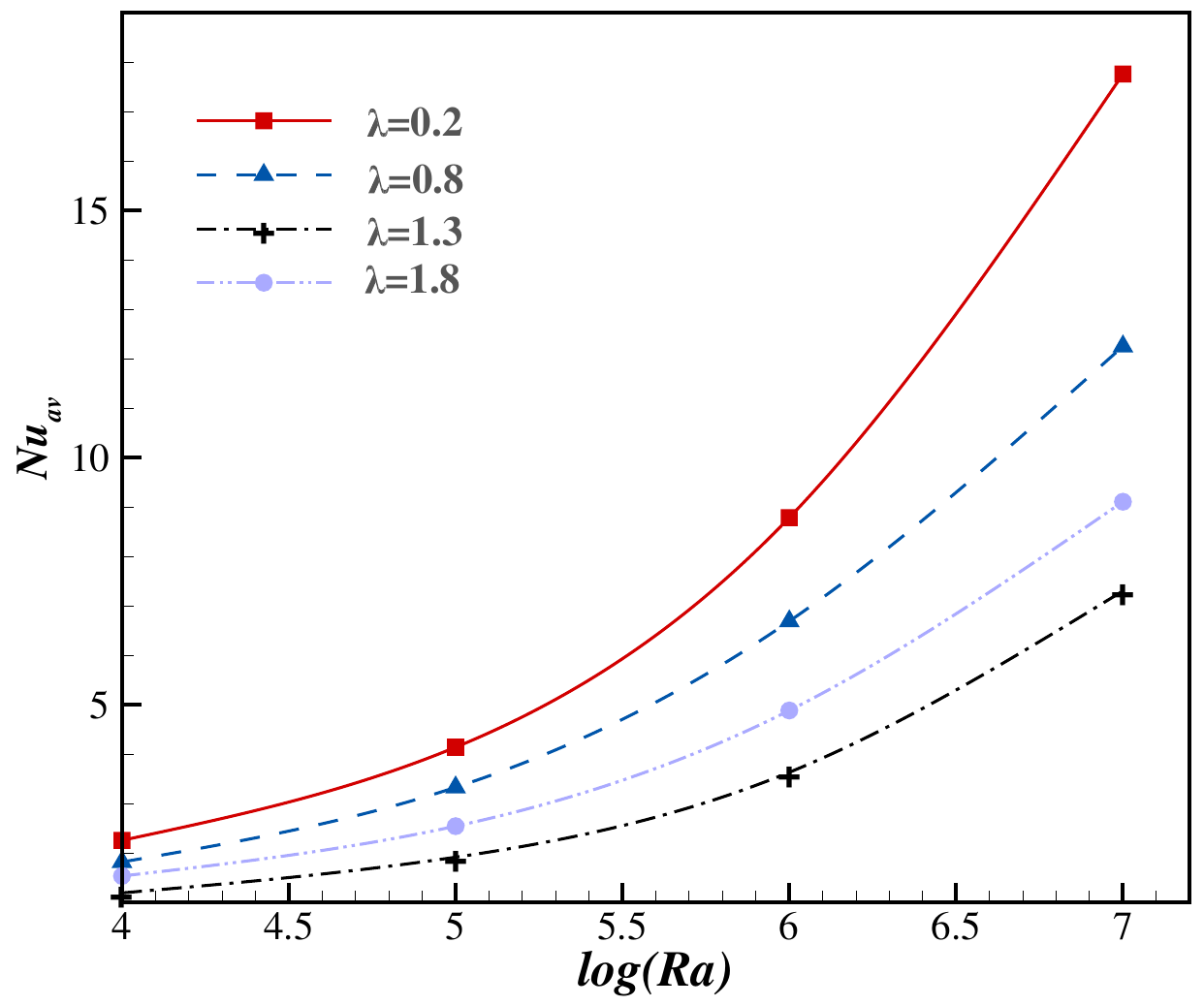}
	}
	\hspace{0.2cm}
	\subfigure[$Sh_{av}$ vs $Ra$]{
		\includegraphics[width=0.35\textwidth]{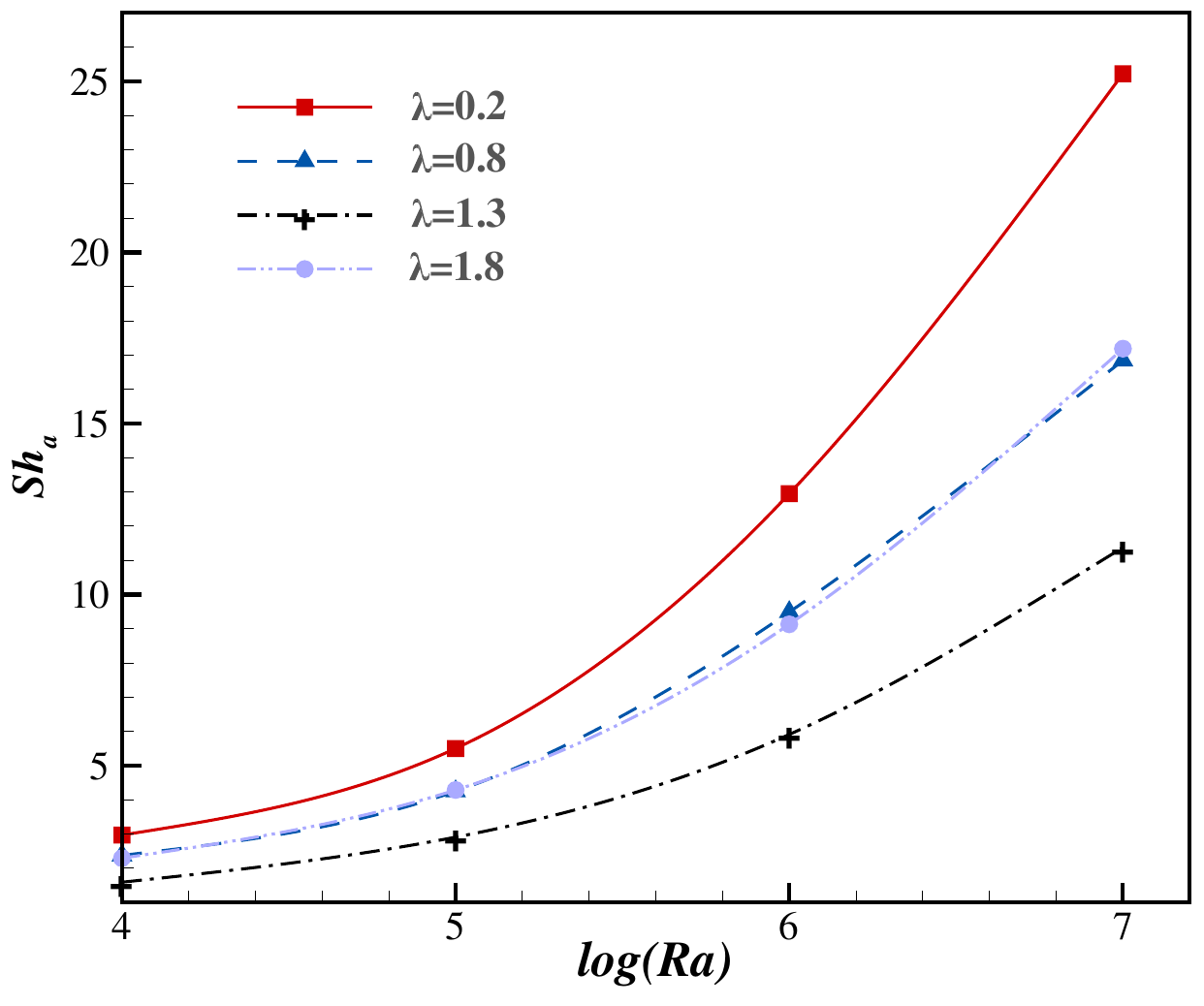}

	}
	\caption{\label{fig:NuShVSRa}Dependence of (a)$Nu_{av}$ and (b)$Sh_{av}$ on Rayleigh number for different $\lambda$.}
\end{figure}

Fig. \ref{fig:NuShVSRa} shows the variation of the heat and mass transfer ($Nu_{av}$ and $Sh_{av}$) with increasing Rayleigh number.
The results show that both $Nu_{av}$ and $Sh_{av}$ increase exponentially with increasing Rayleigh number. In addition, the value of $Sh_{av}$  is always larger than that of $Nu_{av}$, since the heat diffusivity is twice as much as the concentration diffusivity.
\begin{figure}[h]
	\centering
	\subfigure[]{
		\includegraphics[width=0.3\textwidth]{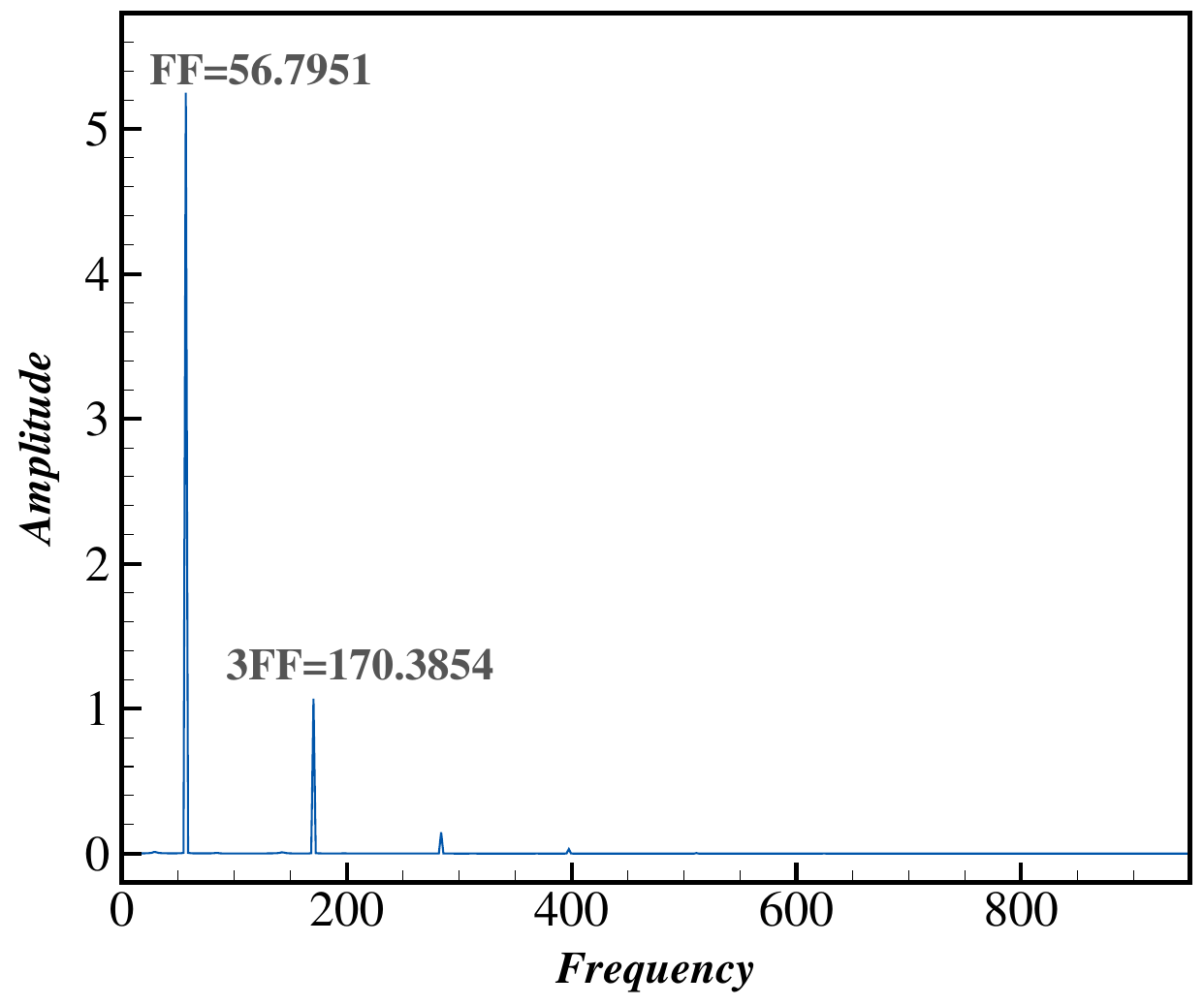}}
	\hspace{0.2cm}
	\subfigure[]{
		\includegraphics[width=0.3\textwidth]{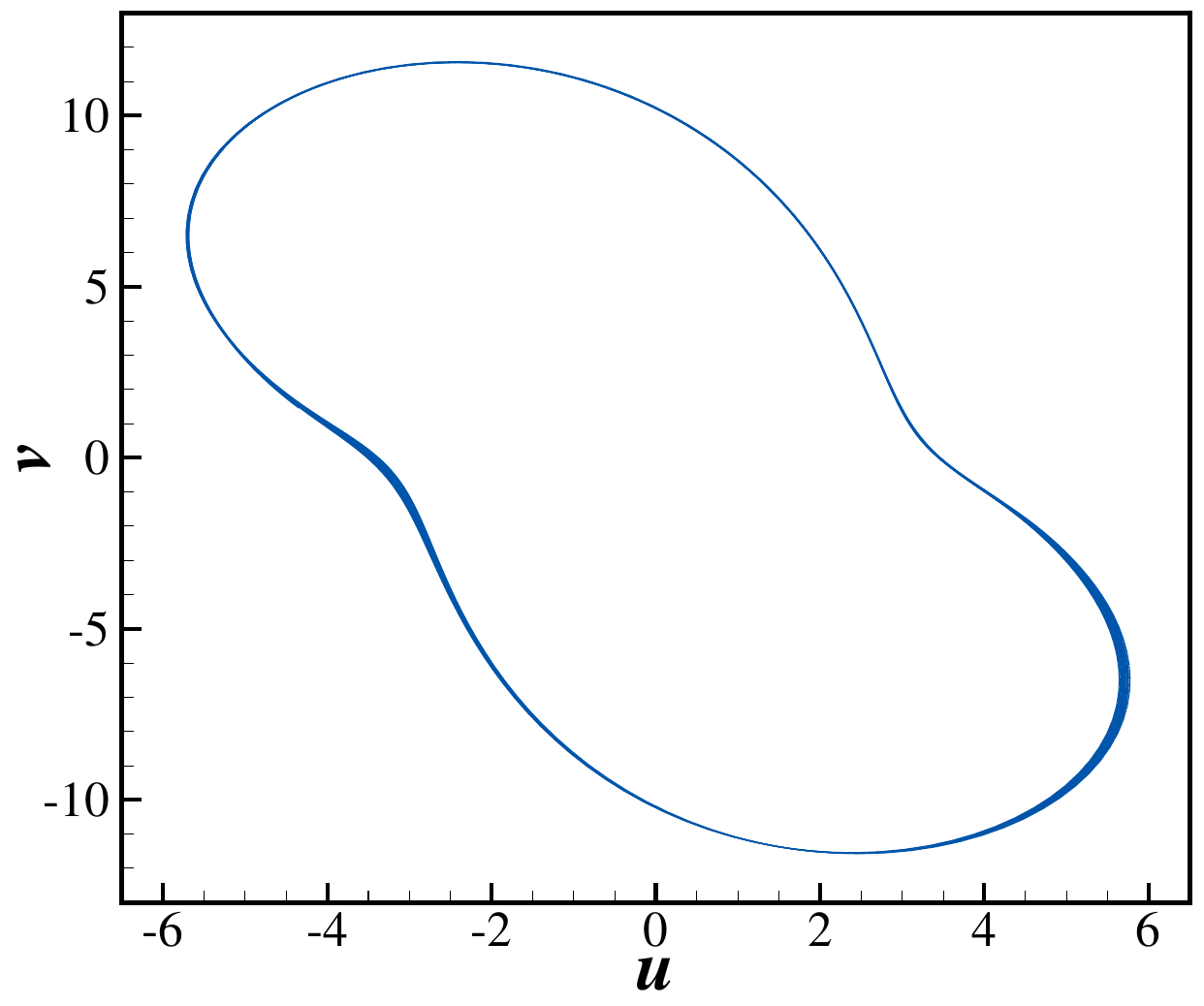}}
	\hspace{0.2cm}
	\subfigure[]{
		\includegraphics[width=0.3\textwidth]{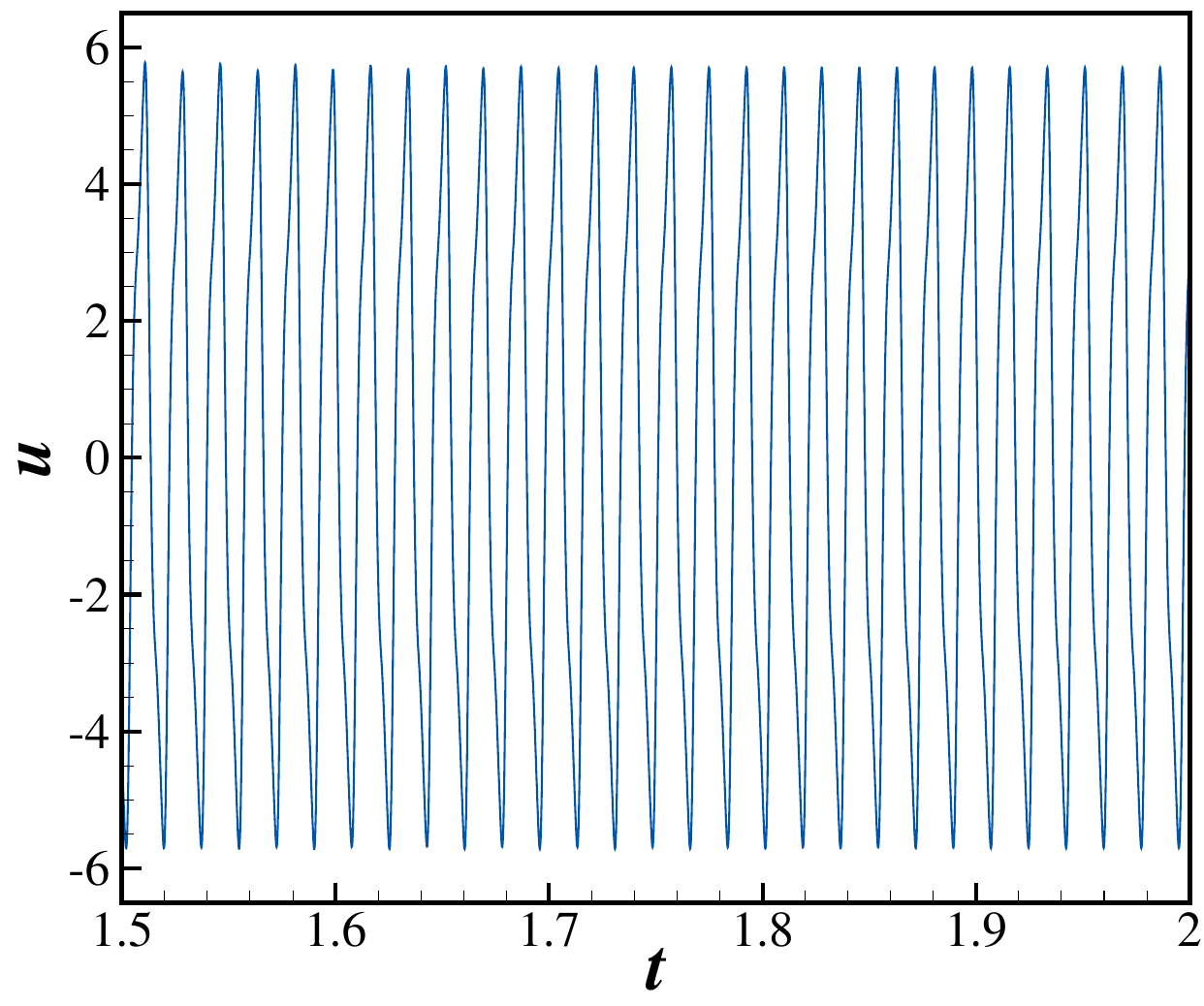}}
	\caption{\label{fig:fq0p8Ra1e6}Properties of the flow at $\lambda=0.8$, $Ra=10^{6}$(MP): (a) Fourier frequency spectrum of the $u$-velocity, (b) phase-space trajectories of $u-v$, (c) time trace of the $u$-velocity.}
\end{figure}
\begin{figure}[htbp]
	\centering
	\subfigure[]{
		\includegraphics[width=0.11\textwidth]{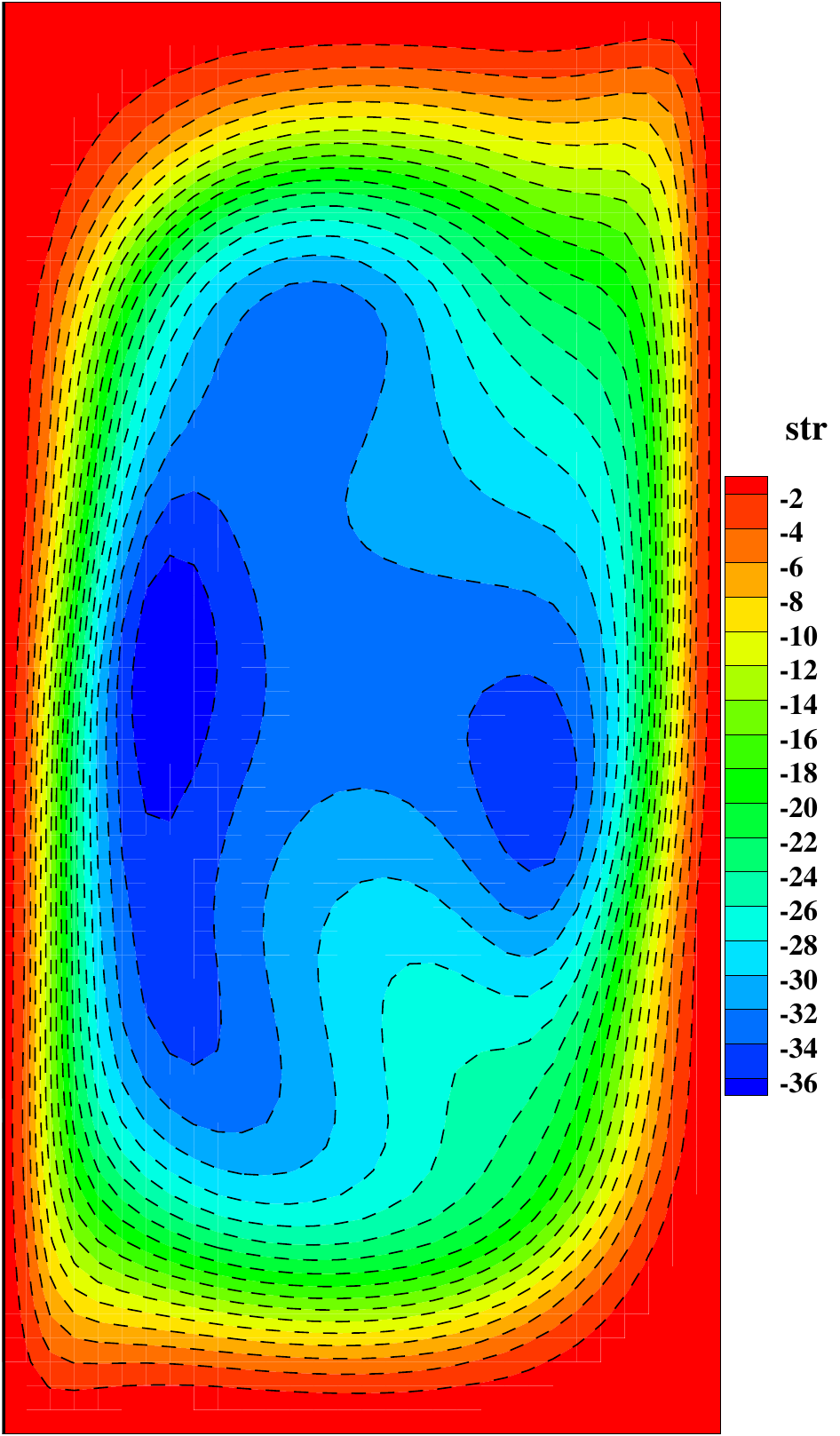}}
	\subfigure[]{
        \includegraphics[width=0.11\textwidth]{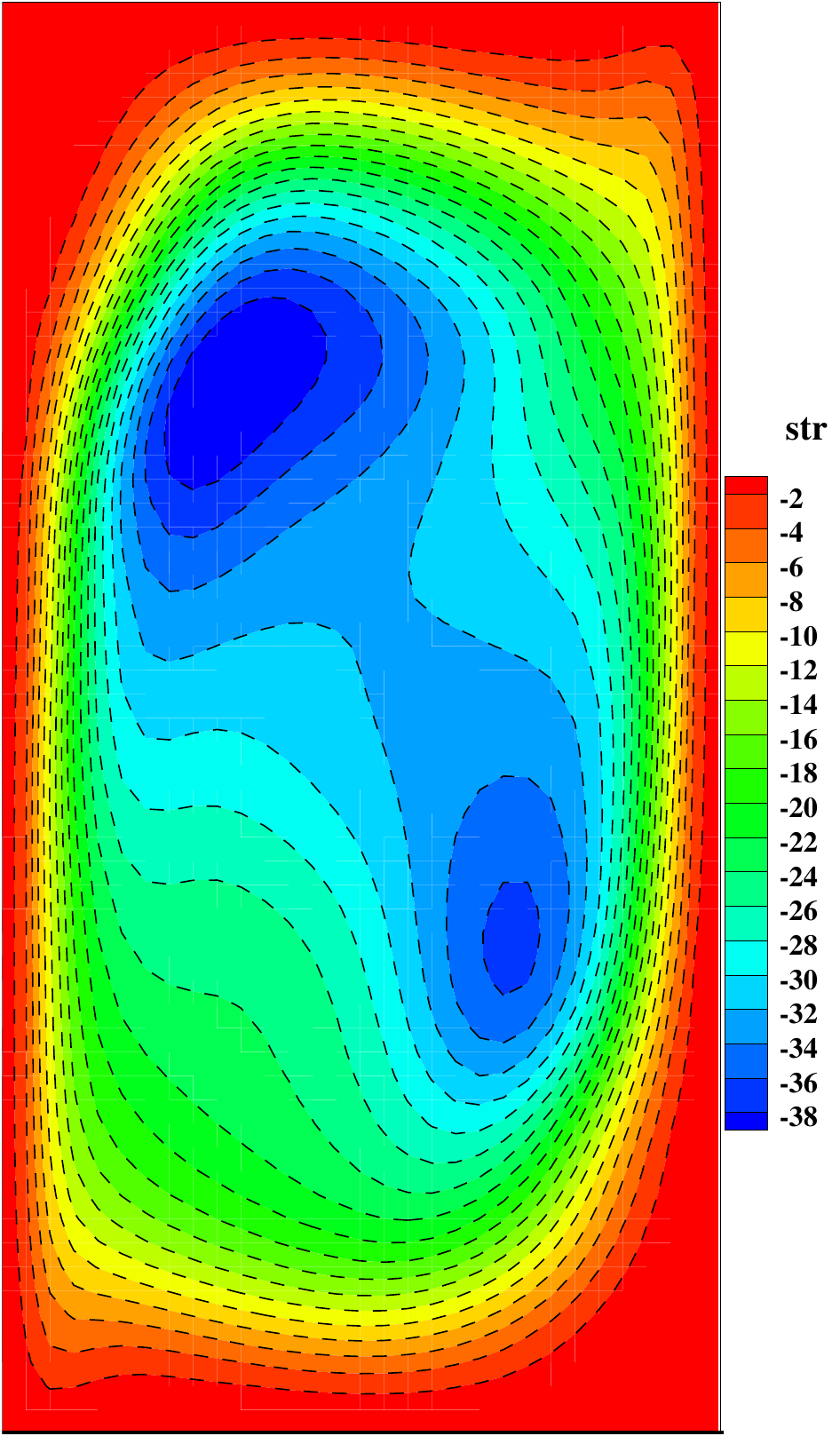}}
	\subfigure[]{
		\includegraphics[width=0.11\textwidth]{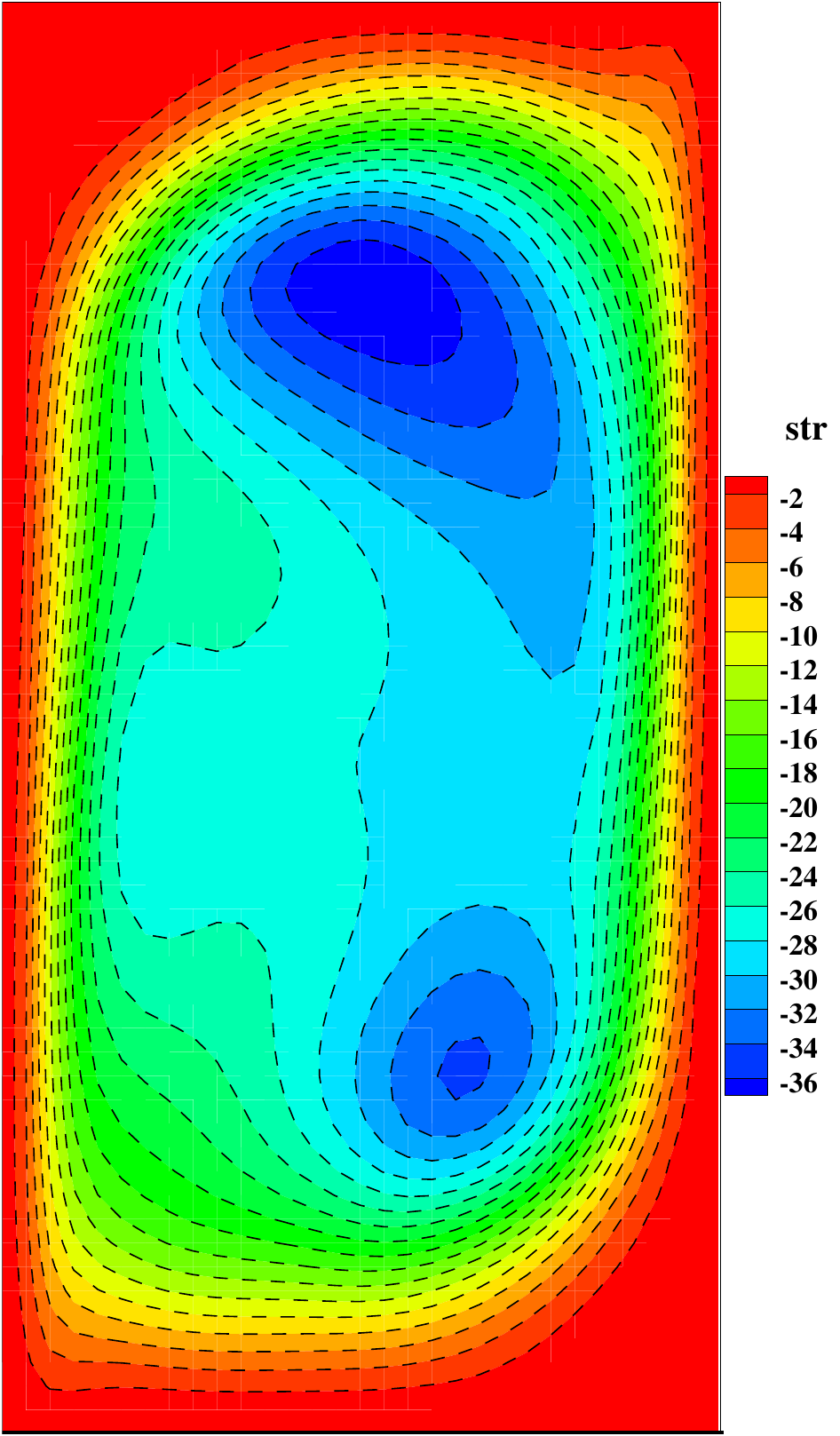}}
	\subfigure[]{
		\includegraphics[width=0.11\textwidth]{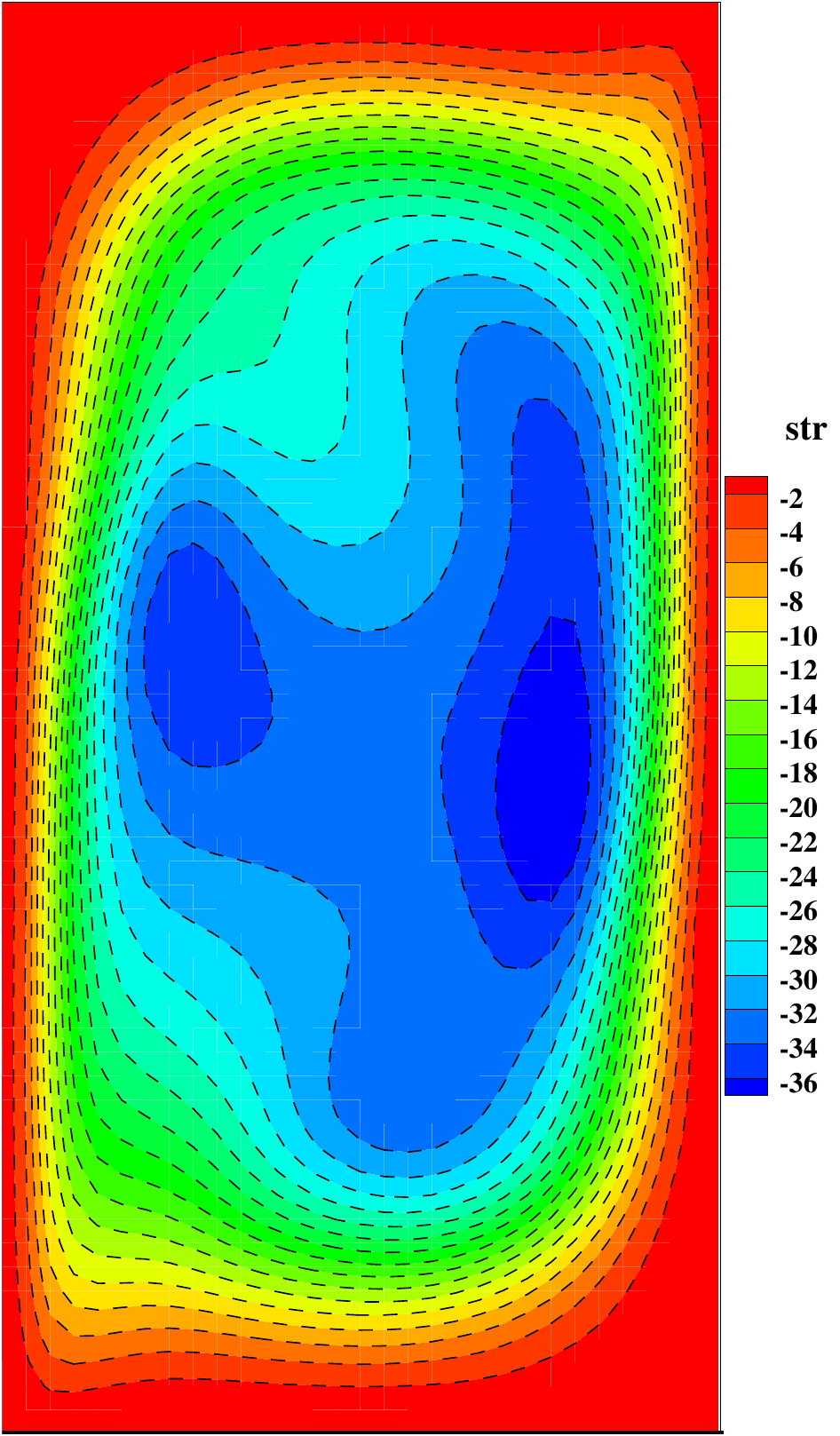}}
    \subfigure[]{
		\includegraphics[width=0.11\textwidth]{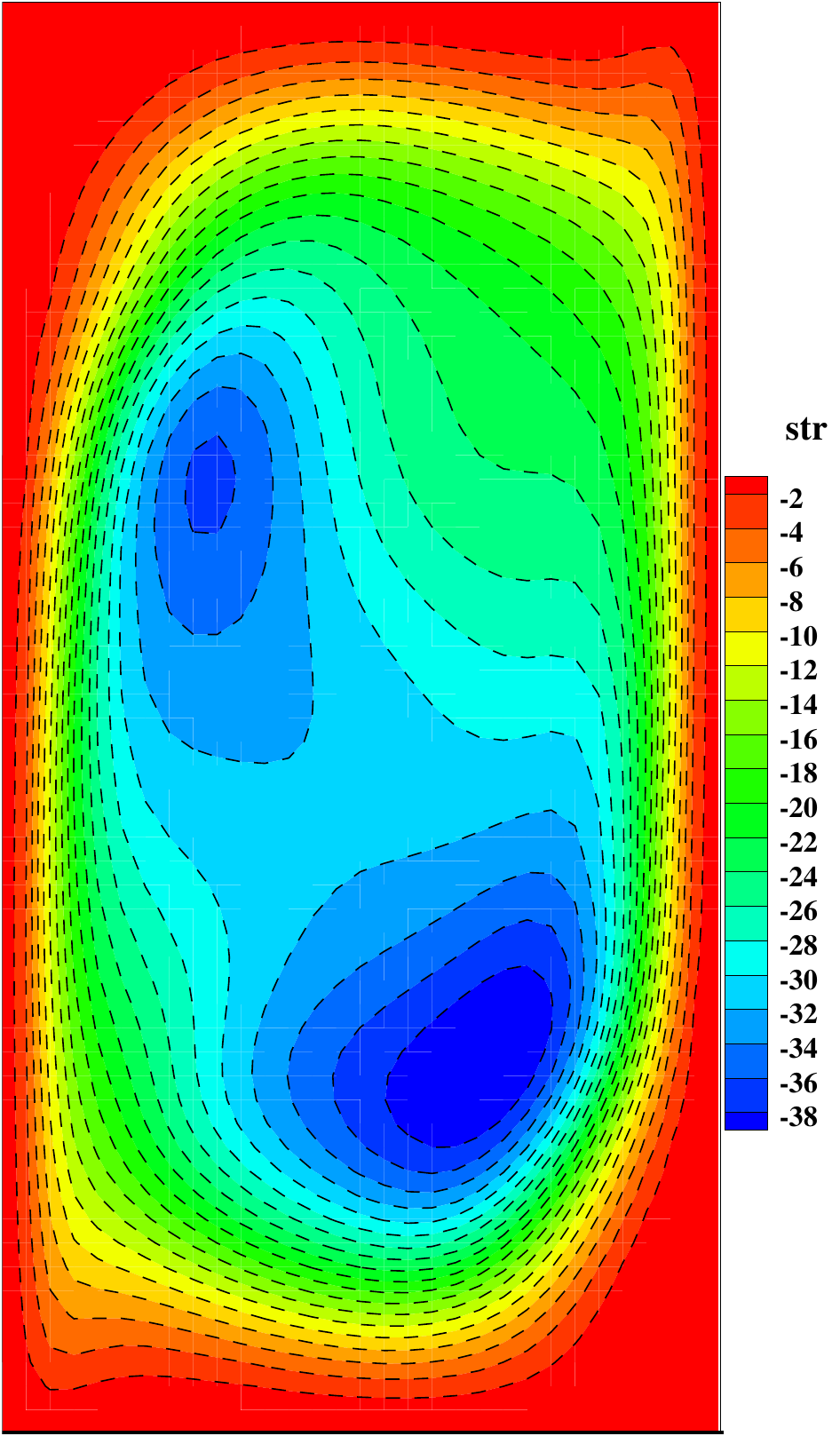}}
	\subfigure[]{
        \includegraphics[width=0.11\textwidth]{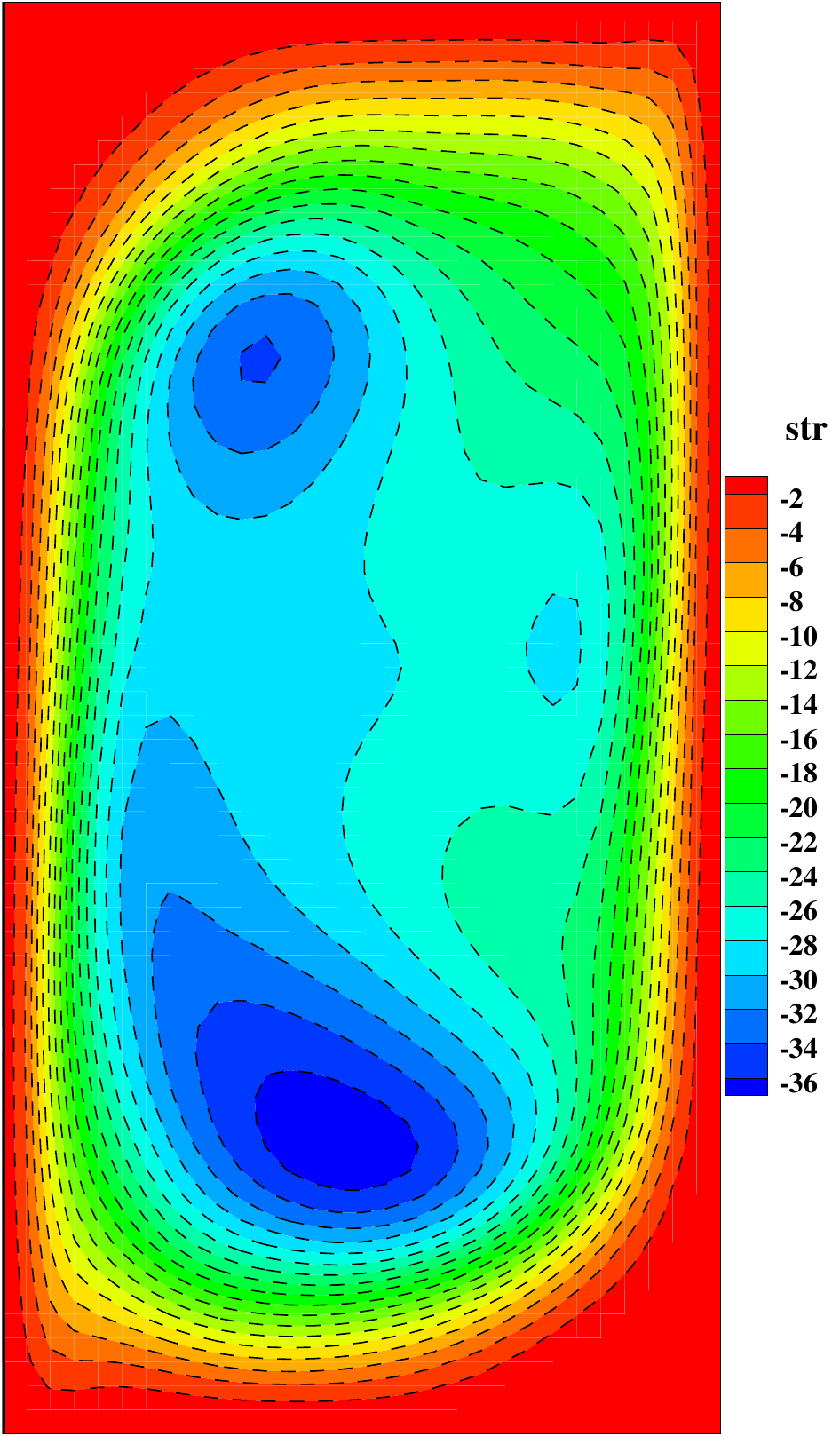}}
	\subfigure[]{
		\includegraphics[width=0.11\textwidth]{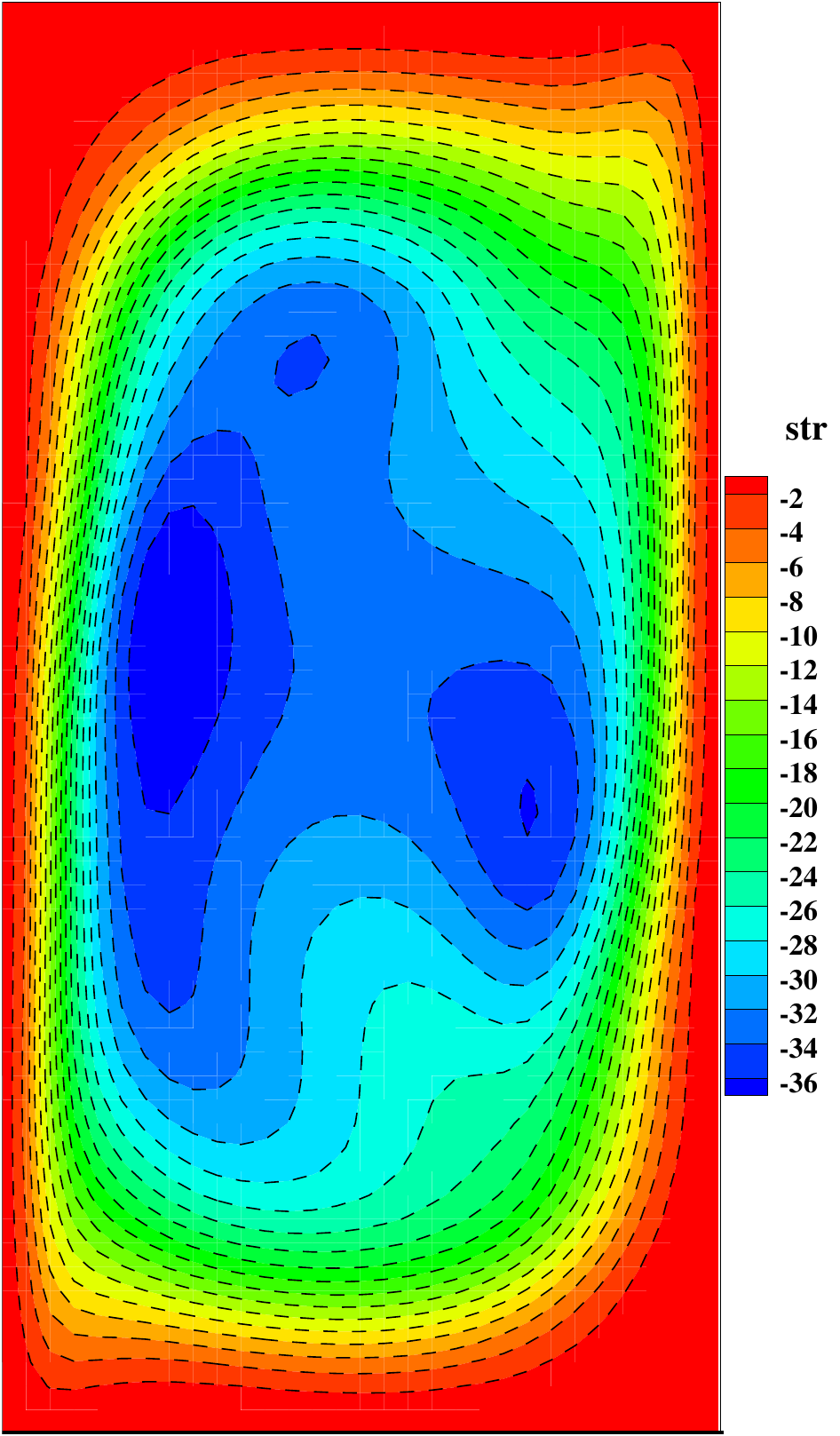}}
	\subfigure[]{
		\includegraphics[width=0.11\textwidth]{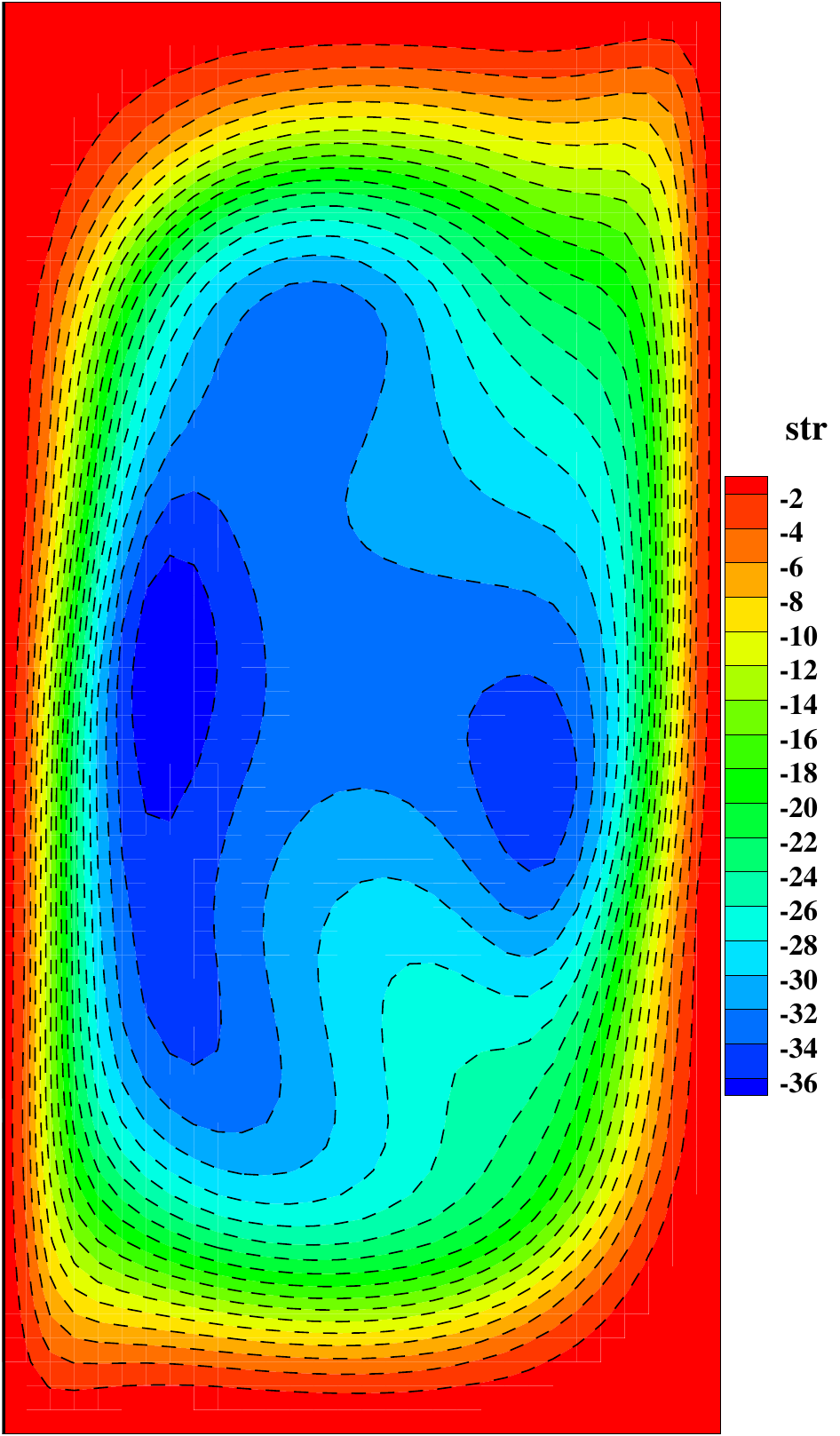}}
	\caption{\label{fig:fq0p8Ra1e6Str} Time history of streamline contours for $\lambda=0.8$, $Ra=10^{6}$. }
\end{figure}

From the Table \ref{table:RaFlow}, we can seen that, under the buoyancy ratio that we consider, the flow changes from the steady flow to a periodic or chaotic flow with increasing Rayleigh number except for $\lambda=0.2$.
When $Ra=10^{6}$, $\lambda=0.8$, the flow is periodic flow. Fig. \ref{fig:fq0p8Ra1e6} shows the time traces of $u$, Fourier frequency spectrum, and phase-space trajectories of $u-v$ at the monitor point $(0.5,\frac{A}{2})$.
The Fourier frequency spectrum of $u$ in Fig. \ref{fig:fq0p8Ra1e6}(a) reveals that the fundamental frequency (FF=56.7951) appears in the flow field initially, interacts with the 3FF wave. Compared to fundamental frequency, the 3FF wave has almost no impact on the flow. Therefor, the time-evolution of $u$ in Fig. \ref{fig:fq0p8Ra1e6}(c) exhibits sinusoidal behavior. The phase trajectory in Fig. \ref{fig:fq0p8Ra1e6}(b) is a closed curve. This implies that the flow is periodic. Fig. \ref{fig:fq0p8Ra1e6Str} presents a complete periodic motion cycle and the minimum period is $0.0176$ at this case.

\begin{figure}[htbp]
	\centering
	\subfigure[]{
		\includegraphics[width=0.3\textwidth]{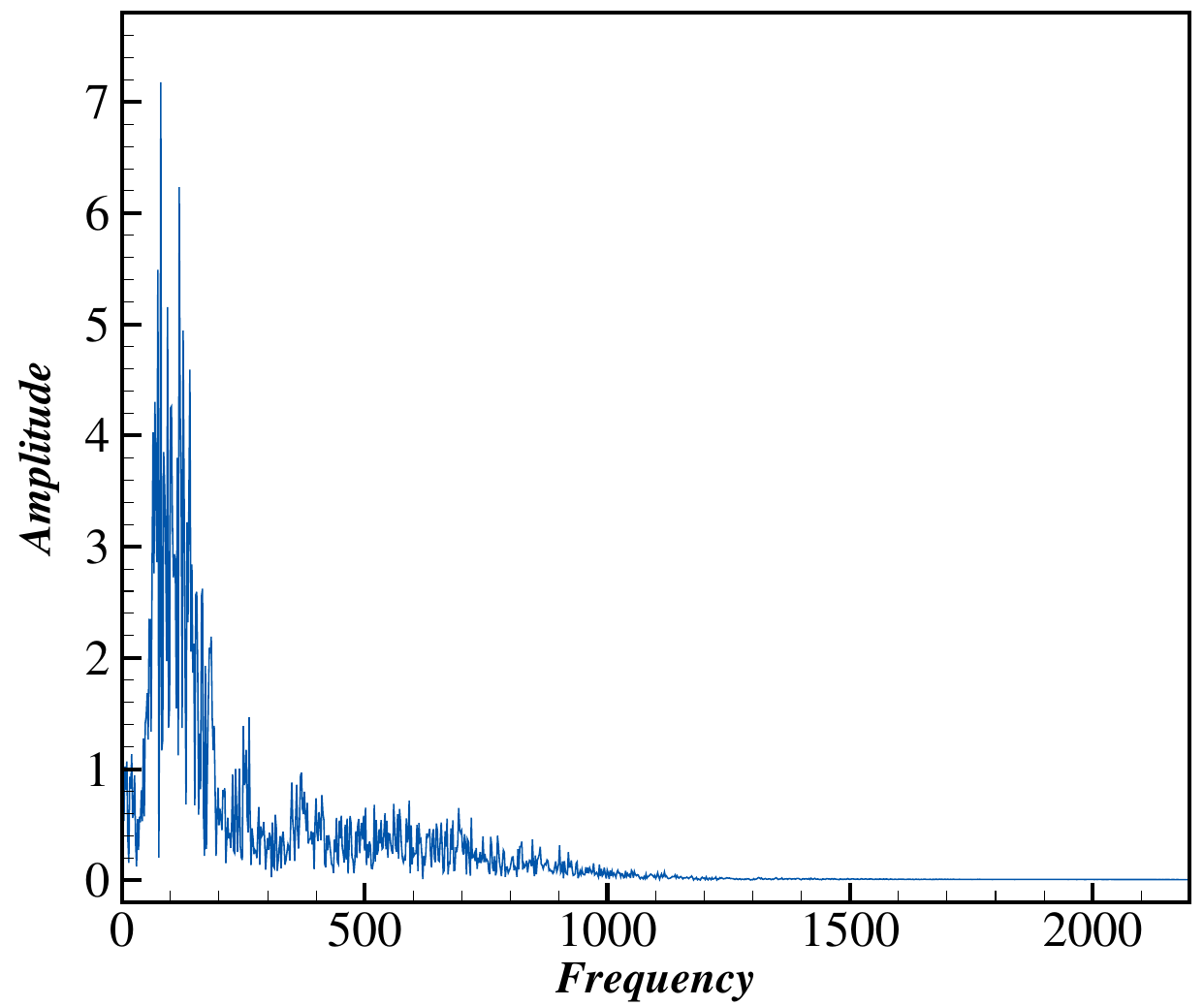}}
	\hspace{0.2cm}
	\subfigure[]{
		\includegraphics[width=0.3\textwidth]{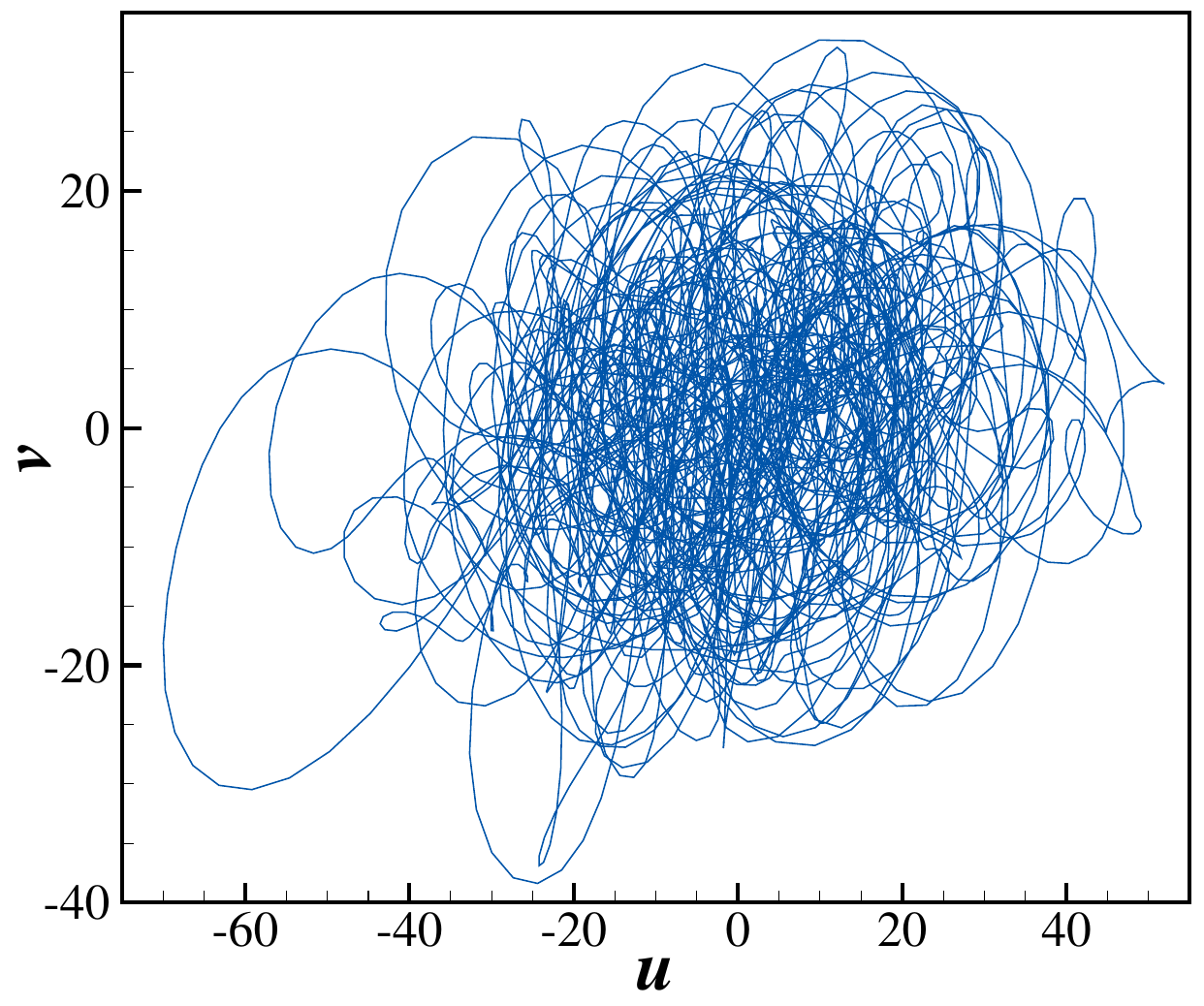}}
	\hspace{0.2cm}
	\subfigure[]{
		\includegraphics[width=0.3\textwidth]{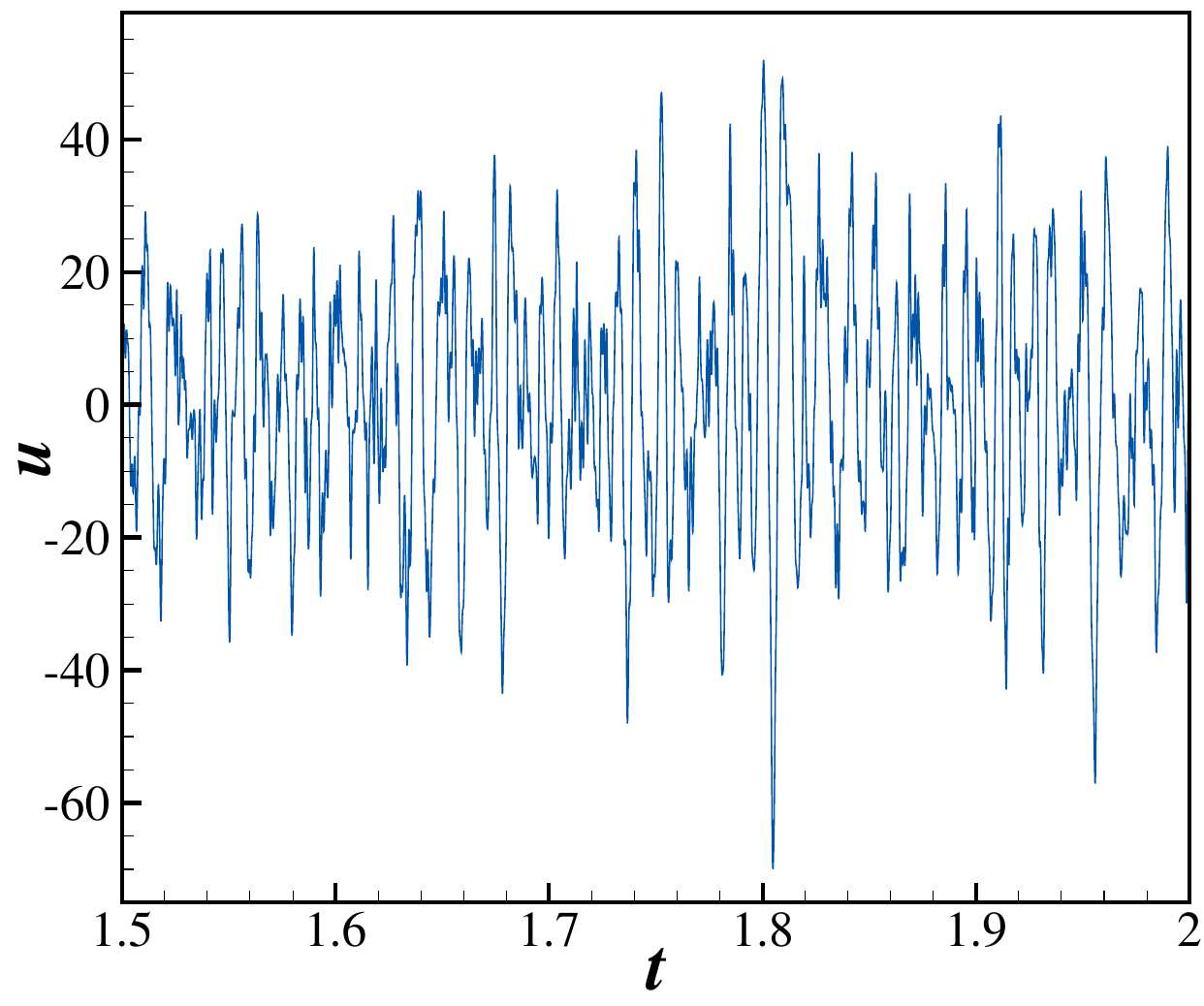}}
	\caption{\label{fig:fq0p8Ra1e7}Properties of the flow at $\lambda=0.8$, $Ra=10^{7}$(CH): (a) Fourier frequency spectrum of the $u$-velocity, (b) phase-space trajectories of $u-v$, (c) time trace of the $u$-velocity.}
\end{figure}
When $Ra=10^{7}$, under the buoyancy ratio that we consider, the flow is chaotic except $\lambda=0.2$. When the flow is chaotic, there is no obvious fundamental frequency, and a large number of sub-harmonic waves are excited in the flow field as shown in Fig. \ref{fig:fq0p8Ra1e7}.

\subsubsection{\label{sec3.2.3} Effect of $A$}
\begin{table}[h]
\centering
\caption{\label{table:AFlow}%
Flow pattern for different aspect ratios at $Ra=10^{5}$, $Pr=1$, and $Le=2$.}
\begin{tabular}{ccccccccc}
\hline
$A$            &$1.0$  &  $2.0$  &  $3.0$  &  $4.0$ &  $5.0$  &  $6.0$  &  $7.0$  &  $8.0$  \\
\hline
$\lambda=0.2$  & S     & S       & S       & S      & S       & S        & S       & S      \\
$\lambda=0.8$  & S     & S       & S       & S      & S       & S        & S       & S      \\
$\lambda=1.3$  & S     & S       & S       & S      & S       & S        & S       & S      \\
$\lambda=1.8$  & S     & S       & S       & S      & S       & S        & S       & S      \\
\hline
\end{tabular}
\end{table}
We also investigate the dependence of the heat and mass transfer on aspect ratio for $Ra=10^{5}$, $Le=2$, $Pr=1$ and $\lambda=0.2, 0.8, 1.3, 1.8$.
The range of Rayleigh number is $A=1.0$--$8.0$. Table \ref{table:AFlow} presents the flow state of double-diffusive convection for each operating condition. We can see that the flow of double-diffusive convection is always steady for different aspect ratios.

\begin{figure}[h]
	\centering
	\subfigure[$Nu_{av}$ vs $A$]{
		\includegraphics[width=0.35\textwidth]{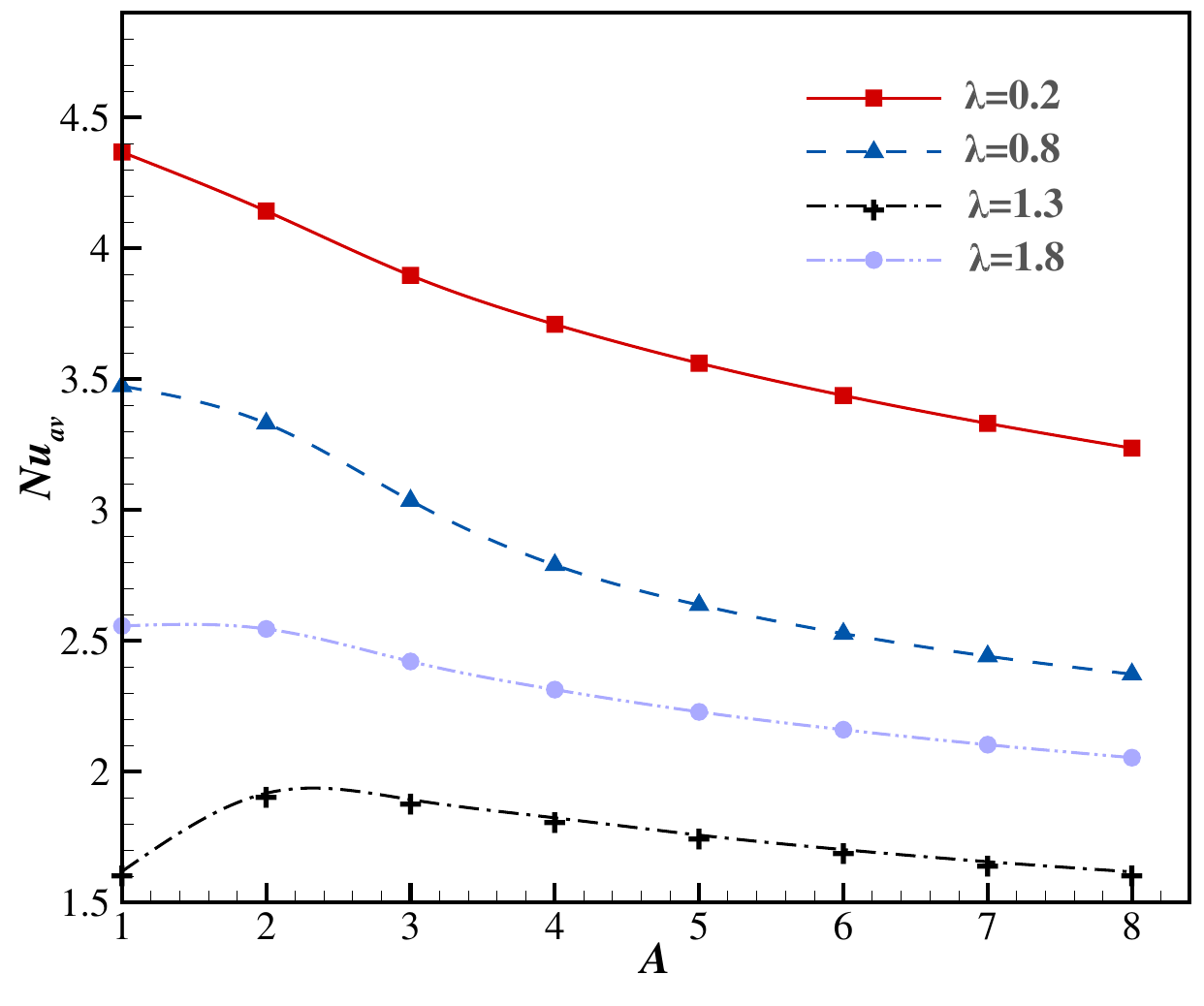}
	}
	\hspace{0.2cm}
	\subfigure[$Sh_{av}$ vs $A$]{
		\includegraphics[width=0.35\textwidth]{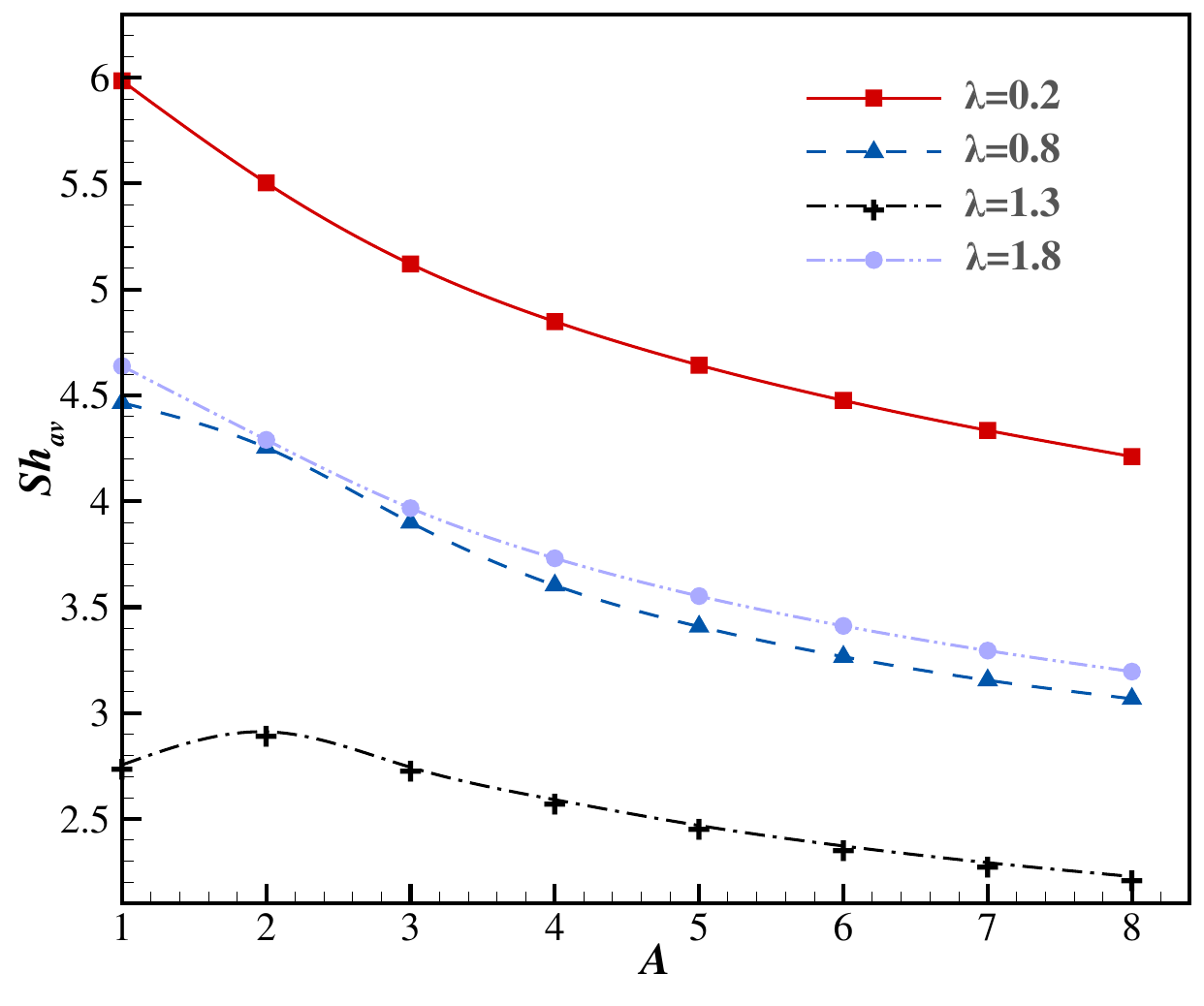}

	}
	\caption{\label{fig:NuShVSA}Dependence of (a)$Nu_{av}$ and (b)$Sh_{av}$ on aspect ratio for different $\lambda$.}
\end{figure}
Fig. \ref{fig:NuShVSA} shows the influence of the aspect ratio $A$ on the heat and mass transfer for different buoyancy ratio.
As the aspect ratio increases, both the $Nu_{av}$ and $Sh_{av}$ exhibit decreasing trends when $\lambda=0.2, 0.8, 1.8$.
However, when $\lambda=1.3$, the variations in  $Nu_{av}$ (and $Sh_{av}$) with respect to $A$ are different from those
when $\lambda=0.2, 0.8, 1.8$.  $Nu_{av}$ and $Sh_{av}$ first increase, then decrease monotonically with increasing $A$. A local maximum of $Nu_{av}$ and $Sh_{av}$ is obtained at $A=2$. The observed changing trends in both $Nu_{av}$ and $Sh_{av}$ can be attributed to the flow field transition from a single primary vortex to a coexisting structure comprising one primary vortex and two secondary vortices.

\begin{figure}[h]
	\centering
	\subfigure[$\lambda=0.2$]{
		\includegraphics[width=0.1\textwidth]{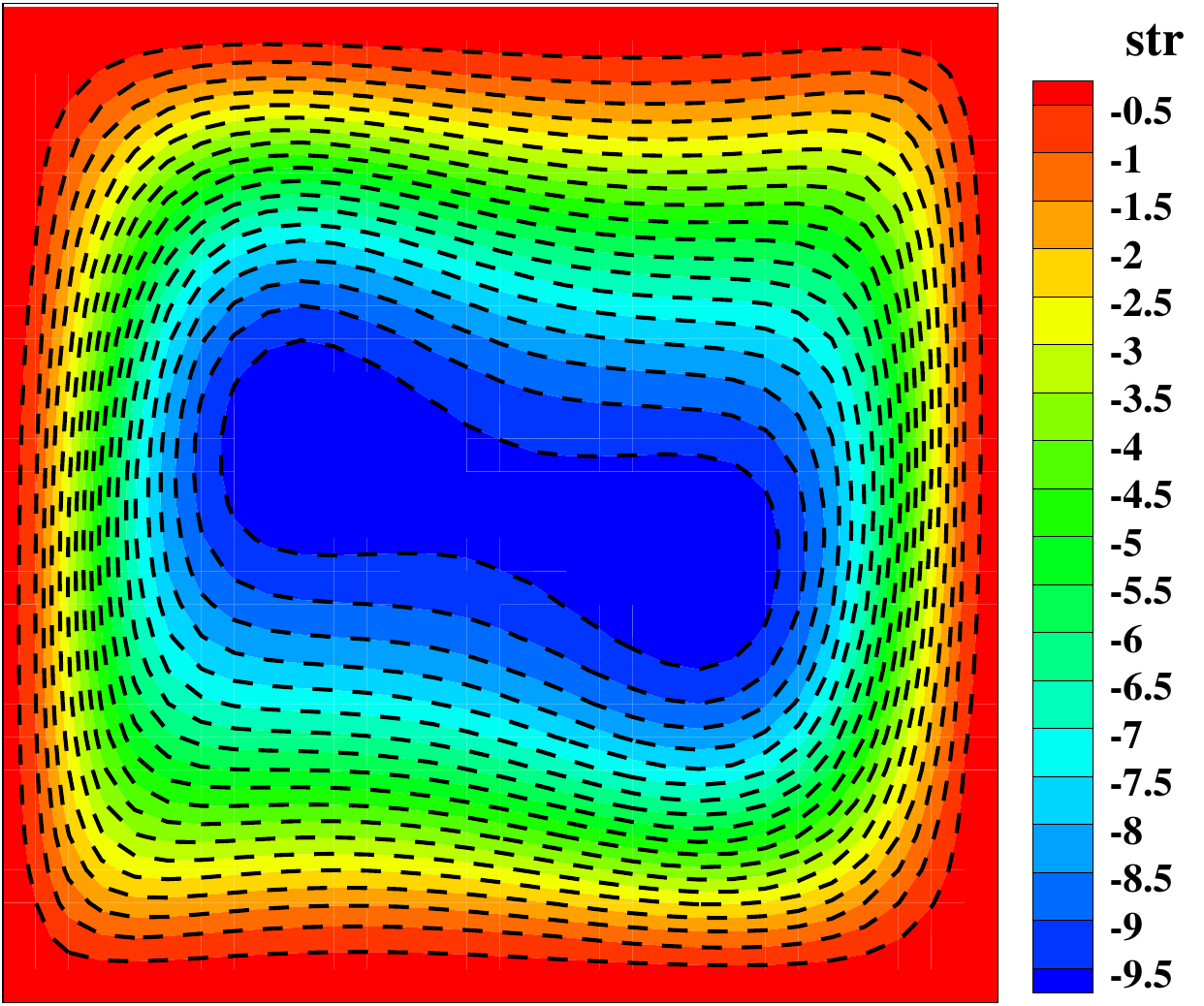}
        \includegraphics[width=0.1\textwidth]{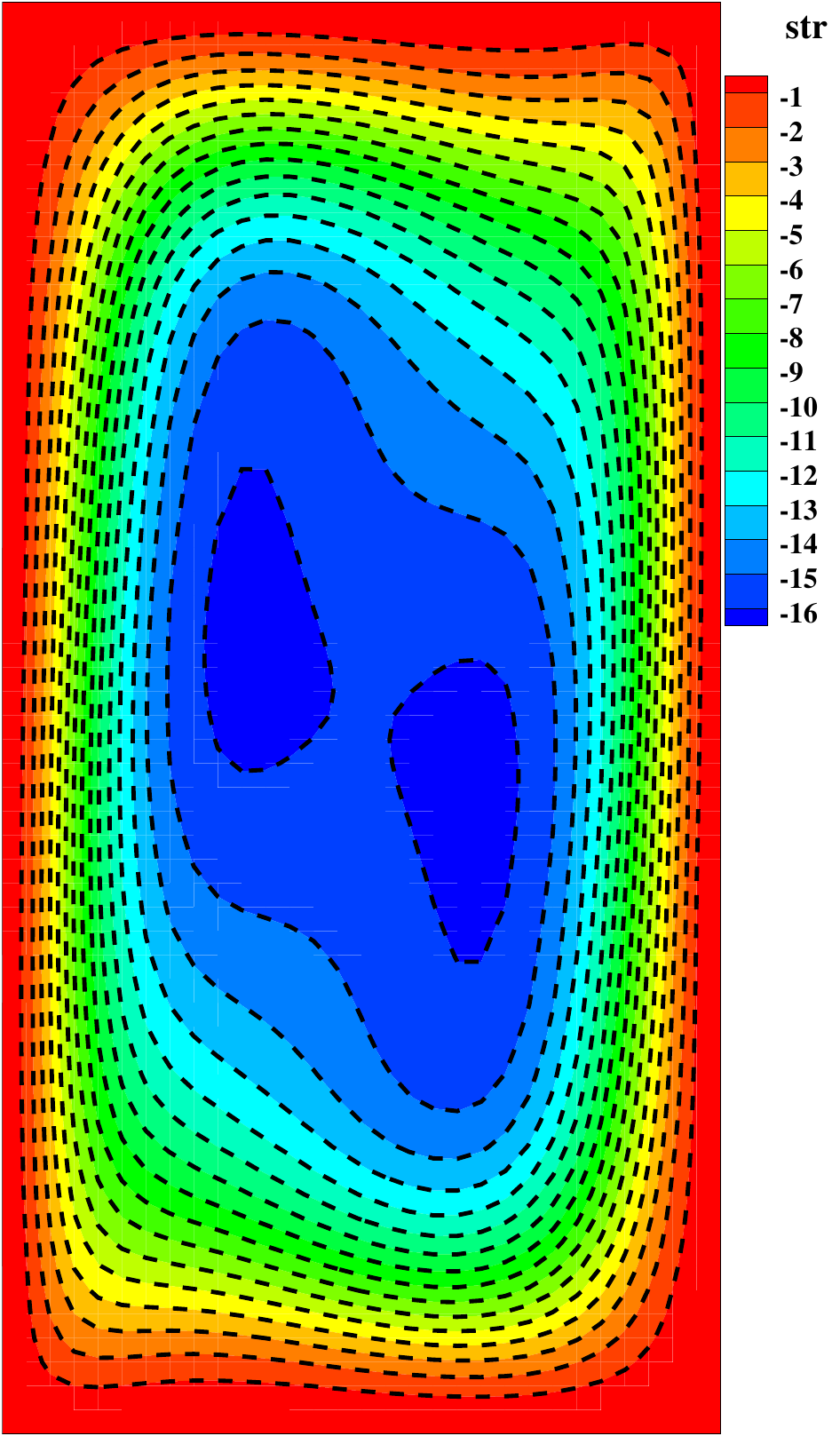}
		\includegraphics[width=0.11\textwidth]{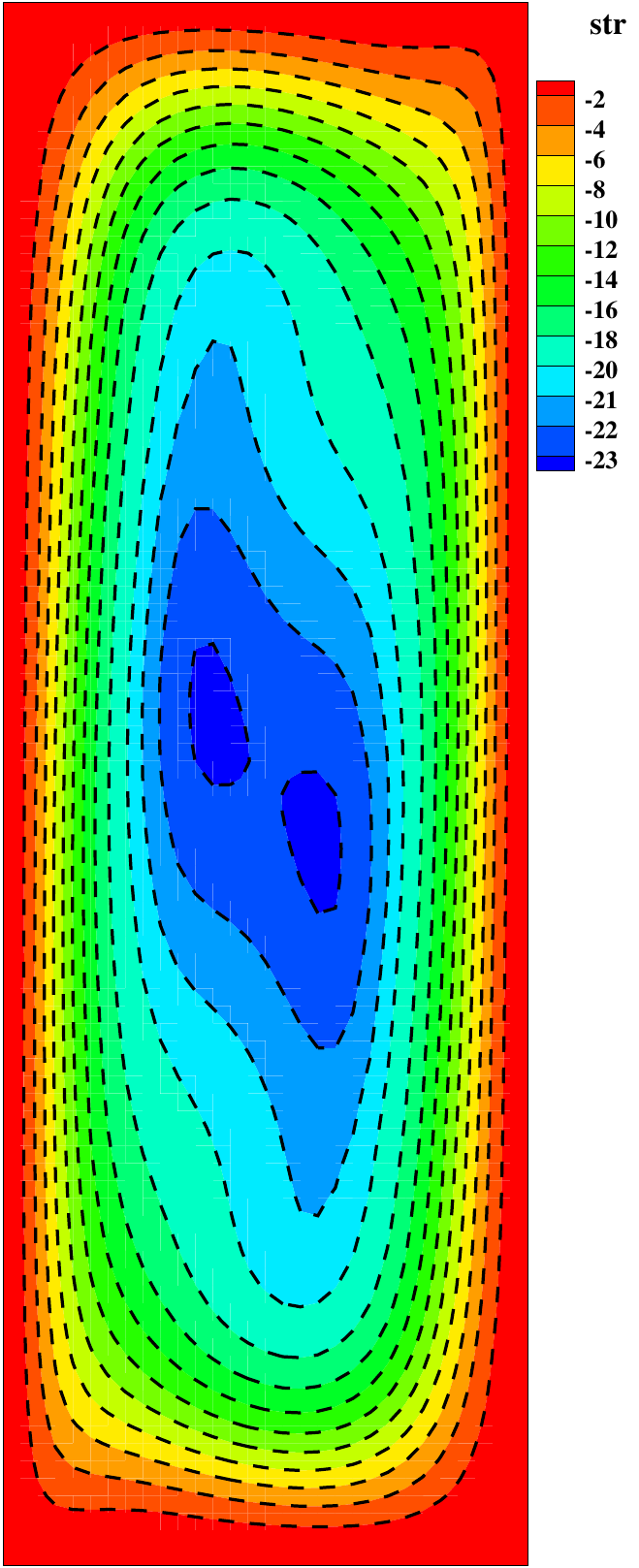}
		\includegraphics[width=0.12\textwidth]{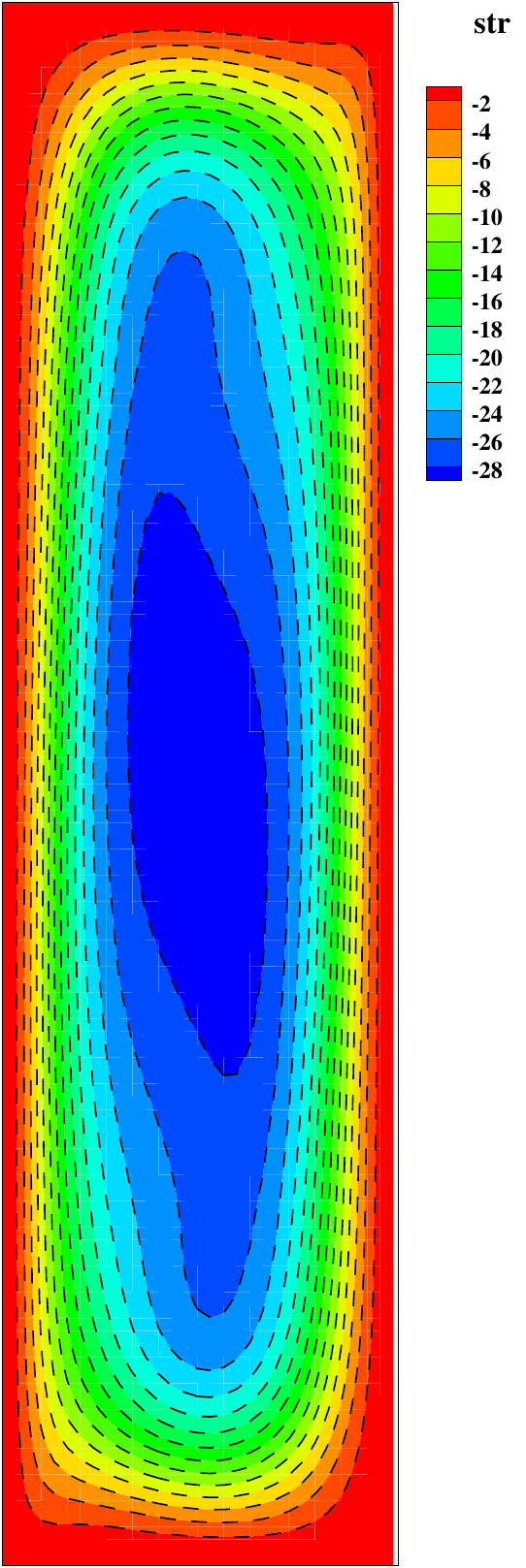}}
    \hspace{0.2cm}
	\subfigure[$\lambda=0.8$]{
		\includegraphics[width=0.1\textwidth]{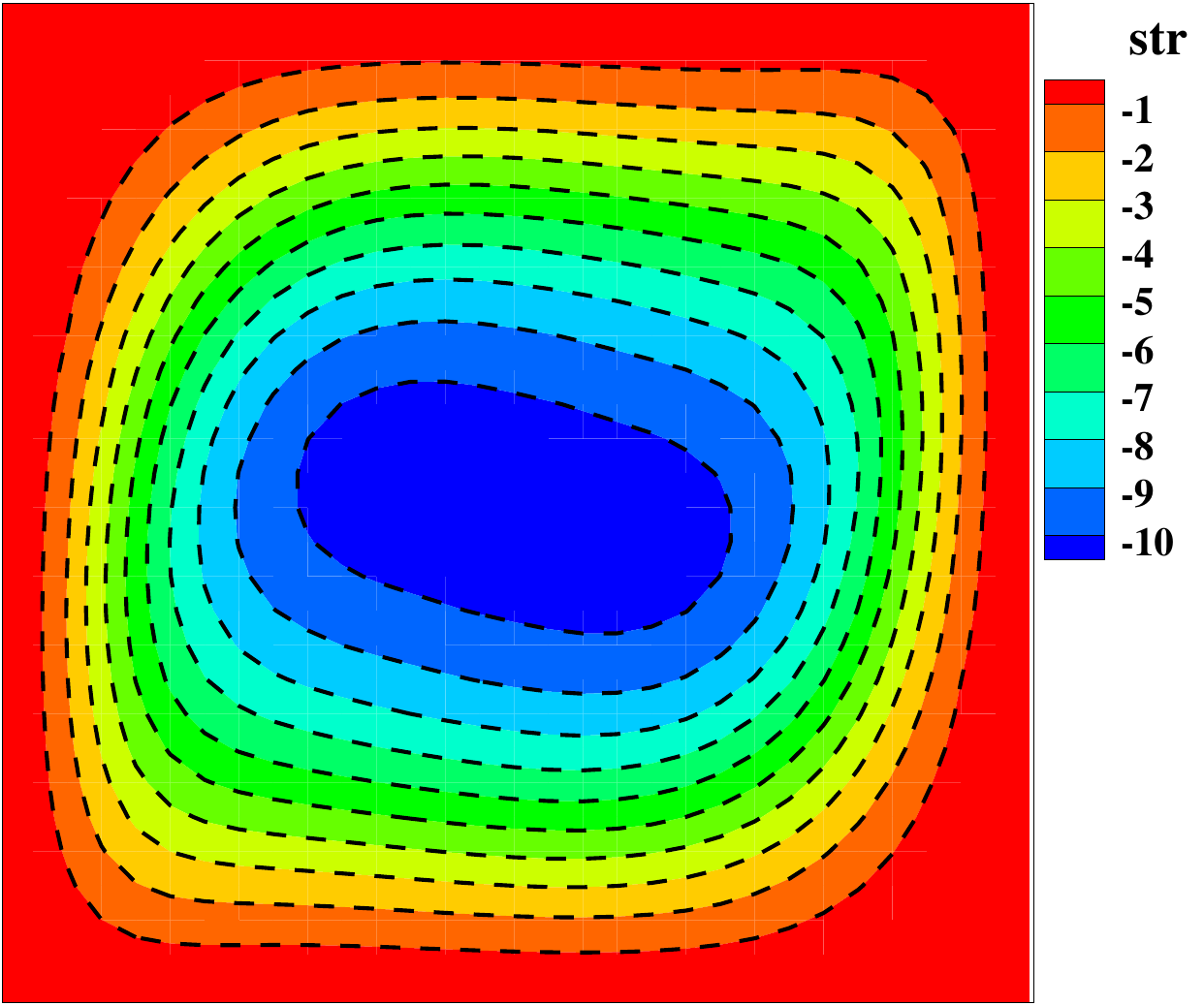}
        \includegraphics[width=0.1\textwidth]{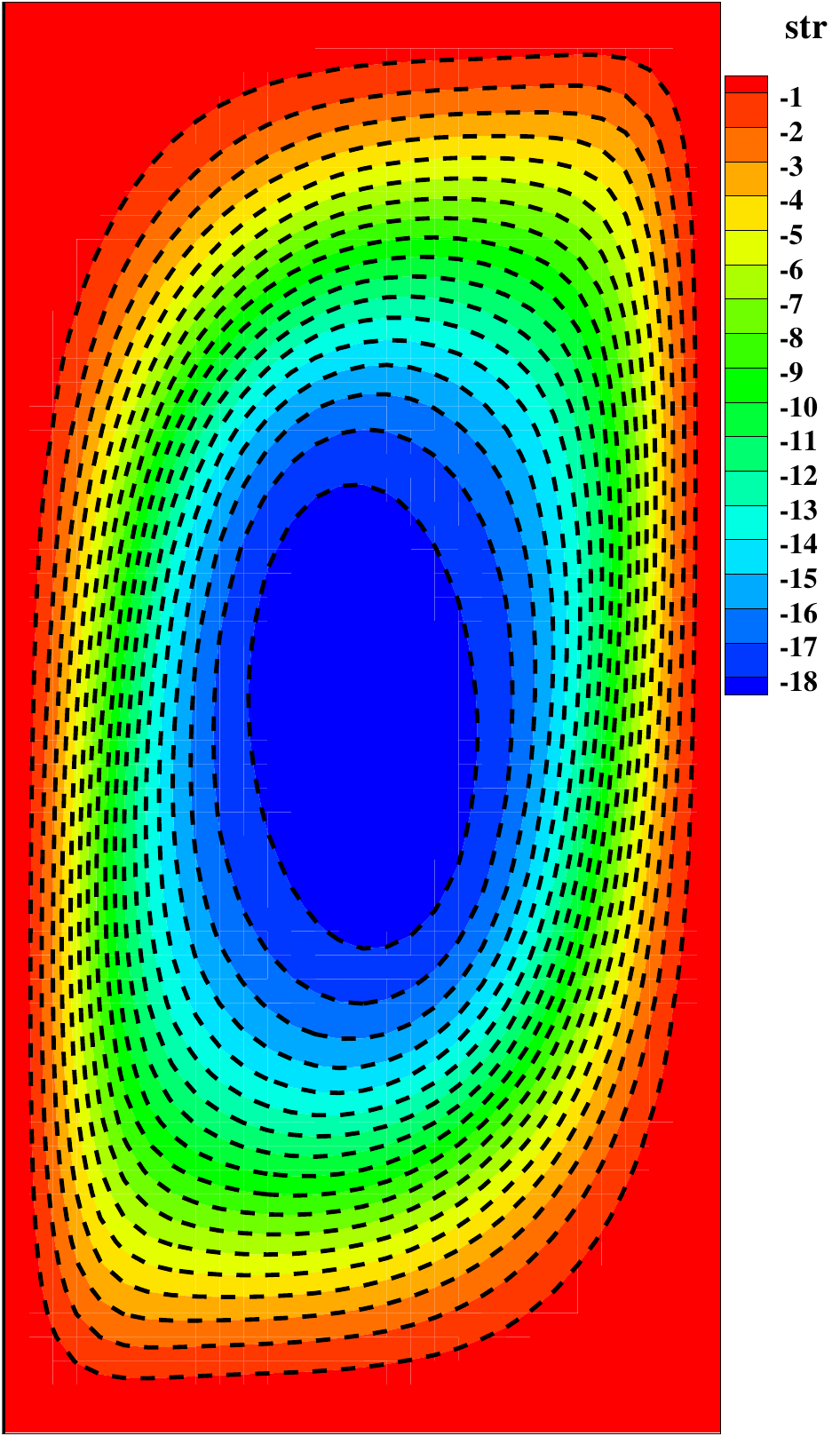}
		\includegraphics[width=0.11\textwidth]{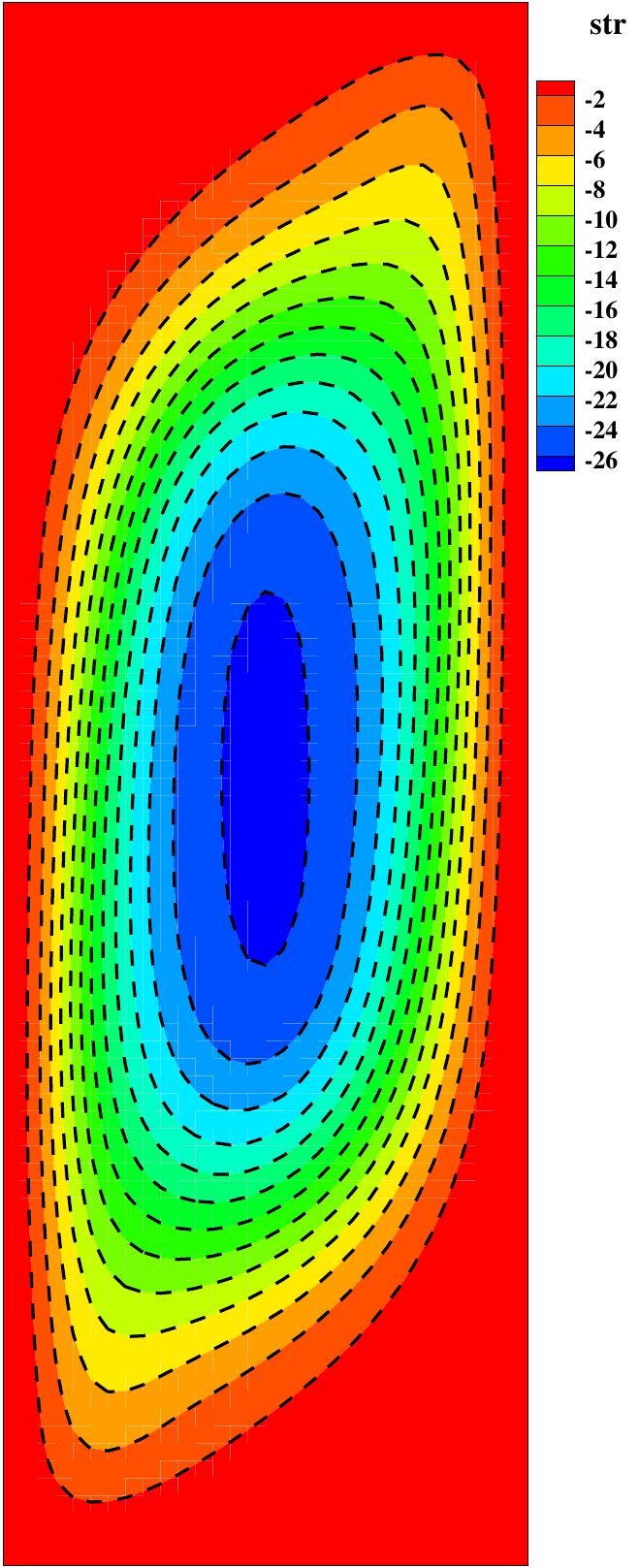}
		\includegraphics[width=0.12\textwidth]{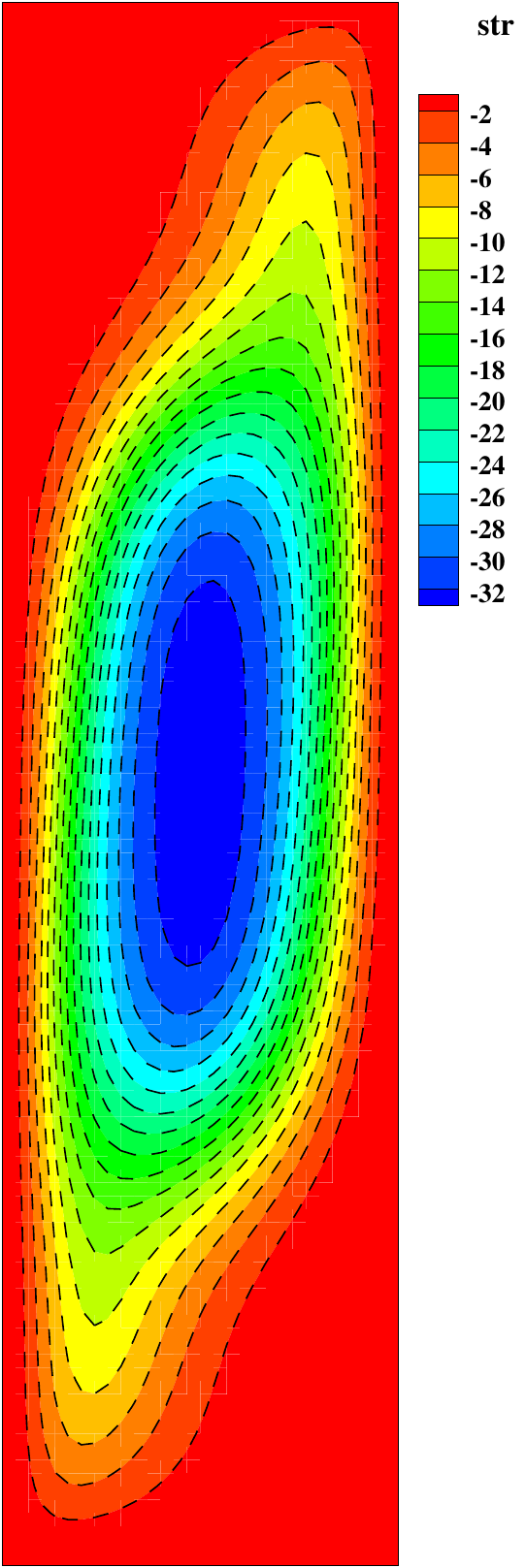}}
    \hspace{0.2cm}
    \subfigure[$\lambda=1.3$]{
		\includegraphics[width=0.1\textwidth]{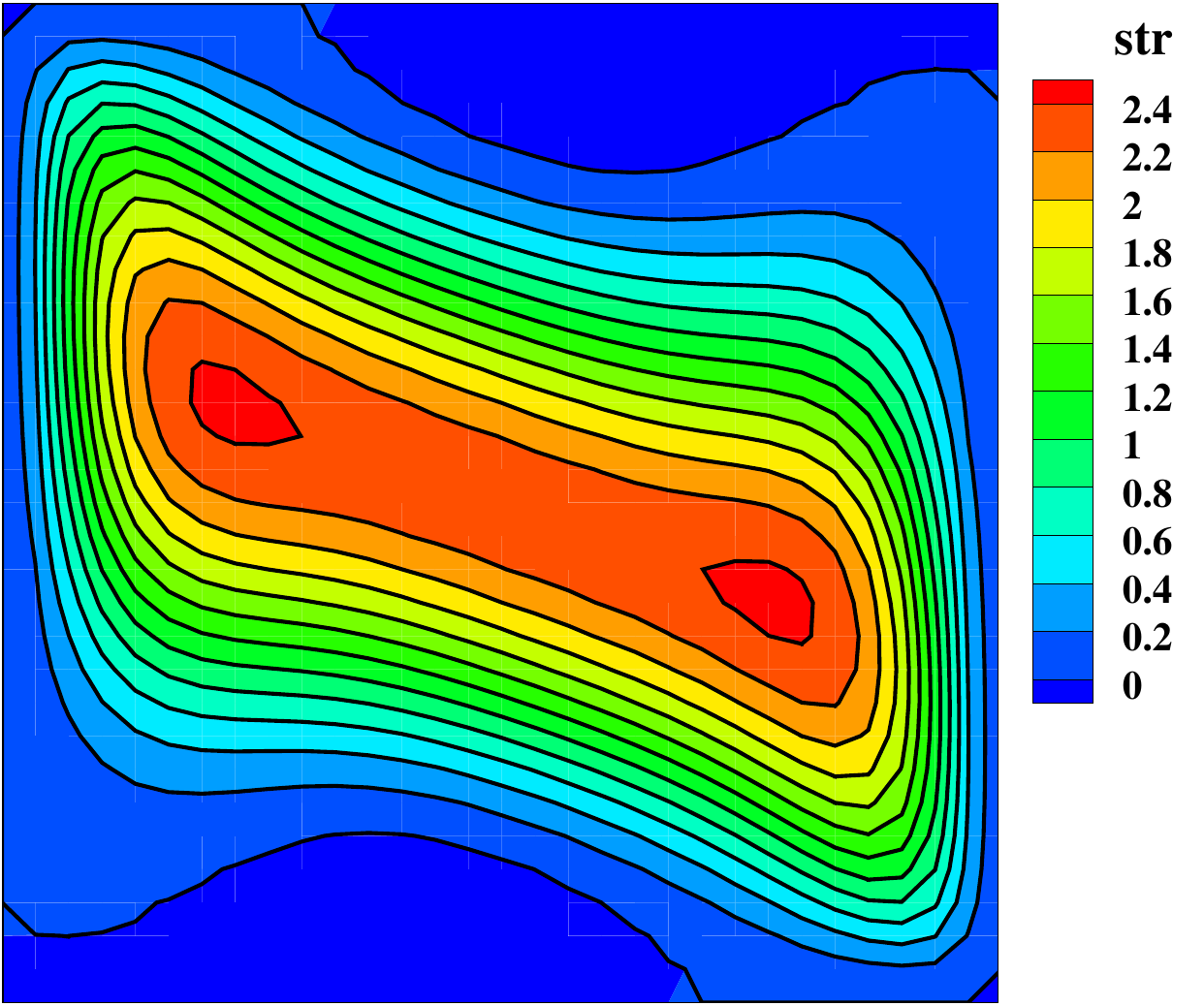}
        \includegraphics[width=0.1\textwidth]{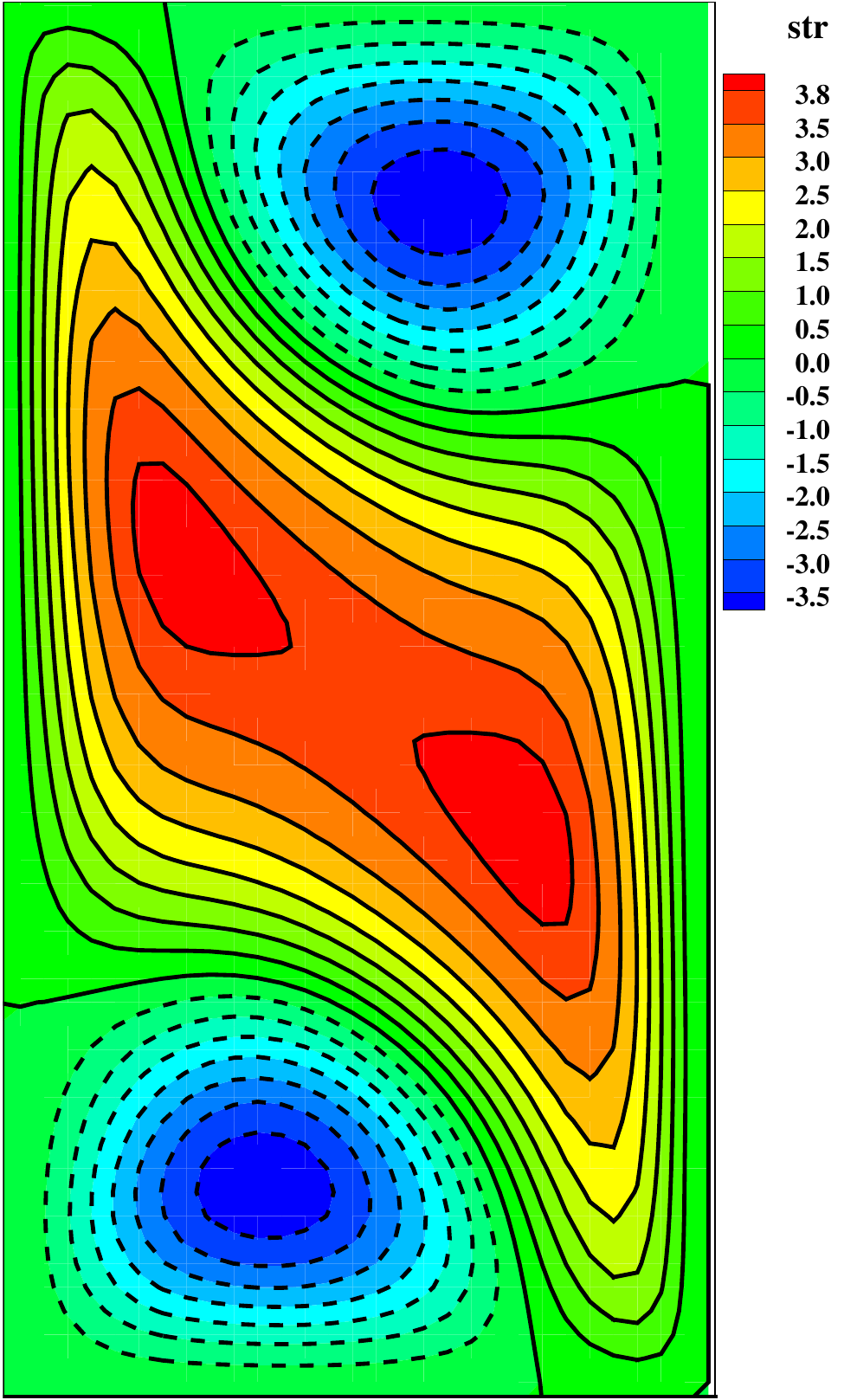}
		\includegraphics[width=0.11\textwidth]{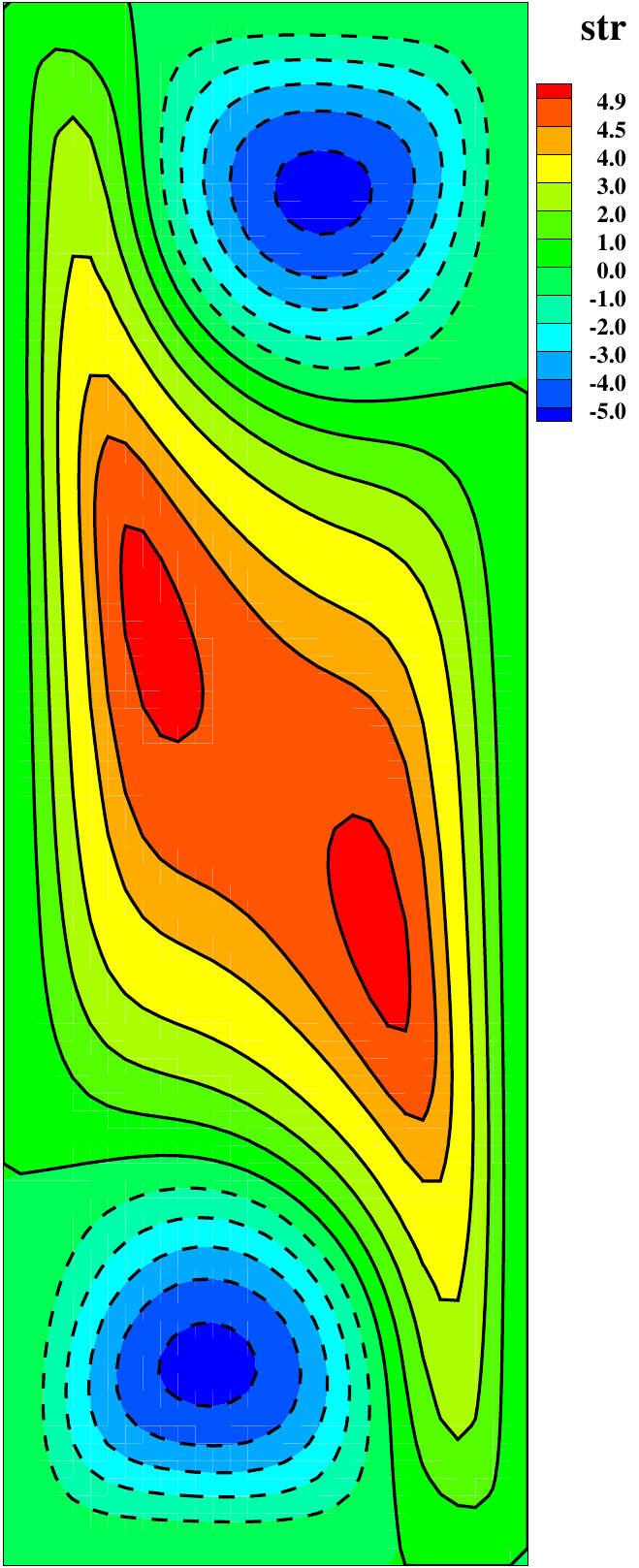}
		\includegraphics[width=0.12\textwidth]{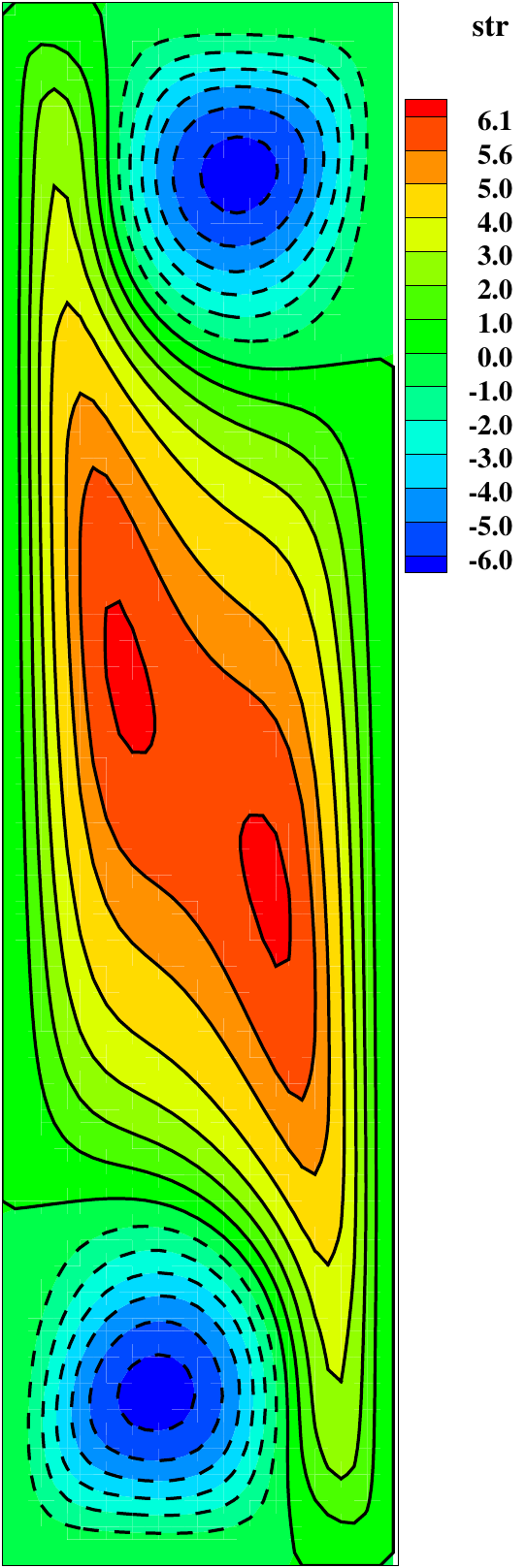}}
    \hspace{0.2cm}
    \subfigure[$\lambda=1.8$]{
		\includegraphics[width=0.1\textwidth]{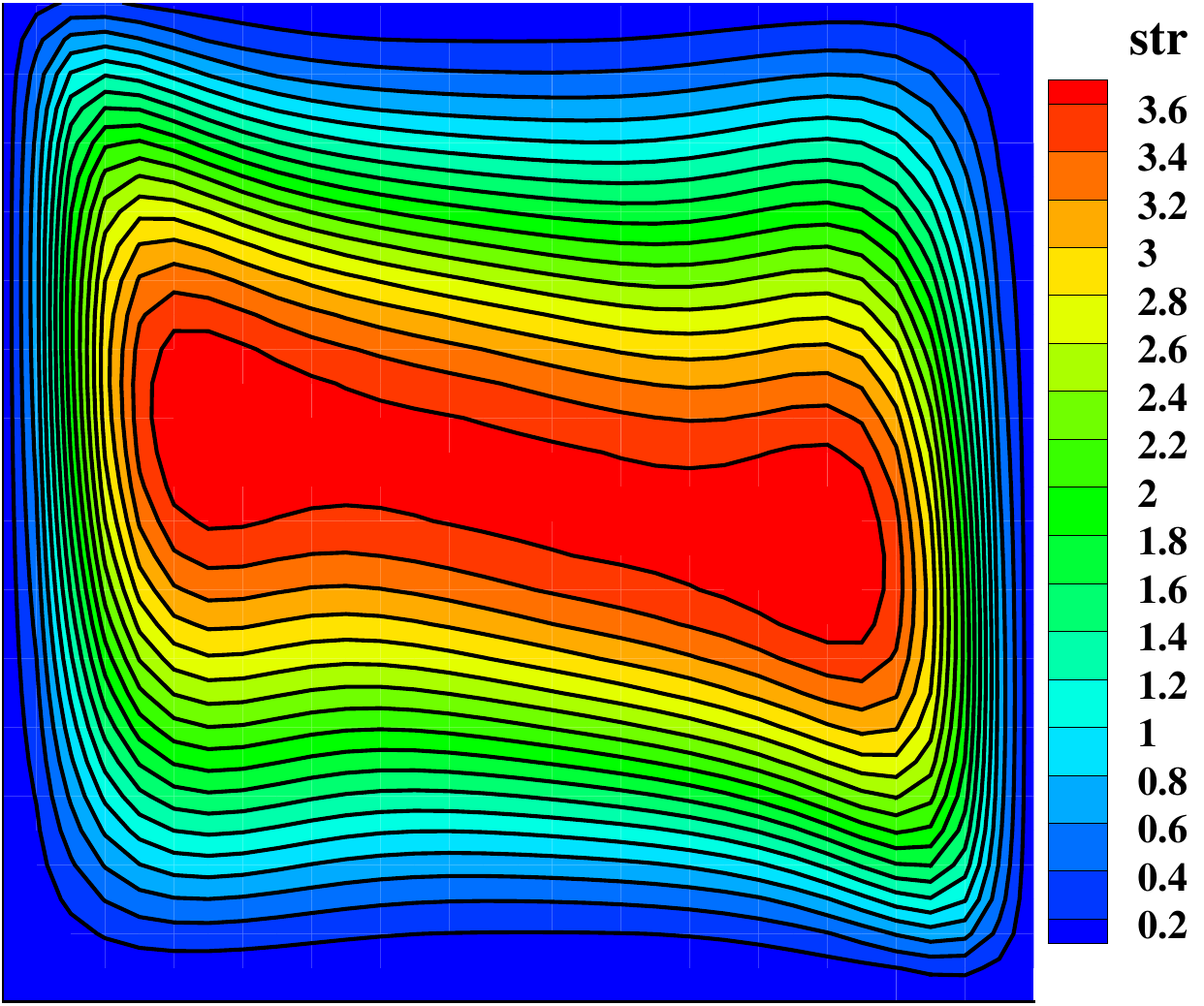}
        \includegraphics[width=0.1\textwidth]{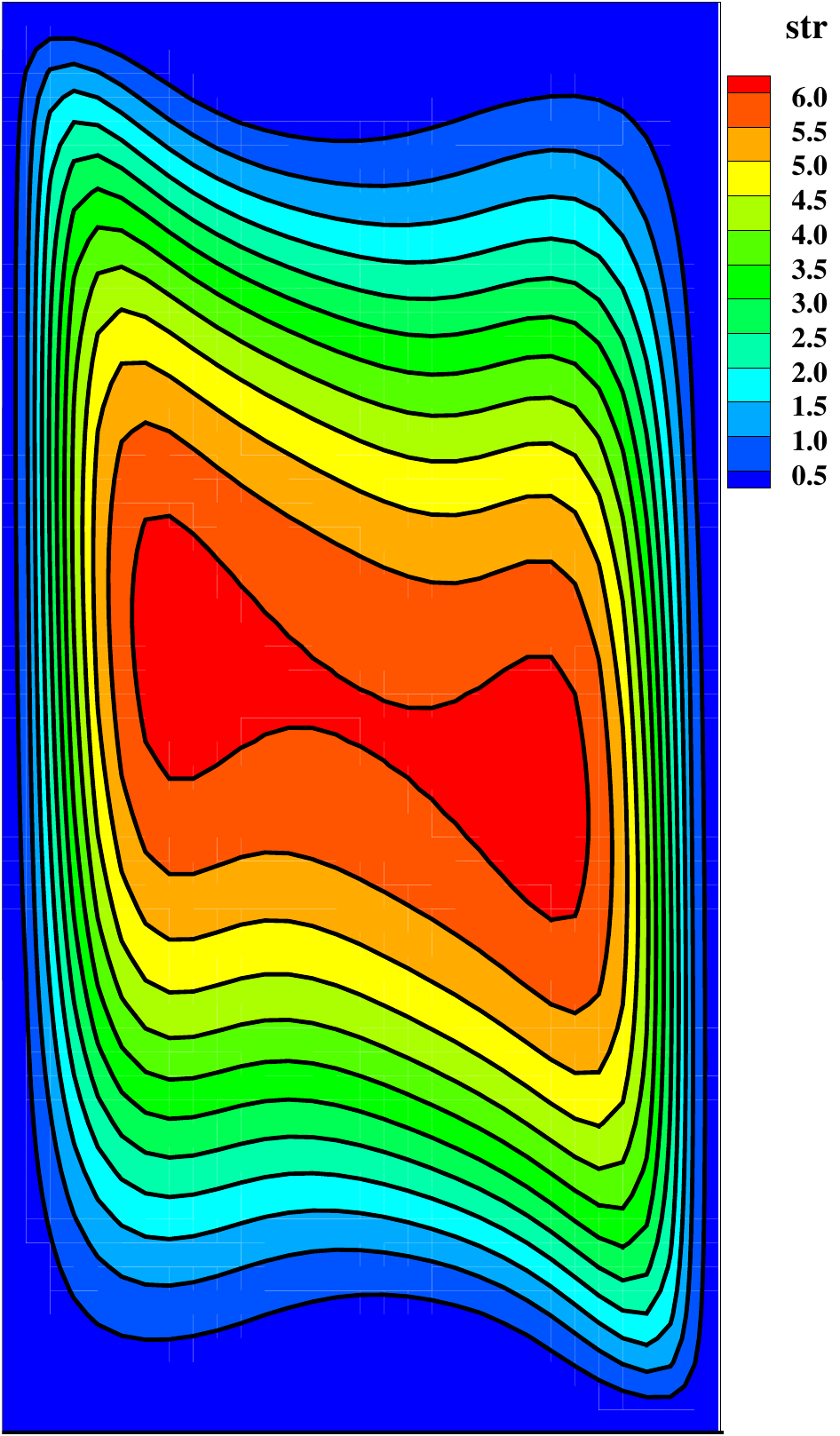}
		\includegraphics[width=0.11\textwidth]{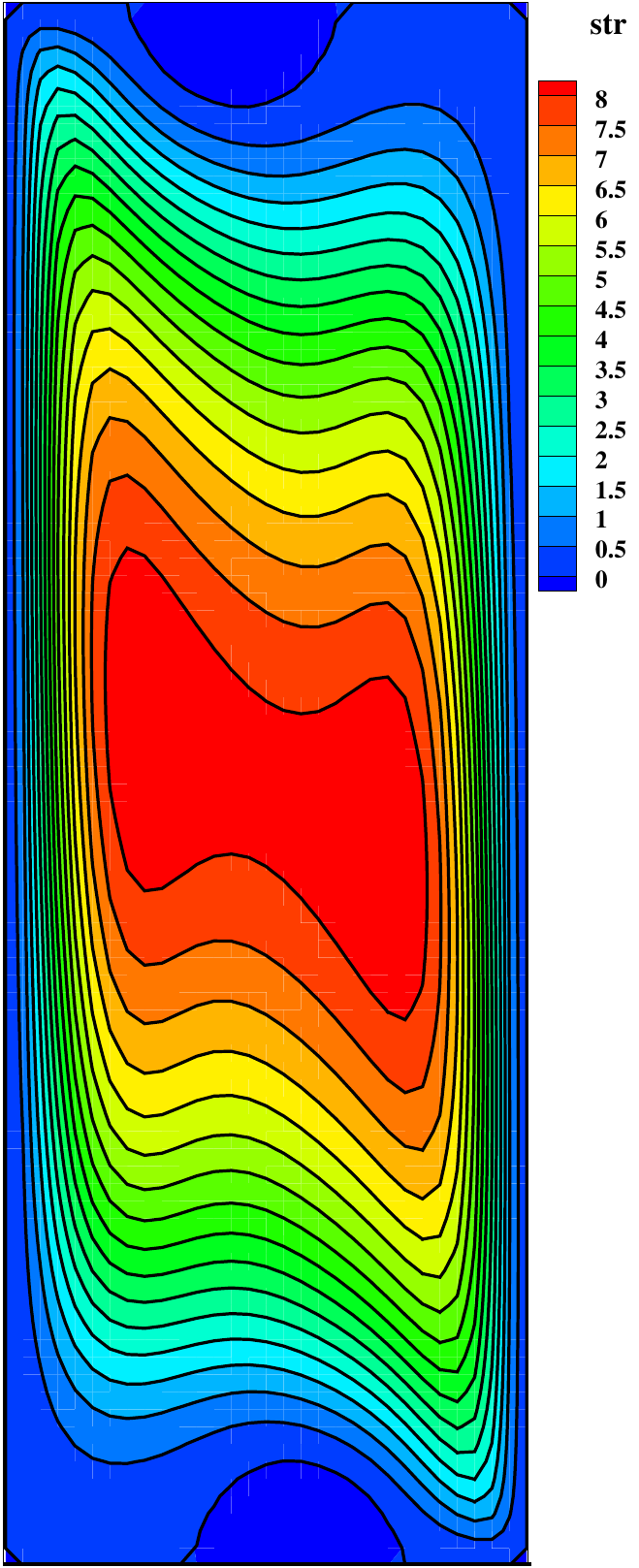}
		\includegraphics[width=0.12\textwidth]{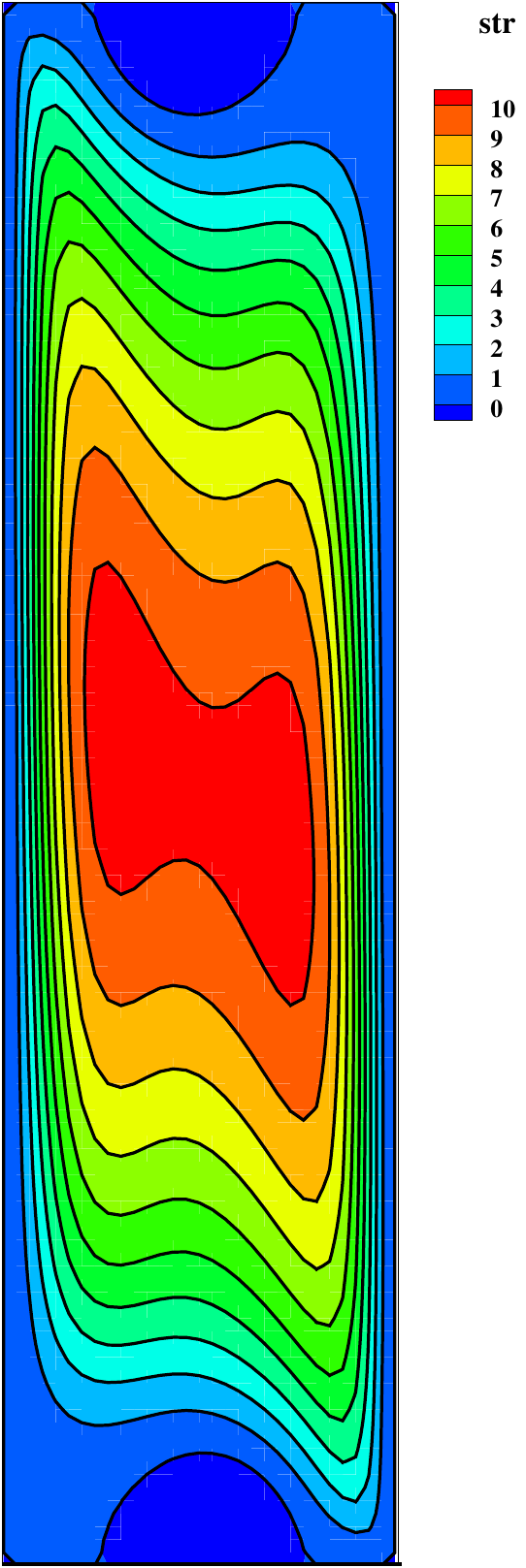}}
	\caption{\label{fig:fq0p2andfq1p3} Solutions at $Ra=10^{5}$, $Pr=1$, and $Le=2$ for different buoyancy ratios: (a) $\lambda=0.2$ and (b) $\lambda=0.8$, $\lambda=1.3$, $\lambda=1.8$. Contour lines of the stream function.}
\end{figure}
The flow patterns of $A=1, 2, 3, 4$ are shown in Fig. \ref{fig:fq0p2andfq1p3} for $\lambda=0.2, 0.8, 1.3, 1.8$.
From the subfigures \ref{fig:fq0p2andfq1p3}(a), (b), (d), there is a primary vortex occupying the main part of the cavity when the aspect ratio changes from $1$ to $4$. The rotation of vortex is clockwise when the thermal buoyancy force is primarily dominated. When the compositional buoyancy force is primarily dominated, the rotation of vortex is counter-clockwise. When the buoyancy ratio is between $5$ and $8$, the flow patterns show no notable change.
When $A=1.3$, the flow patterns exhibit variations with the aspect ratio that differ from the previously discussed cases, as shown
in subfigure \ref{fig:fq0p2andfq1p3} (c).
There is a counter-clockwise rotation primary vortex occupying the main part of the cavity. Then the secondary vortices and primary vortex appear together when $A\geq2$. This variation of flow patterns enhance heat and mass transfer.

\section{\label{sec1:conclusion}Conclusion}
In this work, a high-order compact Hermite scheme with its asymmetric boundary schemes was proposed for solving the governing equations of two dimensional double-diffusive convection. Then the numerical precision of numerical solutions for the proposed method were verified by solving the convection diffusion equation and nonlinear Burgers' equation. The proposed schemes can achieve at least fourth-order or sixth-order accuracy in space. Under identical conditions, our results exhibit higher accuracy compared to existing results.
In the last, the proposed algorithm is used to solving the governing equations of two dimensional double-diffusive convection. The numerical results are in excellent agreement with the benchmark solutions and some of the accurate results available in the literature.
Subsequently, The influence of Rayleigh number, aspect ratio on heat and mass transfer and flow characteristics is investigated.
An increase in the Rayleigh number enhances heat and mass transfer. When Rayleigh number exceeds a critical value, the flow pattern transitions from steady to periodic or chaotic states. The critical value also depends on additional physical parameters, such as buoyancy ratio, etc.
With the increasing of aspect ratio, the heat and mass transfer decrease in the range of aspect ratios considered in this work.
This conclusion is strictly valid for steady flow, while periodic and chaotic regimes may exhibit different behaviors.


\end{document}